\documentclass[11pt,reqno]{amsart}

\title[]{Long time dynamics of a fractional dissipative model of electroconvection in bounded domains}

\author{Elie Abdo}
\address{Department of Mathematics, University of California Santa Barbara, CA 93106, USA.}
\email{elieabdo@ucsb.edu}

\author{Mihaela Ignatova}
\address{Department of Mathematics, Temple University, Philadelphia, PA 19122}
\email{ignatova@temple.edu}

\usepackage[margin=1in]{geometry}
\usepackage{amsmath, amsthm, amssymb}
\usepackage{times}
\usepackage{color}
\usepackage{hyperref}
\newcommand{\pa}{\partial}
\newcommand{\la}{\label}
\newcommand{\fr}{\frac}
\newcommand{\na}{\nabla}
\newcommand{\be}{\begin{equation}}
\newcommand{\ee}{\end{equation}}
\newcommand{\ba}{\begin{array}{l}}
\newcommand{\ea}{\end{array}}

\newtheorem{prop}{Proposition}

\newcommand{\beg}{\begin}

\renewcommand{\l}{\Lambda}

\newtheorem{Thm}{Theorem}

\usepackage{MnSymbol,wasysym}

\newcommand{\N}{\mathbb N}
\newcommand{\C}{\mathbb C}
\newcommand{\R}{\mathbb R}

\def\ZZ{{\mathbb Z}}
\def\RR{{\mathbb R}}
\def\TT{{\mathbb T}}
\def\NN{\mathbb N}
\def\PP{\mathbb P}

\date{\today}
\begin{document}
\begin{abstract} We consider a nonlocal nonlinear model with fractional diffusion motivated by studies of electroconvection phenomena in incompressible viscous fluids. We address the global well-posedness, global regularity and long time dynamics of the model in bounded smooth domains  with Dirichlet boundary conditions. We prove the existence and uniqueness of exponentially decaying in time solutions for $H^1$ initial data regardless of the fractional dissipative regularity. In the presence of time independent body forces in the fluid, we prove the existence of a compact finite dimensional global attractor. In the case of periodic boundary conditions, we prove that the unique smooth solution is globally analytic in time, and belongs to a Gevrey class of functions that depends on the dissipative regularity of the model. 
 \end{abstract} 

\maketitle
\section{Introduction}\la{intro}

Electroconvection is a term associated to the nonlinear dynamics created by the interaction of  fluid flow, ionic transport and electrostatic forces. In certain controlled experimental situations the dynamics are chaotic, similar to classical hydrodynamic transition to turbulence. The analogy to Rayleigh-B\'{e}nard convection \cite{DDMB, mani}, is motivated not only by qualitative observations, but also by the fact that in both systems the fluid is driven by body forces which are the product of a transported scalar and a vector field. In thermal convection the scalar is the temperature and the vector is the gravitational field, while in electroconvection, the scalar is the charge density and the vector is the electric field. Numerical simulations and physical experiments have been used to study electroconvection in thin liquid crystals \cite{BW,DDMB,DDMT}. Electroconvection is of broad interest in electrochemistry, material science and applied physics (see for instance \cite{ZYZZHLYHLBZZ,MLCAW,TR}), but our motivation and focus is on mathematical challenges of long time behavior, in the  important case when physical boundaries are present. 

In \cite{ceiv}, the authors considered an electroconvection model describing the nonlinear time evolution of a surface charge density $q$ in an incompressible viscous fluid, confined to a two dimensional bounded domain $\Omega$, with  velocity $u$ and pressure $p$. The model is described by the system 
\be \la{inmod1}
\pa_t q + u \cdot \na q + \l q = \Delta \Phi,
\ee
\be 
\pa_t u + u \cdot \na u - \Delta u + \na p = -qRq - q\na \Phi + f,
\ee
\be 
\na \cdot u = 0,
\ee 
\be 
q|_{\pa \Omega}  = u|_{\pa \Omega} = 0,
\ee 
\be \la{inmod2}
u(x,0) = u_0(x), q(x,0) = q_0(x),
\ee where $\l$ is the square root of the Dirichlet Laplacian, $R:= \na \l^{-1}$ is the Riesz transform,  $f$ is a time independent body force in the fluid, and $\Phi$ is a time independent potential resulting from a boundary applied  voltage. 

The global regularity of a unique solution to the initial boundary value problem \eqref{inmod1}--\eqref{inmod2} was obtained in \cite{ceiv} for Sobolev $H^2$ initial data based on a two-tier Galerkin approximation. On the torus $\TT^2$ with periodic boundary conditions, we showed in \cite{AI} that the system \eqref{inmod1}--\eqref{inmod2} has a unique strong solution provided that the initial charge density belongs to the Lebesgue space $L^4$  and the initial velocity belongs to the Sobolev space $H^1$. Moreover, we obtained the existence of a finite dimensional global attractor which reduces to a single point in the absence of body forces in the fluid. 

In this paper, we are interested in the long time dynamics of \eqref{inmod1}--\eqref{inmod2} on bounded domains with smooth boundaries. We fix $\alpha \in (0,1]$ and consider the generalized electroconvection model in $\Omega$ described by the system 
\be \la{inmod11}
\pa_t q + u \cdot \na q + \l^{\alpha} q = \Delta \Phi,
\ee
\be 
\pa_t u + u \cdot \na u - \Delta u + \na p = -qRq - q\na \Phi + f,
\ee
\be \la{inmod111}
\na \cdot u = 0,
\ee 
\be 
q|_{\pa \Omega}  = u|_{\pa \Omega} = 0,
\ee 
\be \la{inmod21}
u(x,0) = u_0(x), q(x,0) = q_0(x),
\ee
with a fractional diffusion driven by the operator $\l^{\alpha}$. We address the following three main problems:
\begin{enumerate}
\item[(i)] Global existence and uniqueness of solutions for Sobolev $H^1$ initial data;
\item[(ii)] Long time dynamics in the absence and in the presence of body forces and of voltage applied at the boundary;
\item[(ii)] Global Gevrey regularity of solutions in the case of periodic boundary conditions. 
\end{enumerate}

The system \eqref{inmod11}--\eqref{inmod21} is reminiscent of the dissipative surface quasi-geostrophic (SQG) equation, proposed in \cite{CMT} as a model of hydrodynamic creation of small scales. In SQG the fluid velocity $u$ depends on the scalar $q$ via the relation $u = \na^{\perp} (-\Delta)^{-\fr{1}{2}}q$. In \cite{CW1}, the existence and uniqueness of global smooth solutions were obtained in the subcritical $(\alpha > 1)$ case, whereas the existence of global decaying weak solutions was obtained in the supercritical $(\alpha < 1)$ and critical $(\alpha= 1)$ cases. Global regularity of solutions to the critical SQG equation on $\R^2$ was established in \cite{knv} based on modulus of continuity techniques, in \cite{caf} based on De Giorgi techniques, and in \cite{CV} based on nonlinear maximum principles.  In \cite{cvt}, the authors addressed the long time dynamics of the forced critical SQG in the spatially periodic case and proved the existence of a finite dimensional global attractor. The global well-posedness of the supercritical SQG equation on $\RR^2$ was obtained in \cite{HK} for small initial data in Besov spaces, and the supercritical regularity was studied in \cite{CW} where the authors proved that H\"older continuous solutions of subcritical type are actually $C^{\infty}$ classical solutions of the equation. In \cite{CZV}, it was shown that the solution of the supercritical equation with periodic boundary conditions does not blow up in finite time for all fractional powers $\gamma \ge \gamma_1$ where $\gamma_1$ is a constant depending on the size of the initial data. Recently, the critical SQG equation in bounded smooth domains was addressed in \cite{CI}, \cite{CI3}, and \cite{I} where global interior Lipschitz solutions were constructed, and in \cite{SV} where global H\"older regularity up to the boundary was obtained. In \cite{CN3}, the local well-posedness for the inviscid case and the global existence of strong solutions for small initial data in the supercritical and critical cases are established in bounded smooth domains. The global regularity and time asymptotic behavior of solutions to the critical and supercritical SQG equations on bounded domains are open problems. 

The system \eqref{inmod11}--\eqref{inmod21} has a different smoothness balance than the supercritical SQG equation due to the coupling to the  Navier-Stokes equations which results in a higher spatial regularity for the fluid velocity. However, many challenges arise, not only from the nonlocality  and the nonlinearity of the electric forces driving the fluid velocity, but also and chiefly from the presence of boundaries. 

The existence and uniqueness of solutions of \eqref{inmod11}--\eqref{inmod21} relies on control of the spatial $L^p$ norms of the charge density $q$, which evolve via regular nonlinear advection  $u \cdot \na q$. The need for cancellation of advective terms  in $L^p$ is crucial, and a direct Galerkin approximation procedure does not work. We consider instead a spectral regularization of $(-\Delta)^{-\fr{1}{2}}$, denoted by $(\l^{-1})_{\epsilon}$, that depends on a small positive parameter $\epsilon > 0$, and we define the corresponding truncated Riesz transform $R_{\epsilon} = \na (\l^{-1})_{\epsilon}$. Then we take a regularized version of \eqref{inmod11}--\eqref{inmod21} in which the nonlinear nonlocal electric forces $qRq$ are replaced by $q^{\epsilon}R_{\epsilon}q^{\epsilon}$, and we use Galerkin approximations and compactness arguments to prove that each $\epsilon$-approximate system has global in time regular up to the boundary solutions that may depend badly on $\epsilon$. By making use of convex damping inequalities (Proposition \ref{corcor}), we manage to derive $L^p$ bounds for the family of viscous charge densities $\left\{q^{\epsilon}\right\}_{\epsilon > 0}$, uniform in $\epsilon$. This allows us to obtain good control of $q^{\epsilon} R_{\epsilon} q^{\epsilon}$, uniform in $\epsilon$, due to the boundedness of the Dirichlet Riesz transform on $L^p$ spaces, generating consequently a spatial Sobolev $H^1$ regularity, global in time, for both the charge density and velocity solving \eqref{inmod11}--\eqref{inmod21}. 

The long time dynamics of the forced electroconvection system \eqref{inmod1}--\eqref{inmod2} in periodic domains was addressed in \cite{AI} on the basis of Fourier series techniques employed in the study of commutators 
$$[\l^{s}, u \cdot \na ] q:=  \l^{s}(u \cdot \na q) -  u \cdot  \na \l^{s} q$$ for positive and negative powers $s$, and of interpolation inequalities for fractional powers of the Laplacian. The facts that the fractional Laplacians have explicit representation formulas and 
that the periodic operators $\l^{s}$, defined as Fourier multipliers, commute with differential operators were essentially used in that work. These properties and techniques break down on bounded domains where the nonlocal operators $\l^{s}$ are defined via eigenfunction expansions, in terms of the eigenfunctions of the homogeneous Dirichlet Laplace operator, are not translation invariant, and don't have integral representations with explicit kernels. This gives rise to many technical mathematical challenges and the need for new ideas.  

At each positive time $t$, the forced initial boundary value problem \eqref{inmod11}--\eqref{inmod21} has a well defined solution map $\mathcal{S}(t)$ on 
$$\mathcal{V} = \left\{(q,u) \in H_0^1 (\Omega) \times (H_0^1(\Omega))^2 : \na \cdot u = 0 \right\} ,$$
which is the largest space in which the model \eqref{inmod11}--\eqref{inmod21} has unique solutions. We prove the existence of an absorbing ball
$$\mathcal{B} = \left\{(q, u) \in \mathcal{V} : \|\na q\|_{L^2} + \|\Delta u\|_{L^2} \le R \right\}$$
with a radius $R$ depending only on the forcing terms $f$ and $\Phi$, whose image under $\mathcal{S}(t)$ is a subset of $\mathcal{B}$ itself starting at a time $T:= T(R)$ depending only on $R$ (Proposition \ref{absball}). This requires uniform control of $q$ in $H^1$ and $u$ in $H^2$ starting with Sobolev $H^1$ initial data. The first main challenging elements of the proof are  $L_{t}^{\infty}$ and $L_t^2$ uniform boundedness of the velocity in $H^{1+\epsilon}$  and $H^{2+\epsilon}$ respectively. Due to the presence of electric forces $qRq$ driven by the charge density, fractional product inequalities for small powers of the Stokes operator $A$ are needed to estimate $A^{\epsilon}(\PP(qRq))$ and are established in Proposition \ref{prodfor} below. The second major element of the proof is $L_{t}^{\infty}$ boundedness of the $L^p$ norm of the velocity gradients. These are obtained via a mild formulation of the Navier-Stokes equations and use of Stokes semi-group estimates. The third challenging element of the proof is $L_{t}^{\infty}$ uniform control of the density gradient in $L^2$. For that purpose, we track the time evolution of $\|\na q\|_{L^2}$ via energy estimates and handle the nonlinearity by making use of velocity gradient bounds in $L^p$ and fractional interpolation inequalities that hold on bounded domains. Consequently, the desired property $\mathcal{S}(t) \mathcal{B} \subset \mathcal{B}$ follows after waiting some time. The compactness of the ball $\mathcal{B}$ in the weaker topology of 
$$ \mathcal{H} = \left\{(q,u) \in L^2 (\Omega) \times (L^2(\Omega))^2 : \na \cdot u = 0, q|_{\pa \Omega} = u|_{\pa \Omega}=0 \right\},$$
together with the continuity and injectivity properties of $\mathcal{S}(t)$, yield the existence of a global attractor, compact in the norm of $\mathcal{H}$. 

The passage from the $H^1$ to the $H^{\fr{\alpha}{2}}$ regularity of the charge density is intriguing because the nonlinear transport and dissipation governing the evolution of the $H^{\fr{\alpha}{2}}$ norm of $q$ do not balance each other. Nonetheless, the case $\alpha = 1$ provides a stronger dissipative structure that is exploited to improve the smoothness of the global attractor. In this latter situation, we control the nonlinearity of \eqref{inmod1} by establishing pointwise commutator estimates that are inversely proportional to powers of the distance to the boundary function 
$$d(x) := d(x, \pa \Omega),$$ and we control the weighted vector field $u(x)/d(x)$ and scalar function $q(x)/d(x)$ in $L^p$ via use of Rellich inequalities. This yields good control of the fractional energy norms $\|\l^{\fr{s}{2}} q \|_{L^2}$ for any integer $s \ge 1$, and upgrade of the attractor's regularity via bootstrapping arguments. 

Posed on the two dimensional torus $\TT^2$ with periodic boundary conditions, the system \eqref{inmod11}--\eqref{inmod21} has a unique Gevrey regular solution. In spite of the fractional diffusion governing the system, we show that the time evolution of the Gevrey norm depends on the dissipative structure at hand using Fourier series techniques, Gevrey commutator estimates, and Gevrey cancellation laws. We obtain a local in time control of the Gevrey norm by the Sobolev regularity of the solution and, consequently a global extension in the spirit of the global in time boundedness of solutions in fractional Sobolev spaces.

\section{Main Results}

\subsection{Functional Setting.} Let $\Omega \subset \R^2$ be a bounded domain with smooth boundary. For $1 \le p \le \infty$, we denote by $L^p(\Omega)$ the Lebesgue spaces of measurable functions $f$ from $\Omega$ to $\R$ (or $\RR^2)$ such that 
\be 
\|f\|_{L^p} = \left(\int_{\Omega} \|f\|^p dx\right)^{1/p} <\infty
\ee if $p \in [1, \infty)$ and 
\be 
\|f\|_{L^{\infty}} = {\mathrm{esssup}}_{\Omega}  |f| < \infty
\ee if $p = \infty$. The $L^2$ inner product is denoted by $(\cdot,\cdot)_{L^2}$. 

For $k \in \NN$, we denote by $H^k(\Omega)$ the Sobolev spaces of measurable functions $f$ from $\Omega$ to $\R$ (or $\RR^2)$ with weak derivatives of order $k$ such that  
\be 
\|f\|_{H^k}^2 = \sum\limits_{|\alpha| \le k} \|D^{\alpha}f\|_{L^2}^2 < \infty,
\ee and by $H_{0}^{1}(\Omega)$ the closure of $C_0^{\infty}(\Omega)$ in $H^1(\Omega)$. 

For a Banach space $(X, \|\cdot\|_{X})$ and $p\in [1,\infty]$, we consider the Lebesgue spaces $ L^p(0,T; X)$ of functions $f$  from $X$ to $\R$ (or $\RR^2)$ satisfying 
\be 
\int_{0}^{T} \|f\|_{X}^p dt  <\infty
\ee with the usual convention when $p = \infty$. 

{\bf{Fractional Powers of the Laplacian}}. We denote by $\Delta$ the Laplacian operator with homogeneous Dirichlet boundary conditions. We note that $-\Delta$  is defined on $\mathcal{D}(-\Delta) = H^2(\Omega) \cap H_0^1(\Omega)$, and is positive and self-adjoint in $L^2(\Omega)$. We consider an orthonormal basis of $L^2(\Omega)$ consisting of eigenfunctions $\left\{w_j\right\}_{j=1}^{\infty} \subset H_0^1(\Omega)$ of $-\Delta$ satisfying
\be
-\Delta w_j = \lambda_j w_j
\ee
where the eigenvalues $\lambda_j$ obey $0 < \lambda_1 \leq ... \leq \lambda_j \le ...\rightarrow \infty$. 
For $s \in \R$, we define the fractional Laplacian operator of order $s$, denoted by $\l^s$, by
\be \la{maindef}
\l^s f= \sum_{j=1}^{\infty} \lambda_j^{\fr{s}{2}} (f, w_j)_{L^2} w_j
\ee with domain 
\be 
\mathcal{D}(\l^s) = \left\{f \in L^2(\Omega) : \|\l^{s} f\|_{L^2}^2 := \sum\limits_{j \in \N} \lambda_j^{s}(f, w_j)_{L^2}^2 < \infty \right\}.
\ee 
For $s \in [0,1]$, we identify the domains $\mathcal{D}(\l^s)$ with the usual Sobolev spaces as follows, 
\be 
\mathcal{D}(\l^s) = \begin{cases} H^{s}(\Omega) \hspace{8cm} \mathrm{if} \; s \in [0, \fr{1}{2}) 
\\ H_{00}^{\fr{1}{2}} (\Omega) = \left\{f \in H_0^\fr{1}{2} (\Omega) : f/\sqrt{d(x)} \in L^2(\Omega)\right\} \hspace{2cm} \mathrm{if} \; s = \fr{1}{2}
\\H_{0}^{s} (\Omega) \hspace{8cm} \mathrm{if} \; s \in (\fr{1}{2},1] 
\end{cases}
\ee where $H_0^{s}(\Omega)$ is the Hilbert subspace of $H^s(\Omega)$ with vanishing boundary trace elements. 

{\bf{Stokes Operator.}} We recall some basic notions of the Stokes operator \cite{cfbook}. We denote by $H$ and $V$ the spaces
\be 
H = \left\{v \in (L^2(\Omega))^2: \na \cdot v = 0, v \cdot n|_{\pa \Omega} = 0 \right\}
\ee where $n$ is the outward unit normal to $\pa \Omega$, and 
\be 
V = \left\{v \in (H_0^1(\Omega))^2 : \na \cdot v = 0\right\}.
\ee 
Let $\PP: (L^2(\Omega))^2 \rightarrow H$ be the Leray Hodge projection. We define the Stokes operator, denoted by $A$, as 
\be 
A:= - \PP \Delta
\ee with domain $\mathcal{D}(A) = V \cap (H^2(\Omega))^2$. $A$ is positive, self-adjoint, and injective, and its inverse $A^{-1}$ is compact. We denote the eigenvalues of $A$ by $\mu_j, j = 1,2,...$, and the corresponding eigenfunctions by $\phi_j, j = 1,2...$, and we note that $0 < \mu_1 \le ... \le \mu_j \le ... \rightarrow \infty$. We define the fractional powers of the Stokes operator, denoted by $A^s$, as
\be \la{stoeig}
A^s v= \sum_{j=1}^{\infty} \mu_j^s (v, \phi_j)_{L^2} \phi_j
\ee with domain 
\be 
\mathcal{D}(A^s) = \left\{v \in H : \|A^{s} v\|_{L^2}^2 := \sum\limits_{j \in \N} \mu_j^{2s}(v, \phi_j)_{L^2}^2 < \infty \right\}.
\ee 

{\bf{Periodic Gevrey classes.}} Let $\TT^2 = [0,2\pi]^2$ be the two dimensional torus.

For $s\in\R$, the periodic fractional Laplacian $\l^s$ applied to a mean zero function $f$ is a Fourier multiplier with symbol $|k|^s$, that is, for $f$ given by
\be
f = \sum\limits_{k \in \mathbb{Z}^2 \setminus \left\{0\right\}} f_k e^{ik \cdot x},
\ee and obeying 
\be 
\sum\limits_{k \in \mathbb{Z}^2 \setminus \left\{0\right\}} |k|^{2s} |f_k|^2 < \infty,
\ee
we have
\be
\l^s f = \sum\limits_{k \in \mathbb{Z}^2 \setminus \left\{0\right\}} |k|^s f_k e^{ik \cdot x}.
\ee 
For $\tau> 0$, $s > 0$, we define 
\be 
e^{\tau \l^s} f = \sum\limits_{k \in \mathbb{Z}^2 \setminus \left\{0\right\}} e^{\tau |k|^s} f_k e^{ik \cdot x}
\ee on 
\be 
\mathcal{D}(e^{\tau \l^s}) = \left\{f \in L^2(\TT^2) : \sum\limits_{k \in \mathbb{Z}^2 \setminus \left\{0\right\}} e^{2\tau |k|^{s}} |f_k|^2 < \infty \right\}.
\ee

\subsection{Results.} We prove first the existence and uniqueness of exponentially decaying in time Sobolev $H^1$ solutions to the unforced model \eqref{inmod11}--\eqref{inmod21} where $f= \Phi = 0$:

\beg{Thm} \la{Existence} Suppose $f= \Phi = 0$. Let $u_0 \in \mathcal{D}(A^{\fr{1}{2}})$ and $q_0 \in \mathcal{D}(\l)$. Then the system \eqref{inmod11}--\eqref{inmod21} has a unique solution $(q, u)$ on $[0, \infty)$ with regularity 
\be \la{t11}
q \in L^{\infty}(0,\infty; \mathcal{D}(\l)) \cap L^2(0,\infty; \mathcal{D}(\l^{1 + \fr{\alpha}{2}}))
\ee and 
\be \la{t12}
u \in \left(L^{\infty} (0, \infty; \mathcal{D}(A^{\fr{1}{2}})) \cap L^2(0,\infty; \mathcal{D}(A))\right)^2.
\ee Moreover, there exists a positive constant $\gamma \le 1$ depending on the size of $\Omega$ and the power $\alpha$, such that the following bounds
\be 
\|\l q(t)\|_{L^2}^2 \le \|\l q_0\|_{L^2}^2 e^{C_0 - \gamma t},
\ee 
\be 
\|A^{\fr{1}{2}} u(t)\|_{L^2}^2
\le C_0 e^{-\gamma t},
\ee
\be 
\int_{0}^{t} \|\l^{1 + \fr{\alpha}{2}} q^{\epsilon}(s)\|_{L^2}^2 ds
\le \|\l q_0\|_{L^2}^2 \left(1 + C_0e^{C_0}\right) 
\ee and 
\be 
\int_{0}^{t} \|Au (s)\|_{L^2}^2 ds 
\le  C_0
\ee hold for all $t \ge 0$, where 
\be 
C_0 =  C\left(\|u_0\|_{L^2}^2 + \|q_0\|_{L^2}^4 + 1 \right)^2 e^{C\left(\|u_0\|_{L^2}^2 + C\|q_0\|_{L^2}^4 \right)^2} \left(\|\na u_0\|_{L^2}^2 + C\|\l q_0\|_{L^2}^2 \|q_0\|_{L^2}^2 \right).
\ee 
\end{Thm}

The solutions of the unforced system \eqref{inmod11}--\eqref{inmod21} are smooth and their higher order derivatives decay exponentially in time to $0$ in all Sobolev norms:

\beg{Thm} \la{tt2}  Let $f = \Phi = 0$. Fix an integer $k \ge 2$. Suppose that $q_0 \in \mathcal{D}(\l^k)$ and $u_0 \in \mathcal{D}(A^{\fr{k}{2}})$. Then the unique solution $(q,u)$ to \eqref{inmod11}--\eqref{inmod21} obey 
\be 
q \in L^{\infty}(0,\infty; \mathcal{D}(\l^k)) \cap L^2(0,\infty; \mathcal{D}(\l^{k + \fr{\alpha}{2}}))
\ee and 
\be
u \in \left(L^{\infty} (0, \infty; \mathcal{D}(A^{\fr{k}{2}})) \cap L^2(0,\infty; \mathcal{D}(A^{\fr{k+1}{2}}))\right)^2.
\ee
Moreover, there is a positive constant $\gamma_k$ depending only on $k$, $\alpha$, $\|\l^k q_0\|_{L^2}$ and $\|A^{\fr{k}{2}}u_0\|_{L^2}$ and a positive constant $c$ depending only on the diameter of $\Omega$ and $\alpha$ such that the estimates
\be 
\|\l^{k} q(t) \|_{L^2}^2 + \|A^{\fr{k}{2}} u(t)\|_{L^2}^2 \le \gamma_k e^{-ct}
\ee and 
\be 
\int_{0}^{t} \left(\|A^{\fr{k+1}{2}} u(s)\|_{L^2}^2 + \|\l^{k+\fr{\alpha}{2}}q(s)\|_{L^2}^2 \right) ds \le \gamma_{k}
\ee
hold for any $t \ge 0$. 
\end{Thm} 

Now we address the long time dynamics of the forced system \eqref{inmod11}--\eqref{inmod21} in the presence of body forces in the fluid and a boundary applied voltage. 

We consider the function spaces
\be \la{functio1}
\mathcal{H} = \mathcal{D}(\l^0) \oplus \mathcal{D}(A^0)
\ee and 
\be \la{functio2}
\mathcal{V} = \mathcal{D}(\l) \oplus \mathcal{D}(A^{\fr{1}{2}}).
\ee The boundary value problem \eqref{inmod11}--\eqref{inmod21} gives rise to a solution map 
\be 
\mathcal{S}(t): \mathcal{V} \mapsto \mathcal{V}
\ee defined by 
\be 
\mathcal{S}(t) (q_0, u_0) = (q(t), u(t)),
\ee where $(q(t), u(t))$ is the unique solution of \eqref{inmod11}--\eqref{inmod21} with initial datum $(q_0, u_0)$ at time $t$. For initial datum $\omega_0 = (q_0, u_0)$, we denote by $\omega(t)$ the solution $(q,u)$ at time $t$ corresponding to $\omega_0$. 

The system \eqref{inmod11}--\eqref{inmod21} has a finite dimensional attractor for any $\alpha \in (0,1]$: 

\begin{Thm}\la{att} Let $\alpha \in \left(0, 1\right]$. There exist a time $T > 0$ and a radius $\tilde{R} > 0$ depending only on the body forces $f$, potential $\Phi$, and the power $\alpha$, such that the ball
\be
\mathcal{B} = \left\{(q, u) \in \mathcal{V} : \|\na q\|_{L^2} + \|\Delta u\|_{L^2} \le R \right\}
\ee obeys $\mathcal{S}(t) \mathcal{B} \subset  \mathcal{B}$ for all $t \ge T$. Moreover, the set
 \begin{equation} \la{X}
X = \bigcap_{t > 0} S(t) \mathcal{B}
\end{equation}
satisfies the following properties: 
\begin{enumerate}
\item[(a)] $X$ is compact in $\mathcal H$.
\item[(b)] $S(t)X = X$ for all $t \geq 0$.
\item[(c)] If $Z$ is bounded in $\mathcal V$  in the norm of of $\mathcal V$, and $S(t)Z = Z$ for all $t \geq 0$, then $Z \subset X$. 
\item[(d)] For every $w_0 \in \mathcal V,$
$\lim\limits_{t \to \infty} dist_{\mathcal H} (S(t)w_0, X) = 0$.
\item[(e)] $X$ is connected.
\item[(f)] $X$ has a finite fractal dimension in $\mathcal{H}$, that is there exists a finite real number $M > 0$ depending on the body forces $f$, potential $\Phi$, and  power $\alpha$ such that 
\[ 
\lim\sup_{r\to 0}\fr{\log{N_{\mathcal{H}}(r)}}{\log\left(\fr{1}{r}\right)} \le M 
\] where $N_{\mathcal{H}}(r)$ is the minimal number of balls in $\mathcal{H}$ of radii $r$ needed to cover $X$.  
\end{enumerate}
\end{Thm}

The existence of the global attractor $X$ is based on the compactness of the ball $\mathcal{B}$ in the norm of $\mathcal{H}$ (Proposition \ref{absball}), the instant Lipschitz continuity in $\mathcal{H}$ of the map $\mathcal{S}(t)$ (Proposition \ref{conts}), and the time analyticity of $\mathcal{S}(t)$ (Proposition \ref{backun}). The finite fractal dimensionality follows from the decay of volume elements transported by the flow map (Proposition \ref{vol}).

When $\alpha = 1$, the attractor is compact in $\mathcal{V}$ and is smooth: 

\begin{Thm}\la{attp} Let $\alpha =1$. There exist a time $\tilde T > 0$ and a radius $\tilde{R} > 0$ depending only on the body forces $f$ and potential $\Phi$ such that the ball
\be
\tilde{\mathcal{B}} = \left\{(q, u) \in \mathcal{V} : \|\l^{1+\fr{\alpha}{2}} q\|_{L^2} + \|\Delta u\|_{L^2} \le R \right\}
\ee obeys $\mathcal{S}(t) \tilde{\mathcal{B}} \subset  \tilde{\mathcal{B}}$ for all $t \ge \tilde T$. 
Moreover, the set \begin{equation} 
\tilde{X} = \bigcap_{t > 0} S(t) \tilde{\mathcal{B}}.
\end{equation}
satisfies the following properties: 
\begin{enumerate}
\item[(a)] $\tilde X$ is compact in $\mathcal V$.
\item[(b)] $S(t)\tilde X = \tilde X$ for all $t \geq 0$.
\item[(c)] If $Z$ is bounded in $\mathcal V$  in the norm of of $\mathcal V$, and $S(t)Z = Z$ for all $t \geq 0$, then $Z \subset \tilde X$. 
\item[(d)] For every $w_0 \in \mathcal V,$
$\lim\limits_{t \to \infty} dist_{\mathcal V} (S(t)w_0, \tilde X) = 0$.
\item[(e)] $\tilde X$ is connected.
\item[(f)] $\tilde X$ has a finite fractal dimension in $\mathcal{V}$, that is there exists a finite real number $\tilde M > 0$ depending on the body forces $f$ and potential $\Phi$ such that 
\[ 
\lim\sup_{r\to 0}\fr{\log{N_{\mathcal{V}}(r)}}{\log\left(\fr{1}{r}\right)} \le \tilde M 
\] where $N_{\mathcal{V}}(r)$ is the minimal number of balls in $\mathcal{V}$ of radii $r$ needed to cover $\tilde X$.
\item[(g)] $\tilde{X}$ is smooth, that is for every integer $k>0$, there exists a radius $\rho_k$ depending on the body forces $f$ and potential $\Phi$ and a ball $B_{\rho_k} \subset H^k$ such that the attractor $\tilde{X} \subset B_{\rho_k}$.  
\end{enumerate}
\end{Thm}

In the case of periodic boundary conditions, the system \eqref{inmod11}--\eqref{inmod111} has unique global Gevrey regular solutions for any fractional dissipative regularity:

\beg{Thm} \la{t3}
Suppose $f=\Phi = 0$. Let $m > 2$. Suppose that $u_0 \in H^{\fr{m}{2}+1}(\TT^2)$ and $q_0 \in H^{\fr{m}{2}}(\TT^2)$. Then there exists a time $T_0$ depending only on the size of the initial data, such that the system described by \eqref{inmod11}--\eqref{inmod111} and equipped with periodic boundary conditions has a unique solution $(q(t),u(t))$ on $(0, T_0)$ with the property that 
\be 
t \mapsto e^{\tau(t) \l^{\fr{\alpha}{2}}} (\l^{\fr{m}{2}}q, \l^{\fr{m}{2} + 1} u )
\ee is analytic on $(0, T_0)$, where
\be 
\tau(t) = \min \left\{\fr{t}{4}, 1, T_0\right\}.
\ee Moreover, $(q,u)$ is analytic on $(T_0, \infty)$ with values in $\mathcal{D}(e^{\sigma \l^{\fr{\alpha}{2}}}\l^{\fr{m}{2}}) \times \mathcal{D}(e^{\sigma \l^{\fr{\alpha}{2}}}\l^{\fr{m}{2}+1})$ for some $\sigma >0$. 
\end{Thm}

The body forces $f$ and potential $\Phi$ are taken to be zero in Theorem \ref{t3} for simplicity. The presence of forcing does not affect the existence, regularity, or analyticity of solutions. 

\subsection{Organization of the Paper.} This paper is organized as follows. In section \ref{sec2}, we prove some identities for fractional powers of the homogeneous Dirichlet Laplacian, derive a nonlinear Poincar\'e inequality in $L^p$ based on C\'ordoba-C\'ordoba inequalities, and recall the Brezis-Mironescu interpolation inequality for fractional powers of the Laplacian on bounded domains. In section \ref{sec3}, we establish commutator estimates and fractional product inequalities for small powers of the Stokes operator based on integral representation formulas and kernel estimates. We present the proof of Theorems \ref{Existence}, \ref{tt2}, \ref{att}, \ref{attp}, and \ref{t3} in Sections \ref{s4}, \ref{S4}, \ref{s5}, \ref{s6}, and \ref{s7} respectively. Finally, we state and prove a spectral lemma describing the asymptotic behavior of eigenvalues associated with a vector valued operator in Appendix \ref{s9} and a Gronwall Lemma describing the long time behavior of solutions to a general differential inequality in Appendix \ref{s8}. These lemmas are frequently used in the paper.

\section{Preliminaries} \la{sec2}

\subsection{Properties of the Fractional Powers of the Laplacian.} 

We recall the identity 
\be 
\lambda^{\fr{s}{2}} = c_s \int_{0}^{\infty} t^{-1-\fr{s}{2}} (1-e^{-t\lambda}) dt
\ee that holds for $s \in (0,2)$ and 
\be 
1 = c_s \int_{0}^{\infty} t^{-1-\fr{s}{2}} (1-e^{-t}) dt,
\ee from which we obtain the integral representation 
\be \la{intrep}
(\l^s f)(x) = c_s \int_{0}^{\infty} [f(x) - e^{t\Delta}f(x)]t^{-1 - \fr{s}{2}} dt
\ee for $f \in \mathcal{D}(\l^s)$ and $s \in (0,2)$, where the heat operator $e^{t\Delta}$ is defined as
\be 
(e^{t\Delta}f)(x) = \int_{\Omega} H_D(x, y, t) f(y) dy
\ee with kernel $H_D(x,y, t)$ given by 
\be 
H_D(x,y,t) = \sum_{j=1}^{\infty} e^{-t\lambda_j} w_j(x) w_j(y).
\la{out}
\ee In 2D, the heat kernel $H_D(x,y,t)$ obeys
\be \la{2dheat1}
|H_D(x,y,t)| \le Ct^{-1} e^{\fr{-|x-y|^2}{kt}},
\ee 
\be \la{2dheat2}
|\na_y H_D(x,y,t)| \le Ct^{-\fr{3}{2}} e^{\fr{-|x-y|^2}{kt}},
\ee and
\be \la{2dheat3}
|\na_x H_D(x,y,t)| \le Ct^{-\fr{3}{2}} e^{\fr{-|x-y|^2}{kt}}
\ee for all $(x,y) \in \Omega \times \Omega$ and $t > 0$. Moreover, the following estimates 
\be \la{2dheat5}
\int_{0}^{\infty} t^{-1-\fr{s}{2}} \int_{\Omega} |x-y|^q |(\na_x + \na_y) H_D(x,y,t)| dy  dt
\le Cd(x)^{-s-1 + q},
\ee
\be \la{2dheat6}
\int_{0}^{\infty} t^{-1-\fr{s}{2}} \int_{\Omega} |x-y|^q |\na_x (\na_x + \na_y) H_D(x,y,t)| dy  dt
\le Cd(x)^{-s-2 + q},
\ee and 
\be \la{2dheat7}
\int_{0}^{\infty} t^{-1-\fr{s}{2}} \int_{\Omega} |x-y|^q |\na_y (\na_x + \na_y) H_D(x,y,t)| dy dt
\le Cd(x)^{-s-2 + q}
\ee
hold for any $q \ge 0 $, $s \in (0,2)$ and $ x \in \Omega$. We refer the reader to \cite{CI,CN2} for detailed proofs of analogous estimates.

\beg{prop} The following identities hold:
\begin{enumerate}
\item[(i)] Let $\alpha, \beta, s \in \R$. For $f \in \mathcal{D}(\l^{\alpha}) \cap \mathcal{D}(\l^{\alpha - s}) $ and $g \in \mathcal{D}(\l^{\beta + s}) \cap \mathcal{D}(\l^{\beta})$, we have
\be 
(\l^{\alpha}f, \l^{\beta}g)_{L^2} = (\l^{\alpha - s}f, \l^{\beta + s}g)_{L^2}.
\ee
\item[(ii)] Let $\alpha, \beta \in \R$. For $f \in \mathcal{D}(\l^{\alpha + 1})$ and $g \in \mathcal{D}(\l^{\beta + 1})$, we have 
\be 
(\l^{\alpha + 1}f, \l^{\beta + 1}g)_{L^2} = (\na \l^{\alpha}f, \na \l^{\beta}g)_{L^2}.
\ee
\item[(iii)] Let $s \in (0,1)$. For $\psi \in \mathcal{D}(\l^{s})$, we have 
\be \la{prodfor1}
\|\l^{s} \psi\|_{L^2}^2
= \int_{\Omega} \int_{\Omega} (\psi(x) - \psi(y))^2 K_{s}(x,y) dxdy 
+ \int_{\Omega} \psi(x)^2 B_{s} dx
\ee where the kernels $K_{s}$ and $B_{s}$ are given by  
\be 
K_s(x,y) := \fr{1}{2c_{2s}} \int_{0}^{\infty} H(x,y,t) t^{-1-s} dt
\ee  for all $x \ne y$, and 
\be 
B_s(x,y) = \fr{1}{c_{2s}} \int_{0}^{\infty} \left[1 - e^{t\Delta}1(x)\right] t^{-1-s}dt.
\ee for all $x\in \Omega$. 
\end{enumerate}
\end{prop}

\noindent \textbf{Proof.} 
\begin{enumerate} 
\item[(i)] The proof of (i) follows from the definition \eqref{maindef}.
\item[(ii)] The proof of (ii) follows from the definition \eqref{maindef} and the identity
\be 
(\na w_j, \na w_k)_{L^2} 
= - (w_j, \Delta w_k)_{L^2}
= (w_j, \lambda_k w_k)_{L^2}
= \begin{cases} \lambda_j \hspace{1cm} \mathrm{if \;} j = k
\\ 0   \hspace{1.2cm} \mathrm{if \;} j \ne k
\end{cases}.
\ee
\item[(iii)] The proof of (iii) is based on \cite{CS}. Indeed, we have 
\beg{align} \la{90}
&-c_{2s} \|\l^{s} \psi\|_{L^2}^2
= -c_{2s} (\l^{2s} \psi, \psi)_{L^2} \nonumber
\\&= \int_{0}^{\infty} \int_{\Omega} t^{-1-s} \left[\int_{\Omega} H(x,y,t) \psi(x)\psi(y) dx  - \psi(y)^2 \right] dy dt \nonumber
\\&= \int_{0}^{\infty} \int_{\Omega} \int_{\Omega} t^{-1-s} H(x,y,t) (\psi(x) - \psi(y)) \psi(y) dx dy dt \nonumber
\\&\quad\quad+ \int_{0}^{\infty} \int_{\Omega} t^{-1-s} \psi(y)^2 \left[e^{t\Delta}1(y) - 1 \right] dy dt 
\end{align} in view of the integral representation \eqref{intrep}, and 
\beg{align} \la{91}
-c_{2s} \|\l^{s} \psi\|_{L^2}^2 \nonumber
&= -\int_{0}^{\infty} \int_{\Omega} \int_{\Omega} t^{-1-s} H(x,y,t) (\psi(x) - \psi(y)) \psi(x) dx dy dt
\\&\quad\quad+ \int_{0}^{\infty} \int_{\Omega} t^{-1-s} \psi(y)^2 \left[e^{t\Delta}1(y) - 1 \right] dy dt
\end{align} by interchanging the variables $x$ and $y$ in the first integral in \eqref{90} and using the symmetry of the heat kernel $H_D(x,y,t)$. Adding \eqref{90} and \eqref{91}, we deduce that 
\beg{align} \la{92}
-2c_{2s} \|\l^{s} \psi\|_{L^2}^2 
&= -\int_{0}^{\infty} \int_{\Omega} \int_{\Omega} t^{-1-s} H(x,y,t) (\psi(x) - \psi(y))^2 dx dy dt \nonumber
\\&\quad\quad-2\int_{0}^{\infty} \int_{\Omega} t^{-1-s} \psi(y)^2 \left[1 - e^{t\Delta} 1(y) \right] dydt.
\end{align} Dividing both sides of \eqref{92} by $-2c_{2s}$ and applying Fubini's theorem, we obtain \eqref{prodfor1}.
\end{enumerate}

\beg{rem} The kernels $K_{s}$ and $B_{s}$ obey
\be \la{prodfor2}
0 \le K_{s}(x,y) \le \fr{C_s}{|x-y|^{2 + 2s}}
\ee for all $x \ne y$, and 
\be \la{prodfor3}
B_{s}(x) \ge 0
\ee for all $x\in \Omega$. The estimate \eqref{prodfor2} follows from \eqref{2dheat1}, whereas the nonnegativity of $B_{s}$ follows from the maximum principle. 
\end{rem}

\beg{prop} \la{cru} For any odd integer $m \ge 1$, we have
\be 
\mathcal{D}(\l^m) \cap H^{m+1} = \mathcal{D}(\l^{m+1}).
\ee
\end{prop}

\noindent \textbf{Proof.} The inclusion $\mathcal{D}(\l^{m+1}) \subset \mathcal{D}(\l^m) \cap H^{m+1}$ obviously holds. If $\rho \in \mathcal{D}(\l^m) \cap H^{m+1}$, then $\l^k \rho$ vanishes on the boundary for all even integers $k \le m -1$ and consequently, $\rho \in \mathcal{D}(\l^{m+1})$.

\subsection{Nonlinear Poincar\'e inequality.}  We recall the following 
pointwise inequality in bounded domains \cite{CI}:

\beg{prop} \la{corcor} Let $0 \le s < 2$. There exists a constant $c > 0$ depending only on the domain  $\Omega$ and on $s$, such that, for any $\Phi$, a $C^2$ convex function satisfying $\Phi(0) = 0$, and any function $f \in C_0^{\infty}(\Omega)$, the inequality 
\be 
\Phi'(f) \l^s f - \l^s (\Phi(f)) \ge \fr{c}{d(x)^s} (f \Phi'(f) - \Phi(f))
\ee holds pointwise in $\Omega$.
\end{prop}

For an even integer $p \ge 2$, we let $\Phi(x) = \fr{1}{p} x^p$, and we apply Proposition \ref{corcor} to infer that
\be 
f^{p-1} \l^s f \ge \fr{1}{p} \l^{s} (f^p) + \fr{c}{d(x)^s} \left(1 - \fr{1}{p} \right) f^p
\ee for any $f \in C_0^{\infty}(\Omega)$. Integrating over $\Omega$, we have
\be 
\int_{\Omega} f^{p-1} \l^s f dx \ge \fr{1}{p} \int_{\Omega} \l^s (f^p)  dx + C_{\Omega, s} \left(1 - \fr{1}{p}\right) \|f\|_{L^p}^p
\ee for some positive constant $C_{\Omega,s }$ depending only on the size of $\Omega$ and $s$. In view of the integral representation formula \eqref{intrep}, the maximum principle, and the positivity of $f^p$, we deduce that
\be 
\int_{\Omega} \l^s (f^p)  dx \ge 0.
\ee This yields the $L^p$ nonlinear Poincar\'e inequality 
\be \la{pci}
\int_{\Omega} f^{p-1} \l^s f dx \ge C_{\Omega, s} \left(1 - \fr{1}{p}\right) \|f\|_{L^p}^p.
\ee

\subsection{Brezis-Mironescu interpolation inequality.} We define the fractional spaces $W^{s, p}(\Omega)$ as 
\be 
W^{s, p} (\Omega) = \left\{ v \in L^p(\Omega): \|v\|_{W^{s, p}}= \left(\|v\|_{L^p}^p + \int_{\Omega} \int_{\Omega} \fr{|v(x) - v(y)|^p}{|x-y|^{2+sp}} dxdy \right)^{\fr{1}{p}}<\infty \right\}.
\ee Let $1 \le p, p_1, p_2 \le \infty$ with $p_2 \ne 1$. Let $s, s_1, s_2$ be nonnegative real numbers such that $s_1 \le s_2$. Let $\theta \in (0,1)$ such that $s = \theta s_1 + (1-\theta) s_2$ and $\fr{1}{p} = \fr{\theta}{p_1} + \fr{1-\theta}{p_2}$. Then there exists a positive universal constant $C$ such that the following interpolation inequality
\be 
\|f\|_{W^{s,p}} \le C\|f\|_{W^{s_1, p_1}}^{\theta}\|f\|_{W^{s_2, p_2}}^{1-\theta}
\ee holds for any $f \in W^{s_1, p_1}(\Omega) \cap W^{s_2, p_2} (\Omega)$. We refer the reader to \cite{BM} for a detailed proof.

\subsection{Rellich's inequality} We denote by $W_0^{1,p}(\Omega)$ the closure of the space of smooth compactly supported functions $C_0^{\infty}(\Omega)$ under the norm of $W^{1,p}$. 
For $1 \le p < \infty$ and $f \in W_0^{1,p}(\Omega)$, the following Rellich inequality holds:
\be 
\int_{\Omega} \fr{|f(x)|^p}{d(x)^p} dx \le \int_{\Omega} |\na f(x)|^p dx.
\ee (see \cite{EH} and references therein).

\subsection{Notation.} Throughout the paper, we denote by $C$ a positive constant that depends on the domain $\Omega$ and universal constants. The distance from a point $x \in \Omega$ to the boundary $\pa \Omega$ is denoted by $d(x)$. The notation $[A,B]$ is used to denote the commutator $AB - BA$. 

\section{Commutator and Fractional Product Estimates} \la{sec3}

In this section, we state and prove two-dimensional commutator and fractional product estimates that will be used to control the nonlinearities of our model.

For $0 \le \gamma \le 1$, we denote by $C^{0,\gamma}(\bar \Omega)$ the H\"older space with norm
\be 
\|\tilde{q}\|_{C^{0,\gamma}} = \|\tilde{q}\|_{L^{\infty}} + [\tilde{q}]_{C^{0,\gamma}}
\ee where
\be 
[\tilde{q}]_{C^{0,\gamma}} = \sup\limits_{x \ne y} \fr{|\tilde{q}(x) - \tilde{q}(y)|}{|x-y|^{\gamma}}.
\ee

\beg{prop} \la{coes1} Let $s \in (0,1)$, $\gamma \in [0,1]$, and $s < \gamma$. Suppose $\tilde u \in C^{0, \gamma}$. The operator $[\l^s, \tilde{u}]$ can be uniquely extended from $C_0^{\infty}(\Omega)$ to $L^2(\Omega$) such that
\be \la{coes12}
\|[\l^s, \tilde{u}] \tilde{q}\|_{L^2} \le C[\tilde u]_{C^{0,\gamma}} \|\tilde{q}\|_{L^2}
\ee holds for any $\tilde{q} \in L^2$.
\end{prop}

\noindent \textbf{Proof.} The estimate \eqref{coes12} is a particular case of  Theorem 2.6 in \cite{CN1}.

\beg{prop}\la{coes2} Let $s \in (0,2)$ and $p \in (2, \infty]$. Let $\tilde{q} \in C_0^{\infty}(\Omega)$. Then the estimate 
\be \la{coes21}
|[\na \na, \l^{s}] \tilde{q}(x)| \le C\left(\|\tilde{q}\|_{W^{1,p}} d(x)^{-s-1-\fr{2}{p}} + |\tilde{q}(x)| d(x)^{-s-2}\right)
\ee holds for all $x \in \Omega$.
\end{prop}

\noindent \textbf{Proof.} Using the integral representation formula \eqref{intrep} and integrating by parts, we have 
\be 
|[\na \na, \l^{s}] \tilde{q}(x)|
= c_s\left|\int_{0}^{\infty} t^{-1 - \fr{s}{2}} \int_{\Omega} (\na_x \na_x - \na_y \na_y) H_D(x,y,t) \tilde{q}(y) dy dt \right|,
\ee  which can be bounded as
\beg{align}
|[\na \na, \l^{s}] \tilde q(x)|
&\le c_s\left|\int_{0}^{\infty} t^{-1 - \fr{s}{2}} \int_{\Omega} (\na_x \na_x + \na_x \na_y) H_D(x,y,t) \tilde q(y) dy dt \right| \nonumber
\\&\quad\quad+  c_s\left|\int_{0}^{\infty} t^{-1 - \fr{s}{2}} \int_{\Omega} (\na_y \na_x + \na_y \na_y) H_D(x,y,t) \tilde q(y) dy dt \right|
\end{align} via a direct application of the triangle inequality. Subtracting and adding $\tilde{q}(x)$, the latter inequality yields 
\beg{align} \la{exiu1}
|[\na \na, \l^{s}] \tilde q(x)|
&\le C\int_{0}^{\infty} t^{-1 - \fr{s}{2}} \int_{\Omega} |\na_x (\na_x +  \na_y) H_D(x,y,t)| |\tilde q(y)- \tilde{q}(x)| dy dt \nonumber
\\&\quad\quad+ C|\tilde{q}(x)| \int_{0}^{\infty} t^{-1 - \fr{s}{2}} \int_{\Omega} |\na_x (\na_x + \na_y) H_D(x,y,t)| dy dt \nonumber
\\&\quad\quad\quad+C\int_{0}^{\infty} t^{-1 - \fr{s}{2}} \int_{\Omega} |\na_y (\na_x + \na_y) H_D(x,y,t)| |\tilde q(y) - \tilde{q}(x)| dydt \nonumber
\\&\quad\quad\quad\quad+ C|\tilde{q}(x)|\int_{0}^{\infty} t^{-1 - \fr{s}{2}} \int_{\Omega} |\na_y (\na_x + \na_y) H_D(x,y,t)| dydt.
\end{align} In view of the heat kernel estimate \eqref{2dheat6}, we bound
\beg{align}
&\int_{0}^{\infty} t^{-1 - \fr{s}{2}} \int_{\Omega} |\na_x (\na_x +  \na_y) H_D(x,y,t)| |\tilde q(y)- \tilde{q}(x)| dy dt \nonumber
\\&\quad\quad\le C[q]_{C^{0, 1-\fr{2}{p}}} \int_{0}^{\infty} t^{-1 - \fr{s}{2}} \int_{\Omega} |x-y|^{1-\fr{2}{p}} |\na_x (\na_x +  \na_y) H_D(x,y,t)|  dy dt \nonumber
\\&\quad\quad\le C[q]_{C^{0, 1-\fr{2}{p}}} d(x)^{-s-1-\fr{2}{p}}
\end{align} and
\be 
 \int_{0}^{\infty} t^{-1 - \fr{s}{2}} \int_{\Omega} |\na_x (\na_x + \na_y) H_D(x,y,t)| dy dt 
 \le Cd(x)^{-s-2}.
\ee By making use of the heat kernel estimate \eqref{2dheat7}, we estimate 
\be 
\int_{0}^{\infty} t^{-1 - \fr{s}{2}} \int_{\Omega} |\na_x (\na_x +  \na_y) H_D(x,y,t)| |\tilde q(y)- \tilde{q}(x)| dy dt 
\le C[q]_{C^{0, 1-\fr{2}{p}}} d(x)^{-s-1-\fr{2}{p}}
\ee and 
\be \la{exiu2}
\int_{0}^{\infty} t^{-1 - \fr{s}{2}} \int_{\Omega} |\na_y (\na_x + \na_y) H_D(x,y,t)| dydt
\le Cd(x)^{-s-2}.
\ee Putting \eqref{exiu1}--\eqref{exiu2} together and using the two-dimensional continuous embedding of the Sobolev space $W^{1,p}$ into the H\"older space $C^{0, 1-\fr{2}{p}}(\bar{\Omega})$, we obtain \eqref{coes21}, ending the proof of Proposition \ref{coes2}.

\beg{cor} \la{exiu55} Let $\alpha \in (0,1]$. Let $p \in (2, \infty)$ and $\epsilon > 0$ such that 
\be 
r := 2 - \alpha - 2\alpha \epsilon - \fr{8}{p} - \fr{8\epsilon}{p} >0.
\ee Fix $\tilde{u} \in W_0^{1,\fr{4}{r}} \cap W_0^{1, \fr{16}{\alpha}}$, and define the numbers $p_0 = \max\left\{p, \fr{8}{4-3\alpha} \right\}$ and $r_0 = \max\left\{\fr{4}{r}, \fr{16}{\alpha} \right\}$. The operator $\tilde{u} \cdot [\na \na, \l^{\fr{\alpha}{2}}]$ can be uniquely extended from $C_0^{\infty}(\Omega)$ to $W_0^{1,p_0}$ such that the estimate 
\be \la{exiu5}
\|\tilde{u} \cdot [\na \na, \l^{\fr{\alpha}{2}}] \tilde{q}\|_{L^{\fr{4}{2+\alpha}}} \le C \|\tilde{u}\|_{W^{1,r_0}} \|\tilde{q}\|_{W^{1,p_0}}
\ee holds for any $\tilde{q} \in W_0^{1,p_0}$. 
\end{cor} 

\noindent \textbf{Proof.} Fix $\tilde{q} \in C_0^{\infty}(\Omega)$. In view of Proposition \ref{coes2} with $s= \fr{\alpha}{2}$, we have 
\be \la{exiu6}
|\tilde{u} \cdot [\na \na, \l^{\fr{\alpha}{2}}] \tilde{q}(x)|
\le  C\|\tilde{q}\|_{W^{1,p}} |\tilde{u}(x)| d(x)^{-1-\fr{\alpha}{2} - \fr{2}{p}} + C |\tilde{u}(x)| |\tilde{q}(x)| d(x)^{-2-\fr{\alpha}{2}}.
\ee Since $\tilde{u} \in W_0^{1,\fr{16}{\alpha}}$ and $\tilde{q} \in W_0^{1,\fr{8}{4-3\alpha}}$, we can apply Rellich's inequality to control $\tilde{u}(\cdot) d(\cdot)^{-1}$ and $\tilde{q}(\cdot) d(\cdot)^{-1}$ as follows,
\be 
\|\tilde{u}(\cdot) d(\cdot)^{-1}\|_{L^{\fr{16}{\alpha}}} \le C\|\na \tilde{u}\|_{L^{\fr{16}{\alpha}}}
\ee and 
\be 
\|\tilde{q}(\cdot) d(\cdot)^{-1}\|_{L^{\fr{8}{4-3\alpha}}} \le C\|\na \tilde{q}\|_{L^{\fr{8}{4-3\alpha}}}.
\ee Using H\"older's inequality with exponents $\fr{8}{4-3\alpha}$, $\fr{16}{\alpha}$, and $\fr{16}{9\alpha}$ and the fact that powers of the distance to the boundary function $d(x)^{-\beta}$ are space integrable for $\beta \in [0,1)$, we estimate
\be \la{exiu7}
\|\tilde{u} \tilde{q} d(\cdot)^{-2-\fr{\alpha}{2}}\|_{L^{\fr{4}{2+\alpha}}}
\le \|\tilde{u}(\cdot) d(\cdot)^{-1}\|_{L^{\fr{16}{\alpha}}} \|\tilde{q}(\cdot) d(\cdot)^{-1}\|_{L^{\fr{8}{4-3\alpha}}}  \|d(\cdot)^{-\fr{\alpha}{2}}\|_{L^{\fr{16}{9\alpha}}}
\le C\|\na \tilde{u}\|_{L^{\fr{16}{\alpha}}}\|\na \tilde{q}\|_{L^{\fr{8}{4-3\alpha}}}. 
\ee Another application of Rellich's inequality yields
\be 
\|\tilde{u}(\cdot) d(\cdot)^{-1}\|_{L^{\fr{4}{r}}} \le C\|\na \tilde{u}\|_{L^{\fr{4}{r}}} 
\ee from which we obtain
\be  \la{exiu8}
\|\tilde{u} d(\cdot)^{-1-\fr{\alpha}{2} - \fr{2}{p}}\|_{L^{\fr{4}{2+\alpha}}} \le \|\tilde{u}(\cdot) d(\cdot)^{-1}\|_{L^{\fr{4}{r}}}  \|d(\cdot)^{-\fr{\alpha}{2} - \fr{2}{p}}\|_{L^{\fr{2p}{(\alpha p  + 4)(1+\epsilon)}}} \le C\|\na \tilde{u}\|_{L^{\fr{4}{r}}}
\ee after a direct application of H\"older's inequality with exponents $\fr{4}{r}$ and $\fr{2p}{(\alpha p  + 4)(1+\epsilon)}$. We note that the H\"older exponent $r$ is chosen in such a way that optimizes the value of $p$ for which $ \|d(\cdot)^{-\fr{\alpha}{2} - \fr{2}{p}}\|_{L^{\fr{2p}{(\alpha p  + 4)(1+\epsilon)}}}<\infty$. Finally, we combine \eqref{exiu6}, \eqref{exiu7} and \eqref{exiu8}, use the density of $C_0^{\infty}(\Omega)$ in $W_0^{1, p_0}$, and extend by continuity to obtain \eqref{exiu5} for all $\tilde{q} \in W_0^{1,p_0}$.

\beg{prop} \la{coes3} Let $s \in (0,2)$, $\beta \in [0,1)$, and $p \in (2, \infty]$. Let $\tilde{u} \in  W_0^{1, \fr{2}{\beta}}$ be divergence-free and $\tilde{q} \in C_0^{\infty}(\Omega)$. Then it holds that
\beg{align} \la{coes31}
&|[\na, \l^{s}](\tilde{u} \cdot \na \tilde{q})(x) |  \nonumber
\\&\quad\le C\left(\|u\|_{W^{1,\fr{2}{\beta}}} \|\tilde{q}\|_{W^{1,p}} d(x)^{-s- \beta - \fr{2}{p}} + |\tilde{u}(x)| \|\tilde{q}\|_{W^{1,p}} d(x)^{-s-1 -\fr{2}{p}}  + |\tilde{u}(x)| |\tilde{q}(x)| d(x)^{-s-2}\right)
\end{align}  for a.e. $x \in \Omega$.
\end{prop}

\noindent \textbf{Proof.} The pointwise integral representation formula of the commutator $[\na, \l^{s}](\tilde{u} \cdot \na \tilde{q})$ is given by 
\be 
(\na \l^{s} (\tilde u \cdot \na \tilde q) - \l^{s} \na (\tilde u  \cdot \na \tilde q))(x)
= c_s \int_{0}^{\infty} t^{-1-\fr{s}{2}} \int_{\Omega} (\na_x + \na_y) H_D(x,y,t) \na_y \cdot (\tilde u(y)  \tilde q(y)) dydt,
\ee which, after integration by parts, reduces to  
\be 
(\na \l^{s} (\tilde u \cdot \na \tilde q) - \l^{s} \na (\tilde u  \cdot \na \tilde q))(x)
= -c_s \int_{0}^{\infty} t^{-1-\fr{s}{2}} \int_{\Omega} (\na_y (\na_x + \na_y) H_D(x,y,t)) \cdot \tilde u(y)  \tilde q(y) dydt.
\ee 
Subtracting and adding $\tilde u(x)$ and $\tilde q(x)$ and using the divergence-free condition obeyed by $\tilde{u}$, we obtain
\beg{align} \la{coes32}
&(\na \l^{s} (\tilde u  \cdot \na \tilde q) - \l^{s} \na (\tilde u  \cdot \na \tilde q))(x) \nonumber
\\&= -c_s \int_{0}^{\infty} t^{-1-\fr{s}{2}} \int_{\Omega} (\na_y (\na_x + \na_y) H_D(x,y,t))  \cdot (\tilde u(y) - \tilde u(x)) (\tilde q(y) - \tilde{q}(x))dy dt \nonumber
\\&\quad\quad-  c_s \int_{0}^{\infty} t^{-1-\fr{s}{2}} \int_{\Omega} (\na_y (\na_x + \na_y) H_D(x,y,t))  \cdot (\tilde u(x) \tilde q(y)) dy dt \nonumber
\\&\quad\quad\quad+   c_s \tilde u(x) \cdot \int_{0}^{\infty} t^{-1-\fr{s}{2}} \int_{\Omega} (\na_y (\na_x + \na_y) H_D(x,y,t))  \cdot (\tilde u(x) \tilde q(x)) dy dt \nonumber
\\&:= A_1(x) + A_2(x) + A_3(x).
\end{align} 
In view of the heat kernel estimate \eqref{2dheat7}, we estimate
\beg{align}
|A_1(x)| 
&\le C[\tilde{u}]_{C^{0, 1-\beta}} [\tilde{q}]_{C^{0, 1-\fr{2}{p}}} \int_{0}^{\infty} t^{-1 - \fr{s}{2}} \int_{\Omega} |x-y|^{2-\beta - \fr{2}{p}}|\na_y (\na_x + \na_y) H_D(x,y,t)| dy dt  \nonumber
\\&\le C[\tilde{u}]_{C^{0, 1-\beta}} [\tilde{q}]_{C^{0, 1-\fr{2}{p}}} d(x)^{-s-\beta - \fr{2}{p}},
\end{align}
\beg{align}
|A_2(x)|
&\le C|\tilde{u}(x)| [\tilde{q}]_{C^{0, 1-\fr{2}{p}}}  \int_{0}^{\infty} t^{-1 - \fr{s}{2}} \int_{\Omega} |x-y|^{1 - \fr{2}{p}}|\na_y (\na_x + \na_y) H_D(x,y,t)| dy dt  \nonumber
\\&\le C|\tilde{u}(x)| [\tilde{q}]_{C^{0, 1-\fr{2}{p}}} d(x)^{-s-1-\fr{2}{p}}
\end{align} and 
\beg{align}
|A_3(x)|  \la{Exui10}
&\le  C|\tilde{u}(x)| |\tilde{q}(x)| \int_{0}^{\infty} t^{-1 - \fr{s}{2}} \int_{\Omega} |\na_y (\na_x + \na_y) H_D(x,y,t)| dy dt  \nonumber  \nonumber
\\&\le C|\tilde{u}(x)| |\tilde{q}(x)| d(x)^{-s-2}.
\end{align} Putting \eqref{coes32}--\eqref{Exui10} together and using the continuous embeddings of $W^{\fr{2}{\beta}}$ into $C^{0, 1-\beta} $ and $W^{1,p}$ into $C^{0, 1-\fr{2}{p}}$, we obtain \eqref{coes31}. This finishes the proof of Proposition \ref{coes3}. 

\beg{cor} \la{EXIU} Let $\alpha \in (0,1]$. Let $p \in (2, \infty)$, $\epsilon > 0$, and $\beta>0$ such that 
\be 
r := 2 - \alpha - 2\alpha \epsilon - \fr{8}{p} - \fr{8\epsilon}{p} >0
\ee and 
\be \la{compp1}
\beta < \fr{1}{2} - \fr{\alpha}{4} - \fr{2}{p}.
\ee
Fix $\tilde{u} \in W_0^{1,\fr{4}{r}} \cap W_0^{1, \fr{16}{\alpha}} \cap W_0^{1,\fr{2}{\beta}}$, and define the numbers $p_0 = \max\left\{p, \fr{8}{4-3\alpha} \right\}$ and $r_0 = \max\left\{\fr{4}{r}, \fr{16}{\alpha}, \fr{2}{\beta} \right\}$. The operator $[\na,\l^{\fr{\alpha}{2}}] \tilde{u} \cdot \na$ can be uniquely extended from $C_0^{\infty}(\Omega)$ to $W_0^{1,p_0}$ such that the estimate 
\be \la{exiu5o}
\|[\na, \l^{\fr{\alpha}{2}}](\tilde{u} \cdot \na \tilde{q})\|_{L^{\fr{4}{2+\alpha}}} \le C \|\tilde{u}\|_{W^{1,r_0}} \|\tilde{q}\|_{W^{1,p_0}}
\ee holds for any $\tilde{q} \in W_0^{1,p_0}$. 
\end{cor}

\noindent \textbf{Proof.} Let $\tilde{q} \in C_0^{\infty}(\Omega)$. In view of \eqref{coes31} with $s= \fr{\alpha}{2}$, we have
\beg{align} \la{compp}
&|[\na, \l^{\fr{\alpha}{2}}](\tilde{u} \cdot \na \tilde{q})(x) |  \nonumber
\\&\quad\le C\left(\|u\|_{W^{1,\fr{2}{\beta}}} \|\tilde{q}\|_{W^{1,p}} d(x)^{-\fr{\alpha}{2}- \beta - \fr{2}{p}} + |\tilde{u}(x)| \|\tilde{q}\|_{W^{1,p}} d(x)^{-\fr{\alpha}{2}-1 -\fr{2}{p}}  + |\tilde{u}(x)| |\tilde{q}(x)| d(x)^{-\fr{\alpha}{2}-2}\right).
\end{align} We apply the $L^{\fr{4}{2+\alpha}}$ norm. The second and third terms on the right hand side of \eqref{compp} are estimated as in Corollary \ref{exiu55}. As for the first term, we use the condition \eqref{compp1} which guarantees that $\left(\fr{\alpha}{2}+ \beta + \fr{2}{p} \right) \left(\fr{4}{2+\alpha} \right) < 1 $ and infer that this latter power of the distance to the boundary function is integrable. By the density of $C_{0}^{\infty}(\Omega)$ in $W_0^{1,p_0}$ and extension by continuity, we obtain \eqref{exiu5o}.

\beg{prop} Let $s \in (0,2)$ and $p \in (2, \infty]$. Let $\tilde{q} \in C_0^{\infty}(\Omega)$. Then the estimate 
\be \la{coes21p}
|[\na, \l^{s}] \tilde{q}(x)| \le C\left(\|\tilde{q}\|_{W^{1,p}} d(x)^{-s-\fr{2}{p}} + |\tilde{q}(x)| d(x)^{-s-1}\right)
\ee holds for all $x \in \Omega$.
\end{prop}

\noindent \textbf{Proof.} Using the integral representation formula \eqref{intrep} and integrating by parts, we have 
\be 
|[\na, \l^{s}] \tilde{q}(x)|
= c_s\left|\int_{0}^{\infty} t^{-1 - \fr{s}{2}} \int_{\Omega} (\na_x + \na_y) H_D(x,y,t) \tilde{q}(y) dy dt \right|,
\ee  which, after subtracting and adding $\tilde{q}(x)$, reduces to
\beg{align}
|[\na, \l^{s}] \tilde q(x)|
&\le C\int_{0}^{\infty} t^{-1 - \fr{s}{2}} \int_{\Omega} |(\na_x +  \na_y) H_D(x,y,t)| |\tilde q(y)- \tilde{q}(x)| dy dt \nonumber
\\&\quad\quad+ C|\tilde{q}(x)| \int_{0}^{\infty} t^{-1 - \fr{s}{2}} \int_{\Omega} |(\na_x + \na_y) H_D(x,y,t)| dy dt 
\end{align} In view of the heat kernel estimate \eqref{2dheat5}, we bound
\beg{align}
&\int_{0}^{\infty} t^{-1 - \fr{s}{2}} \int_{\Omega} |(\na_x +  \na_y) H_D(x,y,t)| |\tilde q(y)- \tilde{q}(x)| dy dt \nonumber
\\&\quad\quad\le C[q]_{C^{0, 1-\fr{2}{p}}} \int_{0}^{\infty} t^{-1 - \fr{s}{2}} \int_{\Omega} |x-y|^{1-\fr{2}{p}} |(\na_x +  \na_y) H_D(x,y,t)|  dy dt \nonumber
\\&\quad\quad\le C\|q\|_{W^{1,p}} d(x)^{-s-\fr{2}{p}}
\end{align} and
\be 
 \int_{0}^{\infty} t^{-1 - \fr{s}{2}} \int_{\Omega} |(\na_x + \na_y) H_D(x,y,t)| dy dt 
 \le Cd(x)^{-s-1}.
\ee This gives \eqref{coes21p}.

\beg{cor} \la{exiu55i} Let $\alpha \in (0,1]$. Let $p \in (2, \infty)$ and $\epsilon > 0$ such that 
\be 
r := 2 - \alpha - 2\alpha \epsilon - \fr{8}{p} - \fr{8\epsilon}{p} >0.
\ee Fix $\tilde{u} \in W^{1,\fr{4}{r}} \cap W^{1, \fr{16}{\alpha}}$, and define the numbers $p_0 = \max\left\{p, \fr{8}{4-3\alpha} \right\}$ and $r_0 = \max\left\{\fr{4}{r}, \fr{16}{\alpha} \right\}$. The operator $\na \tilde{u} \cdot [\na, \l^{\fr{\alpha}{2}}]$ can be uniquely extended from $C_0^{\infty}(\Omega)$ to $W_0^{1,p_0}$ such that the estimate 
\be 
\|\na \tilde{u} \cdot [\na, \l^{\fr{\alpha}{2}}] \tilde{q}\|_{L^{\fr{4}{2+\alpha}}} \le C \|\tilde{u}\|_{W^{1,r_0}} \|\tilde{q}\|_{W^{1,p_0}}
\ee holds for any $\tilde{q} \in W_0^{1,p_0}$.
\end{cor} 

\noindent \textbf{Proof.} Let $\tilde{q} \in C_0^{\infty}(\Omega)$. We apply \eqref{coes21p} with $s=\fr{\alpha}{2}$ and obtain
 \be 
|[\na, \l^{\fr{\alpha}{2}}] \tilde{q}(x)| \le C\left(\|\tilde{q}\|_{W^{1,p}} d(x)^{-\fr{\alpha}{2} - \fr{2}{p}} + |\tilde{q}(x)| d(x)^{-1-\fr{\alpha}{2}}\right).
\ee By H\"older and and Rellich inequalities, we have
\beg{align}
\|\na \tilde{u} \cdot [\na, \l^{\fr{\alpha}{2}}] \tilde{q}\|_{L^{\fr{4}{2+\alpha}}} \nonumber
&\le C\|\na \tilde{u}\|_{L^{\fr{4}{r}}}\|\tilde{q}\|_{W^{1,p}} \|d(x)^{-\fr{\alpha}{2}- \fr{2}{p}}\|_{L^{\fr{2p}{(\alpha p  + 4)(1+\epsilon)}}} 
\\&\quad\quad+\|\na \tilde{u} \|_{L^{\fr{16}{\alpha}}} \|\tilde{q}(\cdot)d(\cdot)^{-1}\|_{L^{\fr{8}{4-3\alpha}}} \|d(\cdot)^{-\fr{\alpha}{2}}\|_{L^{\fr{16}{9\alpha}}} \nonumber
\\&\le C \|\tilde{u}\|_{W^{1,r_0}} \|\tilde{q}\|_{W^{1,p_0}}.
\end{align} This completes the proof of Corollary \ref{exiu55i}.

\beg{prop} \la{pr3} Let $m \ge 1$ be an integer. Let $\tilde{u} = (\tilde{u}_1, \tilde{u}_2) \in H^{m+1}$ be a two-dimensional vector field and $\tilde{q} \in \mathcal{D}(\l^{m+\fr{\alpha}{2}})$ be a scalar function. If $m$ is even and $\tilde{u} \cdot \na \tilde{q} \in \mathcal{D}(\l^m)$, then it holds that 
\be \la{pr31}
\|[\l^m, \tilde{u} \cdot \na] \tilde{q}\|_{L^2} \le C\|\tilde{u}\|_{H^{m+1}} \|\l^{m+ \fr{\alpha}{2}} \tilde{q}\|_{L^2}.
\ee  If $m$ is odd and $\tilde{u} \cdot \na \tilde{q} \in \mathcal{D}(\l^{m-1})$, then it holds that 
\be \la{pr32}
\|[\na \l^{m-1}, \tilde{u} \cdot \na] \tilde{q}\|_{L^2} \le C\|\tilde{u}\|_{H^{m+1}} \|\l^{m+ \fr{\alpha}{2}} \tilde{q}\|_{L^2}.
\ee Here, $C$ is a positive constant depending only on $m$ and $\alpha$. 
\end{prop}

\noindent \textbf{Proof.} We present a proof by induction. Suppose $m=2$. Let $\tilde{u} = (\tilde{u}_1, \tilde{u}_2) \in H^{3}$ and $\tilde{q} \in \mathcal{D}(\l^{2+\fr{\alpha}{2}})$  such that $\tilde{u} \cdot \na \tilde{q} \in \mathcal{D}(\l^2)$. Since $\l^2 = -\Delta$, the commutator $(-\Delta) (\tilde{u} \cdot \na \tilde{q}) - \tilde{u} \cdot \na (-\Delta) \tilde{q}$ reduces to $(-\Delta) \tilde{u} \cdot \na \tilde{q} - 2  \na \tilde{u} \cdot  \na \na \tilde{q}$, where
\be 
\na \tilde{u} \cdot \na \na \tilde{q}:= \sum\limits_{i=1}^{2} \na \tilde{u}_i \cdot \na \pa_{x_i}\tilde{q},
\ee hence its $L^2$ norm can be bounded as 
\be \la{Bcom}
\|(-\Delta) (\tilde{u} \cdot \na \tilde{q}) - \tilde{u} \cdot \na (-\Delta) \tilde{q} \|_{L^2}
\le C\left(\|\Delta \tilde{u}\|_{L^q} \|\na \tilde{q}\|_{L^p} + \| \na \tilde{u}\|_{L^q} \|\na \na \tilde{q}\|_{L^p}\right)
\ee for any $p,q \in [1,\infty]$ obeying $\fr{1}{p} + \fr{1}{q} = \fr{1}{2}$. Using the Gagliardo-Nirenberg inequalities, and choosing $p$ so that $\mathcal{D}(\l^{\fr{\alpha}{2}})$ is continuously embedded in $L^p$, we have 
\be
\|(-\Delta) (\tilde{u} \cdot \na \tilde{q}) - \tilde{u} \cdot \na (-\Delta) \tilde{q} \|_{L^2}
\le C\|\tilde{u}\|_{H^3}\|\l^{2 + \fr{\alpha}{2}} \tilde{q}\|_{L^2},
\ee which gives \eqref{pr31} for $m=2$. Suppose that \eqref{pr31} holds for an even integer $m$, any field $\tilde{u} \in H^{m+1}$ and any scalar $\tilde{q} \in \mathcal{D}(\l^{m+\fr{\alpha}{2}})$ with the property $\tilde{u} \cdot \na \tilde{q} \in \mathcal{D}(\l^m)$. We show that 
\be \la{pr33}
\|\l^{m+2} (\tilde{u} \cdot \na \tilde{q}) - \tilde{u} \cdot \na \l^{m+2} \tilde{q}\|_{L^2} \le C\|\tilde{u}\|_{H^{m+3}} \|\l^{m+2 +\fr{\alpha}{2}} \tilde{q}\|_{L^2}
\ee holds for any field $\tilde{u} \in H^{m+3}$ and any scalar $\tilde{q} \in \mathcal{D}(\l^{m+2+\fr{\alpha}{2}})$ with the property $\tilde{u} \cdot \na \tilde{q} \in \mathcal{D}(\l^{m+2})$. Indeed, the commutator in \eqref{pr33} can be written as 
\beg{align}
&(-\Delta)^{\fr{m+2}{2}} (\tilde{u} \cdot \na \tilde{q}) - \tilde{u} \cdot \na (-\Delta)^{\fr{m+2}{2}} \tilde{q} \nonumber
\\&= (-\Delta)^{\fr{m}{2}} \left((-\Delta)\tilde{u} \cdot \na \tilde{q} + \tilde{u} \cdot \na (-\Delta) \tilde{q} - 2 \na \tilde{u} \cdot \na \na \tilde{q} \right)- \tilde{u} \cdot \na (-\Delta)^{\fr{m+2}{2}} \tilde{q} \nonumber
\\&= \left[(-\Delta)^{\fr{m}{2}} (-\Delta \tilde{u} \cdot \na \tilde{q}) - (-\Delta \tilde{u}) \cdot \na (-\Delta)^{\fr{m}{2}} \tilde{q}\right] \nonumber
\\&\quad+ \left[(-\Delta)^{\fr{m}{2}} \left(\tilde{u} \cdot \na (-\Delta) \tilde{q} \right) - \tilde{u} \cdot \na (-\Delta)^{\fr{m}{2}} (-\Delta)\tilde{q} \right]
- 2 (-\Delta)^{\fr{m}{2}} (\na \tilde{u} \cdot \na \na \tilde{q}) - \Delta \tilde{u} \cdot \na (-\Delta)^{\fr{m}{2}} \tilde{q},
\end{align} and its $L^2$ norm is thus bounded as
\beg{align} \la{pr34}
&\|(-\Delta)^{\fr{m+2}{2}} (\tilde{u} \cdot \na \tilde{q}) - \tilde{u} \cdot \na (-\Delta)^{\fr{m+2}{2}} \tilde{q}\|_{L^2} \nonumber
\\&\le C\| \Delta \tilde{u}\|_{H^{m+1}} \|\l^{m + \fr{\alpha}{2}} \tilde{q}\|_{L^2}
+ C\|\tilde{u}\|_{H^{m+1}} \|\l^{m+ \fr{\alpha}{2}} \Delta \tilde{q} \|_{L^2} \nonumber
\\&\quad+ C\|(-\Delta)^{\fr{m}{2}} (\na \tilde{u} \cdot \na \na \tilde{q})\|_{L^2}
+ C\|\Delta \tilde{u} \cdot \na (-\Delta)^{\fr{m}{2}} \tilde{q} \|_{L^2}
\end{align} in view of the induction hypothesis. Since $H^{m}$ is a Banach Algebra and $\mathcal{D}(\l^{m+2})$ is continuously embedded in $H^{m+2}$, we have 
\beg{align} \la{pr35}
&\|(-\Delta)^{\fr{m}{2}} (\na \tilde{u} \cdot \na \na \tilde{q})\|_{L^2}
\le C\|\na \tilde{u}\|_{H^{m}} \|\na \na \tilde{q}\|_{H^{m}} \nonumber
\\&\le C\|\tilde{u}\|_{H^{m+1}} \|\tilde{q}\|_{H^{m+2}}  
\le C\|\tilde{u}\|_{H^{m+1}} \|\l^{m+2} \tilde{q}\|_{L^2} 
\end{align} Putting \eqref{pr34} and \eqref{pr35} together, we obtain \eqref{pr33}. Consequently, \eqref{pr31} holds for all even integers $m$.  The proof of \eqref{pr32} is similar. We omit further details.

\beg{rem} The commutator estimates \eqref{pr31} and \eqref{pr32} are not sharp. Indeed, for any integer $m 
\ge 1$ and $p_1, p_2, q_1, q_2 \in (1, \infty)$ obeying $\fr{1}{p_1} + \fr{1}{q_1} = \fr{1}{p_2} + \fr{1}{q_2} = \fr{1}{2}$, we can show that
\be 
\|[\l^{m}, \tilde{u} \cdot \na] \tilde{q}\|_{L^2} \le C\left[\|\tilde{u}\|_{W^{m, q_1}} \|\tilde{q}\|_{W^{m-1, p_1}} + \|\tilde{u}\|_{W^{m-1, q_2}} \|\tilde{q}\|_{W^{m, p_2}} \right],
\ee holds when $m$ is even, and
\be 
\|[\l^{m-1}\na, \tilde{u} \cdot \na] \tilde{q}\|_{L^2} \le C\left[\|\tilde{u}\|_{W^{m, q_1}} \|\tilde{q}\|_{W^{m-1, p_1}} + \|\tilde{u}\|_{W^{m-1, q_2}} \|\tilde{q}\|_{W^{m, p_2}} \right]
\ee holds when $m$ is odd, by following the induction argument provided above (see for instance \eqref{Bcom} for the base step). The estimate \eqref{pr31} and \eqref{pr32} are adapted to our electroconvection model and are used to prove the $C^{\infty}$ smoothness of solutions. 
\end{rem}

\beg{prop} \la{pr4} Let $m \ge 1$ be an integer. Suppose $v \in \mathcal{D} (A^{\fr{m+1}{2}})$, $\rho \in \mathcal{D}(
\l^m)$, and $F \in \left(H^m\right)^2$. Then it holds that 
\be \la{ellip}
\|A^{\fr{m}{2}} B(v,v)\|_{L^2} \le C\|v\|_{L^{\infty}} \|A^{\fr{m+1}{2}} v\|_{L^2},
\ee and
\be \la{eellip}
\|A^{\fr{m}{2}} \PP(\rho F)\|_{L^2} \le C\left[\|F\|_{L^{\infty}} \|\l^{m} \rho\|_{L^2} + \|F\|_{H^m} \|\rho\|_{L^{\infty}}\right].
\ee
\end{prop}

\noindent \textbf{Proof.} Without loss of generality, assume $m$ is even. By the Helmholtz-Hodge decomposition theorem, the Leray projection of $v \cdot \na v$ can be uniquely decomposed as 
\be 
\PP(v \cdot \na v) = v \cdot \na v + \na \pi
\ee where $\pi$ solves the Poisson equation
\be \la{poisson}
-\Delta \pi = \na \cdot (v \cdot \na v)
\ee with Neumann homogeneous boundary conditions $\fr{\pa \pi}{\pa n} = 0$. Having this decomposition in hand, we bound
\be 
\|A^{\fr{m}{2}} B(v,v)\|_{L^2}
\le C\|\PP(v \cdot \na v)\|_{H^m}
\le C\|v \cdot \na v\|_{H^m} + C\| \na \pi\|_{H^m}
\ee  where we used the estimate $\|A^{\fr{\beta}{2}} \tilde{v}\|_{L^2} \le C\|\tilde{v}\|_{H^{\beta}}$ that holds for any $\tilde{v} \in \mathcal{D}(A^{\fr{\beta}{2}})$ and any $\beta \in \R$ (see \cite{GS}). Since $v$ is divergence-free, we have 
\be 
\|v \cdot \na v\|_{H^m}
= \| \na \cdot (v \otimes  v)\|_{H^m}
\le C\|v \otimes v\|_{H^{m+1}}
\le C\|v\|_{L^{\infty}}\| v\|_{H^{m+1}}
\ee  where the last inequality follows from standard integer Sobolev product estimates. Moreover, the elliptic regularity of the solution to the Poisson equation \eqref{poisson} yields the estimate
\beg{align}
&\|\na \pi\|_{H^m}
\le C\|v \cdot \na v\|_{H^m}
 \le C\|v \otimes v\|_{H^{m+1}}
\le C\|v\|_{L^{\infty}}\| v\|_{H^{m+1}}.
\end{align} In view of the bound $\| v\|_{H^{m+1}} \le C\|A^{\fr{m+1}{2}} v\|_{L^2} $ (see \cite{GS}), we obtain \eqref{ellip}. The proof of \eqref{eellip} is similar and will be omitted.

\beg{prop} \la{prodfor} Let $\delta>0$ be sufficiently small. Let $s \in (\delta, \fr{1}{2} + \delta)$. Let $\rho \in L^{\infty} \cap \mathcal{D}(\l^s)$. Then there exists a constant $C>0$ depending on $s$ and $\delta$ such that the following fractional product inequality 
\be \la{prodfor7}
\|A^{\fr{s-\delta}{2}} \PP (\rho R\rho)\|_{L^2} \le C\|\rho\|_{L^{\infty}} \|\l^{s} \rho\|_{L^2}
\ee holds. 
\end{prop}

\noindent \textbf{Proof.} In view of the bound $\|A^{\fr{\beta}{2}} v\|_{L^2} \le C\|v\|_{H^{\beta}}$ that holds for any $v \in \mathcal{D}(A^{\fr{\beta}{2}})$ and any $\beta \in \R$ (see \cite{GS}), the boundedness of the Leray projector on fractional Sobolev spaces, and the continuous embedding of $\mathcal{D}(\l^{s-\delta})$ into $H^{s-\delta}$ (\cite[Proposition 2.1]{CN3}), we have 
\be 
\|A^{\fr{s-\delta}{2}} \PP (\rho R\rho)\|_{L^2}
\le C\|\PP (\rho R\rho)\|_{H^{s-\delta}}
\le C\|\rho R\rho\|_{H^{s-\delta}}
\le C\|\l^{s-\delta} (\rho R\rho)\|_{L^2}.
\ee We write the Riesz tranform as $R = (R_1, R_2)$, fix $i \in \left\{1,2\right\}$, and use the integral representation \eqref{prodfor1} to estimate 
\beg{align} \la{uyt}
&\|\l^{s-\delta} (\rho R_i\rho)\|_{L^2}^2
= (\l^{s-\delta}(\rho R_i\rho), \l^{s- \delta}(\rho R_i\rho))_{L^2} \nonumber
\\&\quad= \int_{\Omega} \int_{\Omega} (\rho(x)R_i\rho(x) - \rho(y)R_i\rho(y))^2 K_{s-\delta}(x,y) dxdy
+ \int_{\Omega} \rho(x)^2R_i\rho(x)^2 B_{s-\delta}(x) dx \nonumber
\\&\quad\le 2\int_{\Omega} \int_{\Omega} (\rho(x) - \rho(y))^2 R_i\rho(x)^2 K_{s-\delta}(x,y) dxdy \nonumber
\\&\quad\quad\quad+ 2\int_{\Omega} \int_{\Omega} \rho(y)^2 (R_i\rho(x) - R_i\rho(y))^2 K_{s-\delta}(x,y) dxdy
+ \int_{\Omega} \rho(x)^2 R_i\rho(x)^2 B_{s-\delta}(x)dx
\end{align} where the last bound is obtained by adding and subtracting $\rho(y)R_i\rho(x)$ followed by an application of the algebraic inequality $(a+b)^2 \le 2(a^2 + b^2)$. Due to \eqref{prodfor2}. the kernel $K_{s-\delta}$ is bounded from above by a constant multiple of $|x-y|^{-(2+2s-2\delta)}$. Thus, the estimate \eqref{uyt} boils down to
\beg{align} \la{prodfor5}
\|\l^{s-\delta} (\rho R_i\rho)\|_{L^2}^2 
&\le C\int_{\Omega} R_i\rho(x)^2 \int_{\Omega} \fr{|\rho(x)- \rho(y)|^2}{|x-y|^{2 +2s - 2\delta}} dxdy \nonumber
\\&\quad\quad\quad+ C\|\rho\|_{L^{\infty}}^2 \left[\int_{\Omega} \int_{\Omega} (R_i\rho(x) - R_i\rho(y))^2 K_{s-\delta}(x,y) dxdy + \int_{\Omega} R_i\rho(x)^2 B_{s-\delta}(x)dx \right] \nonumber
\\&\quad\le  C\int_{\Omega} R_i\rho(x)^2 \int_{\Omega} \fr{|\rho(x)- \rho(y)|^2}{|x-y|^{2 +2s - 2\delta}} dxdy
+ C\|\rho\|_{L^{\infty}}^2 \|\l^{s-\delta} R_i \rho\|_{L^2}^2.
\end{align} Since $s - \delta \in (0, \fr{1}{2})$, the space $\mathcal{D}(\l^{s-\delta})$ is identified with the usual Sobolev space $H^{s-\delta}$. Thus, the Riesz transform of $\rho$ belongs to $\mathcal{D}(\l^{s-\delta})$ and satisfies the estimate
\be 
\|\l^{s-\delta} R \rho\|_{L^2} 
\le C\|R \rho\|_{H^{s-\delta}} 
= C\|\na \l^{-1} \rho\|_{H^{s-\delta}}
\le C\|\l^{-1} \rho\|_{H^{s-\delta+1}}
\le C\|\l^{s-\delta} \rho\|_{L^2}
\ee where the last inequality follows from the continuous embedding of $\mathcal{D}(\l^{s-\delta+1})$ in $H^{s-\delta+1}$ (\cite[Proposition 2.1]{CN3}). Now we seek good control of the double integral in \eqref{prodfor5}. In fact, an application of H\"older inequality yields
\beg{align} 
\int_{\Omega} R_i\rho(x)^2 \int_{\Omega} \fr{|\rho(x)- \rho(y)|^2}{|x-y|^{2 +2s - 2\delta}} dxdy
&\le \left(\int_{\Omega} \int_{\Omega} R_i\rho(x)^{2p_1} dxdy\right)^{\fr{1}{p_1}} \left(\int_{\Omega} \int_{\Omega} \fr{|\rho(x)- \rho(y)|^{2p_2}}{|x-y|^{(2 +2s - 2\delta)p_2}} dxdy\right)^{\fr{1}{p_2}} \nonumber
\\&\quad\le C_{\Omega} \|R_i\rho\|_{L^{2p_1}}^2  \left(\int_{\Omega} \int_{\Omega} \fr{|\rho(x)- \rho(y)|^{2p_2}}{|x-y|^{2 + 2p_2(1 + s - \delta - 1/p_2)}}\right)^{\fr{1}{p_2}}
\end{align} for any $p_1, p_2 \ge 1$ obeying $\fr{1}{p_1} + \fr{1}{p_2} = 1$. We choose $p_2$ slightly bigger than 1 so that $\fr{1}{p_2} = 1 - \fr{\delta}{2}$. Due to the finiteness of the domain size, it holds that
\be 
\|R_i \rho\|_{L^{2p_1}} \le C\|\rho\|_{L^{2p_1}} \le C_{\Omega} \|\rho\|_{L^{\infty}}.
\ee
Here we used the boundedness of the Dirichlet Riesz transform $R: L^p(\Omega) \rightarrow L^p(\Omega)$ on bounded domains with smooth boundaries for any $p \in (1, \infty)$ (see \cite{JK,DM,S}), a fact that was obtained first in \cite{JK} for bounded Lipschitz domains (with restrictions on the values of $p$) and $C^1$ domains (for any $p \in (1, \infty)$) based on complex interpolation techniques, and later in \cite{S} based on a classical Calderon- Zygmund decomposition approach.
Consequently, we obtain 
\beg{align}
\int_{\Omega} R_i\rho(x)^2 \int_{\Omega} \fr{|\rho(x)- \rho(y)|^2}{|x-y|^{2 +2s - 2\delta}} dxdy
&\le C_{\Omega} \|\rho\|_{L^{\infty}}^2 \left(\int_{\Omega} \int_{\Omega} \fr{|\rho(x)- \rho(y)|^{\fr{4}{2-\delta}}}{|x-y|^{2 + \left(s - \fr{\delta}{2} \right)\left(\fr{4}{2-\delta} \right)}}\right)^{1-\fr{\delta}{2}} \nonumber
\\&= C_{\Omega} \|\rho\|_{L^{\infty}}^2 \|\rho\|_{W^{s - \fr{\delta}{2},\fr{4}{2-\delta}}}^2.
\end{align} Due to the continuous embeddings of $H^{s}$ into $W^{s - \fr{\delta}{2}, \fr{4}{2-\delta}}$ for sufficiently small $\delta$ and  $\mathcal{D}(\l^s)$ into $H^s$,  we infer that 
\be \la{prodfor6}
\int_{\Omega} R_i\rho(x)^2 \int_{\Omega} \fr{|\rho(x)- \rho(y)|^2}{|x-y|^{2 +2s - 2\delta}} dxdy
\le C \|\rho\|_{L^{\infty}}^2 \|\l^s \rho\|_{L^2}^2.
\ee Putting together \eqref{prodfor5} and \eqref{prodfor6} gives the desired product estimate \eqref{prodfor7}.

\section{Existence and Uniqueness of Solutions: Proof of Theorem \ref{Existence}} \la{s4}

For $\alpha \in (0,1]$, we consider the system of equations
\be \la{eq1}
\pa_t q + u \cdot \na q + \l^{\alpha} q = 0,
\ee
\be \la{eq11}
\pa_t u + u \cdot \na u + \na p - \Delta u = - qRq,
\ee
\be \la{eq1t}
\na \cdot u = 0,
\ee on $\Omega \times [0, \infty)$, with homogeneous Dirichlet boundary conditions 
\be \la{eq2}
q|_{\pa \Omega} = u|_{\pa \Omega} = 0,
\ee and with initial data $q_0$ and $u_0$. 

In this section, we address the existence, uniqueness, and long time behavior of solutions to the system described by the equations \eqref{eq1}--\eqref{eq2}. For this objective, we consider the $\epsilon$-approximate system 
\be \la{eq4}
\pa_t q^{\epsilon} + u^{\epsilon} \cdot \na q^{\epsilon} + \l^{\alpha}q^{\epsilon} - \epsilon \Delta q^{\epsilon} = 0,
\ee
\be \la{eq41}
\pa_t u^{\epsilon} + u^{\epsilon} \cdot \na u^{\epsilon} + \na p^{\epsilon} - \Delta u^{\epsilon}
= -q^{\epsilon} \na (\l^{-1})_{\epsilon} q^{\epsilon},
\ee 
\be 
\na \cdot u^{\epsilon} =0
\ee where 
\be 
(\l^{-1})_{\epsilon} \rho = \int_{\epsilon}^{\infty} t^{-\fr{1}{2}} e^{t\Delta} \rho dt,
\ee  with homegeneous Dirichlet boundary conditions 
\be \la{eq6}
q^{\epsilon}|_{\pa \Omega} = u^{\epsilon}|_{\pa \Omega} = 0,
\ee and with initial data $q^{\epsilon}(0) = q_0$ and $u^{\epsilon}(0) = u_0$. For each fixed $\epsilon \in (0,1)$, we prove that the approximate system \eqref{eq4}--\eqref{eq6} has unique global smooth solutions. We need first the following lemma:

\beg{lem} \la{le1} Let $\epsilon > 0$ be fixed. Let $f \in \mathcal{D}(\l^{s-1})$. Then there exist a positive universal constant $C$, independent of $\epsilon$, and a positive constant $C_{s}$ depending only on $s$ and universal constants, such that the inequalities
\be \la{le11}
\|\l^{s} (\l^{-1})_{\epsilon} f\|_{L^2}
\le C\|\l^{s-1} f\|_{L^2}
\ee and 
\be \la{le12}
\|\l^{s} (\l^{-1})_{\epsilon} f\|_{L^2} \le C_s \epsilon^{-s+1} \|f\|_{L^2}
\ee hold for any $s \ge 0$. 
\end{lem}

\noindent \textbf{Proof.} The proof of \eqref{le12} can be found in \cite{I}. We prove the inequality \eqref{le11}. We consider the expansion of $f$
\be 
f = \sum\limits_{j=1}^{\infty} (f, w_j)_{L^2}w_j
\ee in terms of the eigenfunctions $w_j$ of the Laplace operator $-\Delta$. We write the expansion of $(\l^{-1})_{\epsilon}f$ as follows,
\be 
(\l^{-1})_{\epsilon} f = \sum\limits_{j=1}^{\infty} \psi_j w_j,
\ee where the coefficients $\psi_j$ are given by the integral representation
\be 
\psi_j = \int_{\epsilon}^{\infty} t^{-\fr{1}{2}} e^{-t\lambda_j} (f, w_j)_{L^2} dt.
\ee By making the change of variable $t\lambda_j = s$, we have 
\be 
|\psi_j| \le \left(\int_{\epsilon \lambda_j}^{\infty} s^{-\fr{1}{2}} e^{-s} ds \right) \lambda_j^{-\fr{1}{2}} |(f, w_j)_{L^2}|
\le \left(\int_{0}^{\infty} s^{-\fr{1}{2}} e^{-s} ds \right) \lambda_j^{-\fr{1}{2}} |(f, w_j)_{L^2}|
\le  C\lambda_j^{-\fr{1}{2}} |(f, w_j)_{L^2}|
\ee for all $j \ge 1$. Therefore, we have 
\be 
\|\l^{s} (\l^{-1})_{\epsilon}f\|_{L^2}^2 
= \sum\limits_{j=1}^{\infty} \lambda_j^s |\psi_j|^2
\le C\sum\limits_{j=1}^{\infty} \lambda_j^{s-1} (f, w_j)_{L^2}^2
= C\|\l^{s-1}f\|_{L^2}^2,
\ee completing the proof of the Lemma \ref{le1}.

\beg{prop} \la{epsapp} Fix an $\epsilon \in (0,1)$ and an arbitrary time $T>0$. Suppose $u_0 \in \mathcal{D}(A^{\fr{1}{2}})$ and $q_0 \in \mathcal{D}(\l)$. Then the $\epsilon$-approximate system \eqref{eq4}--\eqref{eq6} has a unique solution $(q^{\epsilon}, u^{\epsilon})$ on $[0, T]$ with regularity 
\be \la{prop21}
q^{\epsilon} \in L^{\infty}(0,T; \mathcal{D}(\l)) \cap L^2(0,T; \mathcal{D}(\l^2))
\ee and 
\be \la{prop22}
u^{\epsilon} \in \left(L^{\infty} (0, T; \mathcal{D}(A^{\fr{1}{2}})) \cap L^2(0,T; \mathcal{D}(A))\right)^2
\ee 
\end{prop}

\noindent \textbf{Proof.} For $n \ge 1$, we consider the Galerkin approximants 
\be 
\PP_n \rho = \sum\limits_{j=1}^{n} (\rho, w_j)_{L^2} w_j
\ee and
\be 
\PP_n v = \sum\limits_{j=1}^{n} (v, \phi_j)_{L^2} \phi_j
\ee where $w_j$ and $\phi_j$ are the the eigenfunctions of the homogeneous Dirichlet Laplace operator $-\Delta$ and the Stokes operator $A$ respectively. Here, we abused notation and wrote $\PP_n$ for both projections.

For a fixed $\epsilon > 0$ and $n \ge 1$, we consider the approximate system of ODEs
\be \la{Gal1}
\pa_t q_n^{\epsilon} + \PP_n (u_n^{\epsilon} \cdot \na q_n^{\epsilon}) + \l^{\alpha}q_n^{\epsilon} - \epsilon \Delta q_n^{\epsilon} =  0,
\ee
\be \la{Gal2}
\pa_t u_n^{\epsilon} + Au_n^{\epsilon} + \PP_n (B(u_n^{\epsilon}, u_n^{\epsilon})) = -\PP_n (q_n^{\epsilon} \na (\l^{-1})_{\epsilon} q_n^{\epsilon})
\ee with initial data $q_n^{\epsilon}(0) = \PP_n q_0$ and $u_{n}^{\epsilon} (0) = \PP_n u_0$, and homogeneous Dirichlet boundary conditions $q_n^{\epsilon}|_{\pa \Omega} = u_n^{\epsilon}|_{\pa \Omega} = 0$.  We establish a priori uniform-in-$n$ bounds as follows:

\textbf{Step 1. $L^2$ bounds for the charge density approximants.} We take the scalar product in $L^2$ of the equation \eqref{Gal1} obeyed by the charge density approximants $q_n^{\epsilon}$ with $q_n^{\epsilon}$. We obtain the energy equality 
\be 
\fr{1}{2} \fr{d}{dt} \|q_n^{\epsilon}\|_{L^2}^2 + \|\l^{\fr{\alpha}{2}}q_n^{\epsilon}\|_{L^2}^2 + \epsilon \|\l q_n^{\epsilon}\|_{L^2}^2 = 0, 
\ee where the nonlinear term vanishes due divergence-free condition obeyed by $u_n^{\epsilon}$. Integrating in time from $0$ to $t$, and taking the supremum over $[0,T]$, we conclude that 
\be \la{s1}
\sup\limits_{0 \le t \le T} \|q_n^{\epsilon}(t)\|_{L^2}^2 
+ 2\int_{0}^{T} \left(\|\l^{\fr{\alpha}{2}}q_n^{\epsilon}(t) \|_{L^2}^2 +  \epsilon \|\l q_n^{\epsilon}(t)\|_{L^2}^2  \right) dt
\le 2\|q_0\|_{L^2}^2. 
\ee

\textbf{Step 2. $L^2$ bounds for the velocity approximants.} We take the $L^2$ inner product of the equation \eqref{Gal2} obeyed by the velocity approximants $u_n^{\epsilon}$ with $u_n^{\epsilon}$. The nonlinear term $(\PP_n(B(u_n^{\epsilon}, u_n^{\epsilon})), u_{n}^{\epsilon})_{L^2}$ vanishes due to the self-adjointness of the Leray projector $\PP$ and the divergence-free condition obeyed by $u_{n}^{\epsilon}$. We obtain the energy equation
\be 
\fr{1}{2} \fr{d}{dt} \|u_n^{\epsilon}\|_{L^2}^2 + \| A^{\fr{1}{2}}u_n^{\epsilon}\|_{L^2}^2
= -\int_{\Omega} \PP_n (q_n^{\epsilon} \na (\l^{-1})_{\epsilon} q_n^{\epsilon}) \cdot u_{n}^{\epsilon} dx. 
\ee
Using the fact that $u_n^{\epsilon}$ belongs to the space spanned by the first $n$ eigenfunctions of $A$, choosing $p \in (2,4]$ such that the continuous Sobolev embedding $ \mathcal{D}(\l^{\fr{\alpha}{2}}) \subset L^p$ holds, and using the boundedness of the Riesz transform $R = \na \l^{-1}$ on $L^p$, we estimate 
\beg{align}
&\left|\int_{\Omega} \PP_n (q_n^{\epsilon} \na (\l^{-1})_{\epsilon} q_n^{\epsilon}) \cdot u_{n}^{\epsilon} dx \right|
= \left|\int_{\Omega} q_n^{\epsilon} \na (\l^{-1})_{\epsilon} q_n^{\epsilon} \cdot u_{n}^{\epsilon} dx \right| \nonumber
\\&\le \|q_n^{\epsilon}\|_{L^2} \| \na \l^{-1} \l (\l^{-1})_{\epsilon} q_n^{\epsilon}\|_{L^p} \|u_{n}^{\epsilon}\|_{L^q}
= \|q_n^{\epsilon}\|_{L^2} \| R \l (\l^{-1})_{\epsilon} q_n^{\epsilon}\|_{L^p} \|u_{n}^{\epsilon}\|_{L^q}
\nonumber
\\&\le C\|q_n^{\epsilon}\|_{L^2} \|\l (\l^{-1})_{\epsilon} q_n^{\epsilon}\|_{L^p} \|u_{n}^{\epsilon}\|_{L^q}
\le C\|q_n^{\epsilon}\|_{L^2} \|\l^{\fr{\alpha}{2} + 1}  (\l^{-1})_{\epsilon} q_n^{\epsilon}\|_{L^2} \|u_{n}^{\epsilon}\|_{L^q}
\end{align} where $q$ is the H\"older exponent satisfying $\fr{1}{q} + \fr{1}{p} + \fr{1}{2} = 1.$ In view of Lemma \ref{le1} and the Poincar\'e inequality 
\be 
\|u_{n}^{\epsilon}\|_{L^q} \le C\|\na u_n^{\epsilon}\|_{L^2} = C\|A^{\fr{1}{2}} u_n^{\epsilon}\|_{L^2},
\ee we infer that 
\be \la{s22}
\left|\int_{\Omega} \PP_n (q_n^{\epsilon} \na (\l^{-1})_{\epsilon} q_n^{\epsilon}) \cdot u_{n}^{\epsilon} dx \right|
\le C\|q_n^{\epsilon}\|_{L^2}\|\l^{\fr{\alpha}{2}} q_n^{\epsilon}\|_{L^2} \|A^{\fr{1}{2}} u_n^{\epsilon}\|_{L^2},
\ee from which we obtain the differential inequality
\be \la{s21}
\fr{d}{dt} \|u_n^{\epsilon}\|_{L^2}^2 + \| A^{\fr{1}{2}}u_n^{\epsilon}\|_{L^2}^2 
\le C\|q_n^{\epsilon}\|_{L^2}^2 \|\l^{\fr{\alpha}{2}}q_n^{\epsilon}\|_{L^2}^2
\ee after use of Young's inequality for products. Now we integrate in time from $0$ to $t$, take the supremum over $[0,T]$, and use the uniform bound \eqref{s1} for the charge density approximants $q_n^{\epsilon}$ derived in Step 1 to conclude that
\be \la{s2}
\sup\limits_{0 \le t \le T} \|u_n^{\epsilon}(t) \|_{L^2}^2
+ \int_{0}^{T} \|\na u_n^{\epsilon}(t)\|_{L^2}^2 dt
\le 2\|u_0\|_{L^2}^2 + C\|q_0\|_{L^2}^4.
\ee 

\textbf{Step 3. $H^1$ bounds for the velocity approximants.} We take the scalar product in $L^2$ of the equation \eqref{Gal2} obeyed by $u_n^{\epsilon}$ with $Au_n^{\epsilon}$. We estimate the convective nonlinear term
\be \la{s31}
\left|\int_{\Omega} \PP_n B(u_n^{\epsilon}, u_n^{\epsilon}) Au_n^{\epsilon} dx \right|
\le C\|u_n^{\epsilon}\|_{L^4} \|\na u_n^{\epsilon}\|_{L^4} \|Au_n^{\epsilon}\|_{L^2}
\le C\|u_n^{\epsilon}\|_{L^2}^{\fr{1}{2}} \|\na u_n^{\epsilon}\|_{L^2} \|Au_n^{\epsilon}\|_{L^2}^{\fr{3}{2}}
\ee 
via use of the Ladyzhenskaya interpolation inequality, and the ellipticity of the Stokes operator. As for the electrical forcing nonlinear term, we choose $p \in (2,4]$ so that $\mathcal{D}(\l^{\fr{\alpha}{2}})$ is continuously embedded in $L^p$, apply H\"older's inequality with exponents $p$, $q = \fr{2p}{p-2}$, $2$, use the boundedness of the Riesz transform on $L^q$ and the continuous embedding of $\mathcal{D}(\l^{\fr{2}{p}})$ in $L^q$, and apply Lemma \ref{le1} to obtain 
\beg{align}
&\left|\int_{\Omega} \PP_n (q_n^{\epsilon} \na (\l^{-1})_{\epsilon}q_n^{\epsilon}) Au_n^{\epsilon} dx \right|
\le C\|q_n^{\epsilon}\|_{L^p} \|\na (\l^{-1})_{\epsilon}q_n^{\epsilon}\|_{L^q} \|Au_n^{\epsilon}\|_{L^2} \nonumber
\\&\le C\|\l^{\fr{\alpha}{2}}q_n^{\epsilon}\|_{L^2}\|R \l  (\l^{-1})_{\epsilon}q_n^{\epsilon}\|_{L^q} \|Au_n^{\epsilon}\|_{L^2} 
\le  C\|\l^{\fr{\alpha}{2}}q_n^{\epsilon}\|_{L^2}\|\l (\l^{-1})_{\epsilon}q_n^{\epsilon}\|_{L^q} \|Au_n^{\epsilon}\|_{L^2} \nonumber
\\&\le C\|\l^{\fr{\alpha}{2}}q_n^{\epsilon}\|_{L^2}\|\l^{1 + \fr{2}{p}} (\l^{-1})_{\epsilon}q_n^{\epsilon}\|_{L^2} \|Au_n^{\epsilon}\|_{L^2}
\le C_{\epsilon} \|\l^{\fr{\alpha}{2}}q_n^{\epsilon}\|_{L^2}\|q_n^{\epsilon}\|_{L^2} \|Au_n^{\epsilon}\|_{L^2}
\end{align} where $C_{\epsilon}$ is a positive constant, which does not depend on $n$ but on $\epsilon$, and that blows up as $\epsilon$ approaches $0$. Applying Young's inequality, we obtain the differential inequality
\be 
\fr{d}{dt} \|\na u_n^{\epsilon}\|_{L^2}^2
+ \|Au_n^{\epsilon}\|_{L^2}^2
\le C\|u_n^{\epsilon}\|_{L^2}^2 \|\na u_n^{\epsilon}\|_{L^2}^4
+ C_{\epsilon} \|q_n^{\epsilon}\|_{L^2}^2 \|\l^{\fr{\alpha}{2}} q_n^{\epsilon}\|_{L^2}^2.
\ee In view of Gronwall's inequality and the bounds \eqref{s1} and \eqref{s2} derived in Steps 1 and 2, we conclude that 
\be \la{s3}
\sup\limits_{0 \le t \le T} \|\na u_n^{\epsilon}(t)\|_{L^2}^2 + \int_{0}^{T} \|Au_n^{\epsilon}(t)\|_{L^2}^2 dt
\le Ce^{C(\|u_0\|_{L^2}^4 + \|q_0\|_{L^2}^8)} \left(\|\na u_0\|_{L^2}^2 + C_{\epsilon}\|q_0\|_{L^2}^4 \right).
\ee

\textbf{Step 4. $H^1$ bounds for the charge density approximants.} The $L^2$ norm of $\na q_n^{\epsilon}$ evolves according to the energy equality
\be \la{s41}
\fr{1}{2} \fr{d}{dt} \|\na q_n^{\epsilon}\|_{L^2}^2
+ \|\l^{1+ \fr{\alpha}{2}}q_n^{\epsilon}\|_{L^2}^2
+ \epsilon \|\l^2 q_n^{\epsilon}\|_{L^2}^2
= \int_{\Omega} u_n^{\epsilon} \cdot \na q_n^{\epsilon} \Delta q_n^{\epsilon} dx.
\ee We integrate by parts the nonlinear term, use the homogeneous Dirichlet boundary conditions and divergence-free property obeyed by the velocity approximants $u_n^{\epsilon}$, choose $p$ so that $\mathcal{D}(\l^{\fr{\alpha}{2}})$ is continuously embedded in $L^p$ and $q$ so that $\fr{1}{q} + \fr{1}{p} + \fr{1}{2} = 1$, apply H\"older's inequality, and obtain 
\beg{align} \la{s42}
&\left|\int_{\Omega} u_n^{\epsilon} \cdot \na q_n^{\epsilon} \Delta q_n^{\epsilon} dx \right|
\le \int_{\Omega} |\na u_n^{\epsilon}||\na q_n^{\epsilon}|^2 dx
\le \|\na u_n^{\epsilon}\|_{L^q} \|\na q_n^{\epsilon}\|_{L^p} \|\na q_n^{\epsilon}\|_{L^2} \nonumber
\\&= \|\na u_n^{\epsilon}\|_{L^q} \|R \l q_n^{\epsilon}\|_{L^p} \|\na q_n^{\epsilon}\|_{L^2}
\le C\|\na u_n^{\epsilon}\|_{L^q} \|\l q_n^{\epsilon}\|_{L^p} \|\na q_n^{\epsilon}\|_{L^2} \nonumber
\\&\le C\|\na u_n^{\epsilon}\|_{L^q} \|\l^{1+ \fr{\alpha}{2}} q_n^{\epsilon}\|_{L^2} \|\na q_n^{\epsilon}\|_{L^2}.
\end{align} In view of the Gagliardo-Nirenberg interpolation inequality, we have
\be 
\|\na u_n^{\epsilon}\|_{L^q} 
\le C\|\na u_n^{\epsilon}\|_{L^2}^{\fr{2}{q}} \|\na u_n^{\epsilon}\|_{H^1}^{\fr{q-2}{q}}.
\ee Integrating by parts, and applying the Poincar\'e inequality to the vector field $u_n^{\epsilon}$ that vanishes on the boundary of $\Omega$, we observe that
\be 
\|\na u_n^{\epsilon}\|_{L^2}^2
= -\int_{\Omega} u_n^{\epsilon} \cdot \Delta u_n^{\epsilon} dx
\le \|u_n^{\epsilon}\|_{L^2} \|\Delta u_n^{\epsilon}\|_{L^2}
\le C\|\na u_n^{\epsilon}\|_{L^2}\|\Delta u_n^{\epsilon}\|_{L^2}
\le \fr{1}{2} \|\na u_n^{\epsilon}\|_{L^2}^2 + C\|\Delta u_n^{\epsilon}\|_{L^2}^2,
\ee from which we infer that 
\be 
\|\na u_n^{\epsilon}\|_{L^2} \le C\|\Delta u_n^{\epsilon}\|_{L^2},
\ee and consequently
\be \la{s43}
\|\na u_n^{\epsilon}\|_{L^q} \le C\|\Delta u_n^{\epsilon}\|_{L^2} \le C\|A u_n^{\epsilon}\|_{L^2}.
\ee Putting \eqref{s41}, \eqref{s42} and \eqref{s43} together, we obtain the differential inequality
\be \la{s44}
\fr{d}{dt} \|\na q_n^{\epsilon}\|_{L^2}^2
+ \|\l^{1+ \fr{\alpha}{2}}q_n^{\epsilon}\|_{L^2}^2
+ 2\epsilon \|\l^2 q_n^{\epsilon}\|_{L^2}^2
\le C\|A u_n^{\epsilon}\|_{L^2}^2 \|\na q_n^{\epsilon}\|_{L^2}^2
\ee after applying Young's inequality. By Gronwall's inequality and the bound \eqref{s3}, we conclude that
\beg{align}
&\sup\limits_{0 \le t \le T} \|\na q_n^{\epsilon}(t)\|_{L^2}^2
+ \int_{0}^{T} \left(\|\l^{1 + \fr{\alpha}{2}} q_n^{\epsilon}(t) \|_{L^2}^2 + \epsilon \|\l^2q_n^{\epsilon}(t)\|_{L^2}^2 \right) \nonumber
\\&\quad\quad \le Ce^{(\|\na u_0\|_{L^2}^2 + C_{\epsilon}\|q_0\|_{L^2}^4)e^{C(\|u_0\|_{L^2}^4 + \|q_0\|_{L^2}^8)}} \|\na q_0\|_{L^2}^2,
\end{align} ending the proof of Step 4.

The existence of a solution $(q^{\epsilon}, u^{\epsilon})$ to the $\epsilon$-approximate system \eqref{eq4}--\eqref{eq6} with regularity \eqref{prop21} and \eqref{prop22} is obtained via application of the Aubin-Lions lemma, use of the uniform in $n$ bounds derived in Steps 1, 2, 3 and 4, and passage in the weak limit with use of the lower semi-continuity of the norms. 

As for uniqueness, suppose $(q_1^{\epsilon}, u_1^{\epsilon})$ and $(q_2^{\epsilon}, u_2^{\epsilon})$ are two solutions to the $\epsilon$-approximate system \eqref{eq4}--\eqref{eq6} with regularity \eqref{prop21} and \eqref{prop22}. We denote by $\tilde{q}^{\epsilon}, \tilde{u}^{\epsilon}$, and $\tilde{p}^{\epsilon}$ the differences
\be 
\tilde{q}^{\epsilon} = q_1^{\epsilon} - q_2^{\epsilon}, \tilde{u}^{\epsilon} = u_1^{\epsilon} - u_2^{\epsilon}, \tilde{p}^{\epsilon} = p_1^{\epsilon} - p_2^{\epsilon},
\ee which obey the system 
\be \la{quni}
\pa_t \tilde{q}^{\epsilon} + \l^{\alpha} \tilde{q}^{\epsilon} - \epsilon \Delta \tilde{q}^{\epsilon} = -u_1^{\epsilon} \cdot \na \tilde{q}^{\epsilon} - \tilde{u}^{\epsilon} \cdot \na q_2^{\epsilon},
\ee
\be \la{uuni}
\pa_t \tilde{u}^{\epsilon} - \Delta \tilde{u}^{\epsilon} + \na \tilde{p}^{\epsilon} 
= -u_1^{\epsilon} \cdot \na \tilde{u}^{\epsilon} - \tilde{u}^{\epsilon} \cdot \na u_2^{\epsilon} - q_1^{\epsilon} \cdot \na (\l^{-1})_{\epsilon} \tilde{q}^{\epsilon} - \tilde{q}^{\epsilon} \cdot \na (\l^{-1})_{\epsilon} q_2^{\epsilon},
\ee with homogeneous Dirichlet boundary conditions and vanishing initial data. 
We take the scalar product in $L^2$ of \eqref{quni} and \eqref{uuni} with $q^{\epsilon}$ and $u^{\epsilon}$ respectively. We estimate using Ladyzhenskaya's interpolation inequality, the continuous embeddings of $\mathcal{D}(A^{\fr{1}{2}})$ into $L^{\fr{4}{\alpha}}$ and  $\mathcal{D}(\l^{\fr{\alpha}{2}})$ into $L^{\fr{4}{2-\alpha}}$. We obtain the differential inequalities
\beg{align}
\fr{1}{2} \fr{d}{dt} \|\tilde{q}^{\epsilon}\|_{L^2}^2 + \|\l^{\fr{\alpha}{2}}\tilde{q}^{\epsilon}\|_{L^2}^2
\le \|\tilde{u}^{\epsilon}\|_{L^{\fr{4}{\alpha}}} \|\na q_2^{\epsilon}\|_{L^{\fr{4}{2-\alpha}}} \|\tilde{q}^{\epsilon}\|_{L^2}
\le C\|\na \tilde{u}^{\epsilon}\|_{L^2} \|\l^{1 + \fr{\alpha}{2}} q_2^{\epsilon} \|_{L^2} \|\tilde{q}^{\epsilon}\|_{L^2}
\end{align} and
\beg{align}
&\fr{1}{2} \fr{d}{dt} \|\tilde{u}^{\epsilon}\|_{L^2}^2 + \|\na \tilde{u}^{\epsilon}\|_{L^2}^2
\le \|\tilde{u}^{\epsilon}\|_{L^4}^2 \|\na u_2^{\epsilon}\|_{L^2}
+ \left(\|q_1^{\epsilon}\|_{L^{4}} \|\na (\l^{-1})_{\epsilon} \tilde{q}^{\epsilon}\|_{L^2} 
+ \|\tilde{q}^{\epsilon}\|_{L^2} \|\na (\l^{-1})_{\epsilon} q_2^{\epsilon}\|_{L^4}\right) \|\tilde{u}^{\epsilon}\|_{L^4} \nonumber
\\&\le C\|\tilde{u}^{\epsilon}\|_{L^2} \|\na \tilde{u}^{\epsilon}\|_{L^2} \|\na u_2^{\epsilon}\|_{L^2}
+ C\left(\|q_1^{\epsilon}\|_{L^{4}} + \|\l^{\fr{1}{2}} q_2^{\epsilon}\|_{L^2}\right) \|\tilde{q}^{\epsilon}\|_{L^2} \|\tilde{u}^{\epsilon}\|_{L^2}^{\fr{1}{2}} \|\na \tilde{u}^{\epsilon}\|_{L^2}^{\fr{1}{2}},
\end{align} which, added together, yield the energy inequality
\be \la{co}
\fr{d}{dt} \left(\|\tilde{q}^{\epsilon} \|_{L^2}^2 + \|\tilde{u}^{\epsilon}\|_{L^2}^2 \right) 
\le C\left(\|\l^{1 + \fr{\alpha}{2}}q_2^{\epsilon}\|_{L^{2}}^2 + \|\na u_2^{\epsilon}\|_{L^2}^2 + \|q_1^{\epsilon}\|_{L^{4}}^2 
+ 1 \right)\left(\|\tilde{q}^{\epsilon} \|_{L^2}^2 + \|\tilde{u}^{\epsilon}\|_{L^2}^2 \right) 
\ee after applications of Young's inequality. By Gronwall's inequality, we infer that
\be 
\|\tilde{q}^{\epsilon}(t) \|_{L^2}^2 + \|\tilde{u}^{\epsilon}(t)\|_{L^2}^2 \le \exp \left(C(t) \right) \left(\|\tilde{q}^{\epsilon}(0) \|_{L^2}^2 + \|\tilde{u}^{\epsilon}(0)\|_{L^2}^2  \right),
\ee where
\be 
C(t) = C\int_{0}^{t} \left(\|\l^{1 + \fr{\alpha}{2}}q_2^{\epsilon}(s)\|_{L^{2}}^2 + \|\na u_2^{\epsilon}(s)\|_{L^2}^2 + \|q_1^{\epsilon}(s)\|_{L^{4}}^2 + 1\right) ds.
\ee is finite. This gives the uniqueness of the solutions to \eqref{eq4}--\eqref{eq6}, completing the proof of Proposition \ref{epsapp}. 

Now we prove Theorem \ref{Existence}: 

\noindent \textbf{Proof of Theorem \ref{Existence}.} The proof is divided into several steps. 

\textbf{Step 1. Uniform $L^2$ bounds for $q^{\epsilon}$.} The $L^2$ norm of $q^{\epsilon}$ evolves according to the energy equality
\be \la{t21}
\fr{1}{2} \fr{d}{dt} \|q^{\epsilon}\|_{L^2}^2 + \|\l^{\fr{\alpha}{2}}q^{\epsilon}\|_{L^2}^2 +\epsilon \|\l q^{\epsilon}\|_{L^2}^2 = 0,
\ee from which we obtain the differential inequality
\be 
\fr{d}{dt} \|q^{\epsilon}\|_{L^2}^2 + c_1\|q^{\epsilon}\|_{L^2}^2 \le 0
\ee in view of the Poincar\'e inequality. Multiplying both sides by $e^{c_1t}$, and integrating in time from $0$ to $t$, we infer that
\be \la{t22}
\|q^{\epsilon}(t)\|_{L^2}^2 \le \|q_0\|_{L^2}^2 e^{-c_1t}
\ee for all $t \ge 0$. Integrating \eqref{t21} in time from $0$ to $t$, we also have the bound
\be \la{t230}
\int_{0}^{t} \|\l^{\fr{\alpha}{2}}q^{\epsilon}(s) \|_{L^2}^2 ds \le \fr{1}{2} \|q_0\|_{L^2}^2
\ee for all $t \ge 0$. 

\textbf{Step 2. Uniform $L^2$ bounds for $u^{\epsilon}$.} We take the $L^2$ inner product of the equation \eqref{eq41} obeyed by $u^{\epsilon}$ with $u^{\epsilon}$. Integrating by parts, the nonlinear term $(u^{\epsilon} \cdot \na u^{\epsilon}, u^{\epsilon})_{L^2}$ and the pressure term $(\na p^{\epsilon}, u^{\epsilon})_{L^2}$ vanish due to the divergence-free property of $u^{\epsilon}$. We estimate the nonlinear term $(q^{\epsilon} \na (\l^{-1})_{\epsilon} q^{\epsilon}, u^{\epsilon})_{L^2}$ as in \eqref{s22}, and we conclude that the time derivative of the $L^2$ norm of $u^{\epsilon}$ satisfies the differential inequality
\be \la{t24}
\fr{d}{dt} \|u^{\epsilon}\|_{L^2}^2 + \|\na u^{\epsilon}\|_{L^2}^2
\le C\|q^{\epsilon}\|_{L^2}^2\|\l^{\fr{\alpha}{2}} q^{\epsilon}\|_{L^2}^2,
\ee as shown in \eqref{s21}. In view of the Poincar\'e inequality $c_2 \|u^{\epsilon}\|_{L^2}^2 \le \|\na u^{\epsilon}\|_{L^2}^2$, we obtain the pointwise in time bound
\be 
\|u^{\epsilon}(t)\|_{L^2}^2 \le e^{-\min \left\{c_1, c_2 \right\}t} \left(\|u_0\|_{L^2}^2 + C\|q_0\|_{L^2}^2 \int_{0}^{t} \|\l^{\fr{\alpha}{2}} q^{\epsilon}(s)\|_{L^2}^2 ds \right),
\ee which yields the decaying in time bound 
\be \la{t250}
\|u^{\epsilon}(t)\|_{L^2}^2 \le e^{-\min \left\{c_1, c_2\right\}t} \left(\|u_0\|_{L^2}^2 + C\|q_0\|_{L^2}^4 \right)
\ee due to the bounds \eqref{t22} and \eqref{t230} derived in Step 1. 
Integrating \eqref{t24} in time, we have 
\be \la{t260}
\int_{0}^{t} \|\na u^{\epsilon}(s)\|_{L^2}^2 ds \le \|u_0\|_{L^2}^2 + C\|q_0\|_{L^2}^4.
\ee for all $t \ge 0$.

\textbf{Step 3. Uniform $L^p$ bounds for $q^{\epsilon}$.} For an even integer $p \ge 4$, we take the scalar product in $L^2$ of the equation obeyed by $q^{\epsilon}$ with $(q^{\epsilon})^{p-1}$. The nonlinear term $(u^{\epsilon} \cdot \na q^{\epsilon}, (q^{\epsilon})^{p-1})_{L^2}$ vanishes. Consequently, we obtain the differential equation 
\be 
\fr{1}{p} \fr{d}{dt} \|q^{\epsilon}\|_{L^p}^p
+ \int_{\Omega} (q^{\epsilon})^{p-1} \l^{{\alpha}} (q^{\epsilon}) dx
+ \epsilon \int_{\Omega} (q^{\epsilon})^{p-1} \l^{{2}} (q^{\epsilon}) dx
= 0
\ee which gives the differential inequality 
\be 
 \fr{d}{dt} \|q^{\epsilon}\|_{L^p}^p + C_{\Omega, \alpha} (p-1) \|q^{\epsilon}\|_{L^p}^p \le 0
\ee in view of the Poincar\'e inequality \eqref{pci} for the fractional Laplacian in $L^p$ and the C\'ordoba-C\'ordoba inequality. Therefore, the $L^p$ norm of $q^{\epsilon}$ decays exponentially in time and obeys 
\be \la{t27}
\|q^{\epsilon}(t)\|_{L^p} \le \|q_0\|_{L^p} e^{-\fr{C_{\Omega, \alpha} (p-1)}{p}t}
\ee for all $t \ge 0$.

\textbf{Step 4. Uniform $H^1$ bounds for $u^{\epsilon}$.} The $L^2$ norm of $A^{\fr{1}{2}} u^{\epsilon}$ evolves according to the energy equality
\be \la{t28}
\fr{1}{2} \fr{d}{dt} \|A^{\fr{1}{2}} u^{\epsilon}\|_{L^2}^2
+ \|Au^{\epsilon}\|_{L^2}^2
= (-q^{\epsilon} \na (\l^{-1})_{\epsilon}q^{\epsilon}, Au^{\epsilon})_{L^2}
- (B(u^{\epsilon}, u^{\epsilon}), Au^{\epsilon})_{L^2}.
\ee We estimate 
\be 
|(B(u^{\epsilon}, u^{\epsilon}), Au^{\epsilon})_{L^2}| \le C\|u^{\epsilon}\|_{L^2}^{\fr{1}{2}} \|\na u^{\epsilon}\|_{L^2} \|Au^{\epsilon}\|_{L^2}^{\fr{3}{2}},
\ee as in \eqref{s31}, via applications of Ladyzhenskaya's inequality. Now we choose $p \in (2,4]$ so that $\mathcal{D}(\l^{\fr{\alpha}{2}})$ is continuously embedded in $L^p$, and we let $q$ be the H\"older exponent obeying $\fr{1}{q} + \fr{1}{p} + \fr{1}{2} = 1$. In view of Lemma~\ref{le1}, we have 
\beg{align} \la{t29}
&|(-q^{\epsilon} \na (\l^{-1})_{\epsilon}q^{\epsilon}, Au^{\epsilon})_{L^2}|
\le \|q^{\epsilon}\|_{L^q} \| \na (\l^{-1})_{\epsilon}q^{\epsilon}\|_{L^p} \|Au^{\epsilon}\|_{L^2} \nonumber
\\&\le C\|q^{\epsilon}\|_{L^q} \| \l^{1+ \fr{\alpha}{2}} (\l^{-1})_{\epsilon}q^{\epsilon}\|_{L^2} \|Au^{\epsilon}\|_{L^2}
\le C\|q^{\epsilon}\|_{L^q} \| \l^{\fr{\alpha}{2}}q^{\epsilon}\|_{L^2} \|Au^{\epsilon}\|_{L^2}
\end{align} Putting \eqref{t28}--\eqref{t29} together, and applying Young's inequality for products, we obtain 
\be \la{t210}
\fr{d}{dt} \|A^{\fr{1}{2}} u^{\epsilon}\|_{L^2}^2
+ \|Au^{\epsilon}\|_{L^2}^2
\le C\|q^{\epsilon}\|_{L^{\tilde{q}}}^2\| \l^{\fr{\alpha}{2}}q^{\epsilon}\|_{L^2}^2
+ C\|u^{\epsilon}\|_{L^2}^{2} \|A^{\fr{1}{2}} u^{\epsilon}\|_{L^2}^4 
\ee where $\tilde{q}$ is the smallest even integer greater than or equal to $q$. Finally, we bound the dissipation from below by $\min\left\{c_2, \fr{2C_{\Omega, \alpha}(\tilde{q}-1)}{\tilde{q}} \right\} \|A^{\fr{1}{2}} u^{\epsilon}\|_{L^2}$ using the Poincar\'e inequality $c_2 \|A^{\fr{1}{2}} u^{\epsilon}\|_{L^2}^2 \le \|Au^{\epsilon}\|_{L^2}^2$, multiply both sides of the resulting differential inequality by the integrating factor 
\be 
e^{\min\left\{c_2, \fr{2C_{\Omega, \alpha}(\tilde{q}-1)}{\tilde{q}} \right\}t - C\int_{0}^{t} \|u^{\epsilon}(s)\|_{L^2}^{2} \|A^{\fr{1}{2}} u^{\epsilon}(s)\|_{L^2}^2 ds},
\ee integrate in time from $0$ to $t$, and use the time decaying estimate \eqref{t27} to conclude that
\be 
\|A^{\fr{1}{2}} u^{\epsilon}(t)\|_{L^2}^2
\le e^{- \min\left\{c_2, \fr{2C_{\Omega, \alpha}(\tilde{q}-1)}{\tilde{q}} \right\}t + C\int_{0}^{t} \|u^{\epsilon}\|_{L^2}^{2} \|A^{\fr{1}{2}} u^{\epsilon}\|_{L^2}^2 ds} \left(\|\na u_0\|_{L^2}^2 + C\int_{0}^{t} \|q_0\|_{L^{\tilde{q}}}^2 \|\l^{\fr{\alpha}{2}}q^{\epsilon}\|_{L^2}^2 ds \right),
\ee which reduces to 
\be \la{t212}
\|A^{\fr{1}{2}} u^{\epsilon}(t)\|_{L^2}^2
\le e^{- \min\left\{c_2, \fr{2C_{\Omega, \alpha}(\tilde{q}-1)}{\tilde{q}} \right\}t} e^{C\left(\|u_0\|_{L^2}^2 + C\|q_0\|_{L^2}^4 \right)^2} \left(\|\na u_0\|_{L^2}^2 + C\|q_0\|_{L^{\tilde{q}}}^2 \|q_0\|_{L^2}^2 \right)
\ee as a consequence of the bounds \eqref{t230}, \eqref{t250} and \eqref{t260}. Moreover, the $L^2$ norm of $Au^{\epsilon}$ is square integrable in time and obeys
\be \la{t211}
\int_{0}^{t} \|Au^{\epsilon}(s)\|_{L^2}^2 ds 
\le  C\left(\|u_0\|_{L^2}^2 + \|q_0\|_{L^2}^4 + 1 \right)^2 e^{C\left(\|u_0\|_{L^2}^2 + C\|q_0\|_{L^2}^4 \right)^2} \left(\|\na u_0\|_{L^2}^2 + C\|q_0\|_{L^{\tilde{q}}}^2 \|q_0\|_{L^2}^2 \right)
\ee for all $t \ge 0$. 

\textbf{Step 5. Uniform $H^1$ bounds for $q^{\epsilon}$.} As shown for the Galerkin approximants in \eqref{s44}, the time derivative of the $L^2$ norm of $\na q^{\epsilon}$ satisfies the differential inequality
\be \la{t214}
\fr{d}{dt} \|\na q^{\epsilon}\|_{L^2}^2 + \|\l^{1 + \fr{\alpha}{2}} q^{\epsilon}\|_{L^2}^2
\le C\|Au^{\epsilon}\|_{L^2}^2 \|\na q^{\epsilon}\|_{L^2}^2.
\ee Here we have implicitly used the cancellation law
\be 
\int_{\Omega} u^{\epsilon} \na \na q^{\epsilon} \na q^{\epsilon} = 0
\ee that holds due to Sobolev $H^2$ regularity of $q^{\epsilon}$.  From \eqref{t214}, we obtain 
\be 
\fr{d}{dt} \|\na q^{\epsilon}\|_{L^2}^2 + \left(c_1 - C\|Au^{\epsilon}\|_{L^2}^2 \right) \|\na q^{\epsilon}\|_{L^2}^2
\le 0
\ee after bounding the dissipation from below using the Poincar\'e inequality. We multiply by the integrating factor $e^{c_1t - \int_{0}^{t} \|Au^{\epsilon}\|_{L^2}^2 ds}$, integrate in time from $0$ to $t$, and infer that 
\be \la{t213}
\|\na q^{\epsilon}\|_{L^2}^2 \le \|\na q_0\|_{L^2}^2 e^{C_0} e^{-c_1t}
\ee for all $t \ge 0$, where $C_0$ is a constant depending only on the initial data and is given explicitly by 
\be 
C_0 =  C\left(\|u_0\|_{L^2}^2 + \|q_0\|_{L^2}^4 + 1 \right)^2 e^{C\left(\|u_0\|_{L^2}^2 + C\|q_0\|_{L^2}^4 \right)^2} \left(\|\na u_0\|_{L^2}^2 + C\|q_0\|_{L^{\tilde{q}}}^2 \|q_0\|_{L^2}^2 \right).
\ee Integrating \eqref{t214} in time from $0$ to $t$, we obtain 
\be \la{uuu}
\int_{0}^{t} \|\l^{1 + \fr{\alpha}{2}} q^{\epsilon}(s)\|_{L^2}^2 ds
\le \|\na q_0\|_{L^2}^2 \left(1 + C_0e^{C_0}\right) 
\ee for all $t \ge 0$.  

\beg{rem} Compared to the existence result obtained in \cite{ceiv} which was proved only for $\alpha=1$, Theorem \ref{Existence} requires less regularity on the initial data ($H^1$), and yields furthermore exponential decay in time.
\end{rem}

\section{Higher Regularity: Proof of Theorem \ref{tt2} } \la{S4}

In this section, we bootstrap the regularity of solutions and show that the charge density $q$ and velocity $u$ satisfying \eqref{eq1}--\eqref{eq2} decay exponentially in time in all Sobolev spaces and, consequently, in all H\"older spaces. 

Theorem \ref{tt2} is a direct consequence of the following proposition:

\beg{prop} \la{tt22} Let $\epsilon > 0$. Fix an integer $k \ge 2$. Suppose that $q_0 \in \mathcal{D}(\l^k)$ and $u_0 \in \mathcal{D}(A^{\fr{k}{2}})$. Assume there is a positive constant $\Gamma_{k}$ depending only on the $H^{k-1}$ norms of the initial data and $k$ such that the bounds 
\be \la{tt21}
\|\l^{{k-1}} q^{\epsilon}(t)\|_{L^2}^2 \le \Gamma_{k} e^{-c_1 t}
\ee and
\be \la{tt22}
\int_{0}^{t} \|A^{\fr{k}{2}} u^{\epsilon}(s)\|_{L^2}^2 ds \le \Gamma_{k} 
\ee hold for any $t \ge 0$. Then there is a positive constant $\Gamma_k'$ depending only on the $H^{k}$ norms of the initial data and $k$ and a constant $c>0$ depending only on the size of $\Omega$ and $\alpha$ such that the bounds 
\be 
\|\l^{k} q^{\epsilon}(t) \|_{L^2}^2 + \|A^{\fr{k}{2}} q^{\epsilon} (t)\|_{L^2}^2 \le \Gamma_k' e^{-ct}
\ee and 
\be 
\int_{0}^{t} \|A^{\fr{k+1}{2}} u^{\epsilon}(s)\|_{L^2}^2 ds \le \Gamma_{k}'
\ee
hold for any $t \ge 0$. 
\end{prop}

\noindent \textbf{Proof.} For any arbitrary positive time $T$, the charge density and velocity approximants $q^{\epsilon}$ and $u^{\epsilon}$ belong to the spaces $L^{\infty}(0, T; \mathcal{D} (\l^{k}))\cap L^2(0, T; \mathcal{D}(\l^{k + 1}))$ and $L^{\infty}(0, T; \mathcal{D}(A^{\fr{k}{2}}))\cap L^2(0, T; \mathcal{D}(A^{\fr{k+1}{2}}))$ respectively, a fact that can be shown by performing energy estimates on the Galerkin regularized system \eqref{Gal1}--\eqref{Gal2} and passage in the weak limit by use of the Banach Alaoglu theorem. We establish decaying in time bounds which do not depend on $\epsilon$.  

We start by showing that the spatial $H^{k+1}$ norm of $u^{\epsilon}$ is finite in time over $[0, \infty)$ and obeys
\be  \la{tttt2}
\int_{0}^{\infty} \|u^{\epsilon}(t)\|_{H^{k+1}}^2 dt \le \tilde{\Gamma}_{k}
\ee for some positive constant $\tilde{\Gamma}_k$ depending only on $k$ and the initial data. Indeed, the $L^2$ norm of $A^{\fr{k}{2}} u^{\epsilon}$ evolves according to the energy equality
\beg{align} \la{high3}
&\fr{1}{2} \fr{d}{dt} \|A^{\fr{k}{2}}u^{\epsilon}\|_{L^2}^2 
+ \|A^{\fr{k}{2} + \fr{1}{2}} u^{\epsilon}\|_{L^2}^2 \nonumber
\\&\quad\quad= - \int_{\Omega} A^{\fr{k}{2} - \fr{1}{2}} B(u^{\epsilon}, u^{\epsilon}) \cdot A^{\fr{k}{2} + \fr{1}{2}} u^{\epsilon} dx
- \int_{\Omega} A^{\fr{k}{2} - \fr{1}{2}} \PP(q^{\epsilon} \na (\l^{-1})_{\epsilon}q^{\epsilon}) \cdot A^{\fr{k}{2} + \fr{1}{2}} u^{\epsilon} dx.
\end{align} We estimate the nonlinear term in $u^{\epsilon}$ as follows,
\beg{align} \la{high1}
&\left|\int_{\Omega} A^{\fr{k}{2} - \fr{1}{2}} B(u^{\epsilon}, u^{\epsilon}) \cdot A^{\fr{k}{2} + \fr{1}{2}} u^{\epsilon} dx\right|
\le \fr{1}{4} \|A^{\fr{k}{2} + \fr{1}{2}} u^{\epsilon}\|_{L^2}^2
+ C\|A^{\fr{k}{2} - \fr{1}{2}} B(u^{\epsilon}, u^{\epsilon})\|_{L^2}^2 \nonumber
\\&\quad\quad\le \fr{1}{4} \|A^{\fr{k}{2} + \fr{1}{2}} u^{\epsilon}\|_{L^2}^2
+ C\|u^{\epsilon}\|_{L^{\infty}}^2 \|A^{\fr{k}{2}} u^{\epsilon}\|_{L^2}^2 
\le \fr{1}{4} \|A^{\fr{k}{2} + \fr{1}{2}} u^{\epsilon}\|_{L^2}^2
+ C\|A^{\fr{k}{2}} u^{\epsilon}\|_{L^2}^4 
\end{align} where the second inequality holds due to the fractional product estimate \eqref{ellip} and the last inequality follows from the continuous embedding of $\mathcal{D}(A^{\fr{k}{2}})$ in $L^{\infty}$ when $k$ is strictly greater than 1. As for the nonlinear term in $q^{\epsilon}$, we have
\beg{align} \la{high2}
&\left|\int_{\Omega} A^{\fr{k}{2} - \fr{1}{2}} \PP(q^{\epsilon} \na (\l^{-1})_{\epsilon}q^{\epsilon}) \cdot A^{\fr{k}{2} + \fr{1}{2}} u^{\epsilon} dx\right|
\le \fr{1}{4} \|A^{\fr{k}{2} + \fr{1}{2}} u^{\epsilon}\|_{L^2}^2 
+ C\|A^{\fr{k}{2} - \fr{1}{2}} \PP(q^{\epsilon} \na (\l^{-1})_{\epsilon}q^{\epsilon}) \|_{L^2}^2 \nonumber
\\&\quad\quad\le \fr{1}{4} \|A^{\fr{k}{2} + \fr{1}{2}} u^{\epsilon}\|_{L^2}^2 
+ C\left(\|\l^{k-1} q^{\epsilon}\|_{L^2}^2 \|\na (\l^{-1})_{\epsilon} q^{\epsilon}\|_{L^{\infty}}^2 + \|q^{\epsilon}\|_{L^{\infty}}^2  \|\na (\l^{-1})_{\epsilon} q^{\epsilon}\|_{H^{k-1}}^2   \right) \nonumber
\\&\quad\quad\le  \fr{1}{4} \|A^{\fr{k}{2} + \fr{1}{2}} u^{\epsilon}\|_{L^2}^2 
+ C\|\l^{1+\fr{\alpha}{2}}q^{\epsilon}\|_{L^2}^2  \|\l^{k-1} q^{\epsilon}\|_{L^2}^2
\end{align} where the second inequality follows from Proposition \ref{pr4}, and the last inequality uses the continuous embedding of $\mathcal{D}(\l^{1 + \fr{\alpha}{2}})$ in $L^{\infty}$ and the uniform boundedness of $\na (\l)^{-1}_{\epsilon}$ in Sobolev spaces (see Lemma \ref{le1}). Collecting the bounds \eqref{high1} and \eqref{high2}, and inserting them in \eqref{high3}, we obtain the differential inequality
\be \la{t23}
\fr{d}{dt} \|A^{\fr{k}{2}} u^{\epsilon}\|_{L^2}^2
+ \|A^{\fr{k}{2} + \fr{1}{2}} u^{\epsilon}\|_{L^2}^2
\le C\|\l^{1 + \fr{\alpha}{2}} q^{\epsilon}\|_{L^2}^2 \|\l^{k-1} q^{\epsilon}\|_{L^2}^2
+ C\|A^{\fr{k}{2}} u^{\epsilon}\|_{L^2}^4, 
\ee which reduces to 
\be 
\fr{d}{dt} \|A^{\fr{k}{2}} u^{\epsilon}\|_{L^2}^2
+ (c_2 - C\|A^{\fr{k}{2}} u^{\epsilon}\|_{L^2}^2) \|A^{\fr{k}{2} } u^{\epsilon}\|_{L^2}^2
\le C\|\l^{1 + \fr{\alpha}{2}} q^{\epsilon}\|_{L^2}^2 \|\l^{k-1} q^{\epsilon}\|_{L^2}^2
\ee by making use of the Poincar\'e inequality. Multiplying both sides by the integrating factor 
\be 
e^{\min \left\{c_1,c_2\right\} t - C\int_{0}^{t} \|A^{\fr{k}{2}} u^{\epsilon} (s)\|_{L^2}^2 ds},
\ee integrating in time from $0$ to $t$, using the hypotheses \eqref{tt21} and \eqref{tt22}, and exploiting the square integrability in time of $\|\l^{1+\fr{\alpha}{2}}q^{\epsilon} \|_{L^2}$ obtained in \eqref{uuu}, we infer that 
\be 
\|A^{\fr{k}{2}} u^{\epsilon}(t) \|_{L^2}^2 \le \Gamma_{k,1} e^{-\min\left\{c_1, c_2\right\}t}
\ee for any $t \ge 0$. Integrating \eqref{t23} in time from $0$ to $t$, we also conclude that  
\be \la{t25}
\int_{0}^{t} \|A^{\fr{k}{2} + \fr{1}{2}} u^{\epsilon}(s)\|_{L^2}^2 ds 
\le \Gamma_{k,2}
\ee for any $t \ge 0$. Here $\Gamma_{k,1}$ and $\Gamma_{k,2}$ are positive constants which do not depend on $\epsilon$ nor on the time $t$ but only on the initial data and the order of regularity $k$. Since $\pa \Omega$ is smooth, $\mathcal{D}(A^{\fr{k}{2} + \fr{1}{2}})$ is continuously embedded in $H^{k+1} \cap \mathcal{D}(A^{\fr{1}{2}})$, yielding  consequently the desired estimate \eqref{tttt2}.

The evolution of the spatial $L^2$ norm of $\l^k q^{\epsilon}$ is described by the ODE
\be 
\fr{1}{2} \fr{d}{dt} \|\l^k q^{\epsilon}\|_{L^2}^2 
+ \|\l^{k + \fr{\alpha}{2}} q^{\epsilon}\|_{L^2}^2
+ \epsilon \|\l^{k+1} q^{\epsilon}\|_{L^2}^2 
= \int_{\Omega} \l^{k-1} (u^{\epsilon} \cdot \na q^{\epsilon}) \l^{k+1} q^{\epsilon} dx.
\ee 
Suppose $k$ is even. Since $q^{\epsilon} \in \mathcal{D}(\l^{k+1})$, we have $u^{\epsilon} \cdot \na q^{\epsilon} \in \mathcal{D}(\l^{k-1})$ in view of the equation \eqref{eq4} obeyed by $q^{\epsilon}$. From Proposition \ref{cru}, we conclude that $u^{\epsilon} \cdot \na q^{\epsilon} \in \mathcal{D}(\l^{k})$. We apply the commutator estimate \eqref{pr31} and estimate
\beg{align}
&\left|\int_{\Omega} \l^k (u^{\epsilon} \cdot \na q^{\epsilon}) \l^k q^{\epsilon} dx \right|
= \left|\int_{\Omega} \left[\l^k (u^{\epsilon} \cdot \na q^{\epsilon}) - u^{\epsilon} \cdot \na \l^k q^{\epsilon} \right] \l^k q^{\epsilon} dx \right| \nonumber
\\&\quad\quad\le \|\l^k q^{\epsilon}\|_{L^2} \|\l^k (u^{\epsilon} \cdot \na q^{\epsilon}) -  u^{\epsilon} \cdot \na \l^k q^{\epsilon}  \|_{L^2} \nonumber
\\&\quad\quad\le C\|\l^k q^{\epsilon}\|_{L^2} \| u^{\epsilon}\|_{H^{k+1}} \|\l^{k + \fr{\alpha}{2}} q^{\epsilon}\|_{L^2} \nonumber
\\&\quad\quad\le \fr{1}{2} \|\l^{k + \fr{\alpha}{2}} q^{\epsilon}\|_{L^2}^2
+ C\| u^{\epsilon}\|_{H^{k+1}}^2 \|\l^k q^{\epsilon}\|_{L^2}^2.
\end{align} Now suppose that $k$ is odd. Then it holds that $u^{\epsilon} \cdot \na q^{\epsilon} \in \mathcal{D}(\l^{k-1})$, thus 
\beg{align}
&\left|\int_{\Omega} \l^{k-1} (u^{\epsilon} \cdot \na q^{\epsilon}) \l^{k+1} q^{\epsilon} dx\right|
= \left|\int_{\Omega} \na \l^{k-1} (u^{\epsilon} \cdot \na q^{\epsilon}) \cdot \na \l^{k-1} q^{\epsilon} dx \right|
\nonumber
\\&\quad\quad= \left|\int_{\Omega} \left[\na \l^{k-1} (u^{\epsilon} \cdot \na q^{\epsilon}) - u^{\epsilon} \cdot \na  \na \l^{k-1} q^{\epsilon} \right]\cdot \na \l^{k-1} q^{\epsilon} dx \right| \nonumber
\\&\quad\quad\le C\|\l^k q^{\epsilon}\|_{L^2} \|u^{\epsilon}\|_{H^{k+1}} \|\l^{k + \fr{\alpha}{2}} q^{\epsilon}\|_{L^2} \nonumber
\\&\quad\quad\le \fr{1}{2} \|\l^{k + \fr{\alpha}{2}} q^{\epsilon}\|_{L^2}^2
+ C\|u^{\epsilon}\|_{H^{k+1}}^2 \|\l^k q^{\epsilon}\|_{L^2}^2, 
\end{align}  in view of the estimate \eqref{pr32}.  This yields the differential inequality
\be \la{t26}
\fr{d}{dt} \|\l^k q^{\epsilon}\|_{L^2}^2 
+ \|\l^{k + \fr{\alpha}{2}} q^{\epsilon}\|_{L^2}^2 
\le C\|u^{\epsilon}\|_{H^{k+1}}^2 \|\l^k q^{\epsilon}\|_{L^2}^2, 
\ee which implies that 
\be 
\fr{d}{dt} \|\l^k q^{\epsilon}\|_{L^2}^2 
+ \left(c_1 - C\|u^{\epsilon}\|_{H^{k+1}}^2\right) \|\l^k q^{\epsilon}\|_{L^2}^2 
\le 0.
\ee We multiply by the integrating factor $e^{c_1t - C\int_{0}^{t} \|u^{\epsilon}\|_{H^{k+1}}^2 ds}$, integrate in time from $0$ to $t$, make use of \eqref{t25}, and conclude that 
\be 
\|\l^k q^{\epsilon}(t)\|_{L^2}^2  \le \Gamma_{k,3} e^{-c_1t}
\ee for any $t \ge 0$. Also here $\Gamma_{k,3}$ is a positive constant depending only on the initial data and $k$. We have thus completed the proof of Theorem \ref{tt22}.

\section{Global Attractor: Proof of Theorem \ref{att}} \la{s5}

In this section, we address the long time dynamics of the forced model \eqref{inmod11}--\eqref{inmod21}. We take the potential $\Phi$ to be zero for simplicity (see Remark \ref{potzero} below).

For $\alpha \in (0,1]$, we consider the forced system 
\be \la{system}
\begin{cases}
\pa_t q + u \cdot \na q + \l^{\alpha} q = 0
\\ \pa_t u + u \cdot \na u + \na p - \Delta u = - qRq + f, 
\\ \na \cdot u = 0
\end{cases}
\ee in the presence of smooth time independent divergence-free body forces $f$ in the fluid. The system \eqref{system} is posed on $\Omega \times [0, \infty)$, with homogeneous Dirichlet boundary conditions and initial data $q_0$ and $u_0$. We address the long time behavior of solutions.

We recall the spaces $\mathcal{H}$ and $\mathcal{V}$ defined respectively in \eqref{functio1} and \eqref{functio2} and  the solution map 
\be 
\mathcal{S}(t): \mathcal{V} \mapsto \mathcal{V}
\ee associated with \eqref{system} and defined by 
\be 
\mathcal{S}(t) (q_0, u_0) = (q(t), u(t)),
\ee where $\omega(t):= (q(t), u(t))$ is the unique solution of \eqref{system} with initial datum $\omega_0 := (q_0, u_0)$ at time $t$.

\subsection{Existence of an Absorbing Ball.} 

We start by proving the existence of a ball $B_{\rho}$, compact in $\mathcal{H}$, such that the image of $B_{\rho}$ under $\mathcal{S}(t)$ lies in $B_{\rho}$ itself for all large times.

\beg{prop} \la{absball}  Suppose $\alpha \in \left(0, 1\right]$ and $f \in \mathcal{D}(A^{\fr{1}{2}})$. Then there exists a radius $\rho > 0$ depending only on the $H^1$ norm of $f$ and some universal constants such that for each $\omega_0 = (q_0, u_0) \in \mathcal{V}$, there exists a time $T_0$ depending only on $\|\na q_0\|_{L^2}$ and $\|\na u_0\|_{L^2}$ and universal constants such that 
\be \la{absball1}
\mathcal{S}(t) \omega_0 \in \mathcal{B}_{\rho} := \left\{(q, u) \in \mathcal{V}: \|\na q\|_{L^2} + \|\Delta u\|_{L^2} \le \rho \right\}
\ee for all $t \ge T_0$.
\end{prop}

\noindent \textbf{Proof.} Fix $(q_0, u_0) \in \mathcal{V}$. In view of the first three steps established in the proof of Theorem \ref{Existence}, and in the presence of body forces in the fluid, we have the following bounds 
\be \la{absbal1}
\|q(t)\|_{L^2}^2 \le \|q_0\|_{L^2}^2 e^{-ct},
\ee 
\be \la{profor5}
\int_{s}^{t} \|\l^{\fr{\alpha}{2}} q (\tau)\|_{L^2}^2 d\tau \le \|q_0\|_{L^2}^2 e^{-cs},
\ee 
\be 
\|u(t)\|_{L^2}^2 \le \left(\|u_0\|_{L^2}^2 + C\|q_0\|_{L^2}^4 \right)e^{-ct} + \|f\|_{L^2}^2,
\ee  
\be 
\int_{s}^{t} \|\na u(\tau)\|_{L^2}^2 d\tau \le \left(\|u_0\|_{L^2}^2 + C\|q_0\|_{L^2}^4 \right)e^{-cs} + \|f\|_{L^2}^2 + \|f\|_{L^2}^2 (t-s),
\ee and
\be \la{absbal2}
\|q(t)\|_{L^p} \le \|q_0\|_{L^p}e^{-\fr{c (p-1)}{p} t}
\ee for any $t \ge 0$, any $s \in [0,t]$, and any even integer $p \ge 4$. The constant $c$  depends only on the size of the domain $\Omega$ and the power $\alpha$. Based on \eqref{absbal1}--\eqref{absbal2}, we deduce the existence of a time $t_0$ depending on $\|\omega_0\|_{\mathcal{V}}$ and a radius $R$ depending only on $\|f\|_{L^2}$ such that the bound
\be \la{profor111}
\|q(t)\|_{L^2}^2 + \|u(t)\|_{L^2}^2 + \int_{t}^{t+1} \left(\|\l^{\fr{\alpha}{2}}q(s)\|_{L^2}^2 + \|\na u(s)\|_{L^2}^2 \right) ds \le R
\ee holds for all $t \ge t_0$.

\textbf{Step 1. Velocity $H^1$ bounds.} The analogous energy inequality of \eqref{t210} in the presence of forces $f$ is given by  
\be \la{absbal4}
\fr{d}{dt} \|A^{\fr{1}{2}} u\|_{L^2}^2
+ \|Au\|_{L^2}^2
\le C\|q\|_{L^{\tilde{q}}}^2\| \l^{\fr{\alpha}{2}}q\|_{L^2}^2
+ C\|u\|_{L^2}^{2} \|A^{\fr{1}{2}} u\|_{L^2}^4  + C\|f\|_{L^2}^2
\ee where $\tilde{q}$ is some large even integer.  In view of \eqref{absbal2}, there exists a time $t_1 \ge t_0$ depending on $\tilde{q}$ and  $\|\omega_0\|_{\mathcal{V}}$ such that 
\be \la{profor112}
\|q(t)\|_{L^{\tilde q}}  \le 1
\ee for all $t \ge t_1$. Due to \eqref{profor111}, the differential inequality \eqref{absbal4} is of the same type as \eqref{gron1} with $y = \|A^{\fr{1}{2}} u\|_{L^2}^2$, $F_1 = \|q\|_{L^{\tilde{q}}}^2\| \l^{\fr{\alpha}{2}}q\|_{L^2}^2$, $F_2 = \|u\|_{L^2}^{2}$, $n=2$, $C_1 = C\|f\|_{L^2}$, and $C_2 = C_3 = C$, for some positive universal constant $C$. As a consequence of Lemma \ref{gron}, we infer the existence of a time $t_2 \ge t_1$ depending only on $\|\omega_0\|_{\mathcal{V}}$ such that 
\be \la{absbal88}
\|\na u(t)\|_{L^2}^2 + \int_{t}^{t+1} \|\Delta u(s)\|_{L^2}^2 ds \le R_1
\ee for all $t \ge t_2$, where $R_1$ is a positive constant depending only on the body forces and universal constants.   

\textbf{Step 2. Charge Density $L^{\infty}$ bounds.} There exists a time $t_3 \ge t_2$ such that $\l^{1+\fr{\alpha}{2}}q(t_3)$ is bounded in $L^2$ by some constant depending on $\|\na q_0\|_{L^2}, \|\na u_0\|_{L^2},$ and $\|f\|_{L^2}$, a fact that follows by repeating the energy calculations obtained in Steps 4 and 5 of Theorem \ref{Existence} but in the presence of body forces in the fluid. At this specific time $t_3$, the charge density is $L^{\infty} \cap H^1$ regular due to continuous Sobolev embeddings. By making use of the $L^p$ estimates \eqref{absbal2}, we have 
\be \la{profor7}
\|q(t)\|_{L^{\infty}} \le \|q(t_3)\|_{L^{\infty}}e^{-c(t-t_3)}
\ee for all $t \ge t_3$. From \eqref{profor7}, we deduce the existence of a time $t_4 \ge t_3$ such that 
\be \la{proo}
\|q(t)\|_{L^{\infty}} \le 1
\ee for all $t \ge t_4$.

\textbf{Step 3. Velocity $\mathcal{D}(A^{\fr{1}{2} + \fr{\alpha}{4} - \delta})$ bounds.} Let $\delta > 0$ be sufficiently small. In view of \eqref{absbal88}, there exists a time $t_5 \ge t_4$ such that $\|Au(t_5)\|_{L^2}$ is bounded by some constant depending on $\|\omega_0\|_{\mathcal{V}}$ and $f$. We study the evolution of the norm $\|A^{\fr{1}{2} + \fr{\alpha}{4} - \delta}u\|_{L^2}^2$ starting at time $t_5$. Indeed, we have
\beg{align} \la{profor1}
&\fr{1}{2} \fr{d}{dt} \|A^{\fr{1}{2} + \fr{\alpha}{4} - \delta}u\|_{L^2}^2
+ \|A^{1 + \fr{\alpha}{4} - \delta}u\|_{L^2}^2 \nonumber
\\&= - \int_{\Omega} A^{\fr{\alpha}{4} - \delta} \PP(qRq) \cdot A^{1+\fr{\alpha}{4} - \delta} u dx
- \int_{\Omega} A^{\fr{\alpha}{4} - \delta} \PP(u \cdot \na u) \cdot  A^{1+\fr{\alpha}{4} - \delta} u dx
- \int_{\Omega} A^{\fr{\alpha}{4} - \delta}f \cdot  A^{1+\fr{\alpha}{4} - \delta} u dx,
\end{align} from which we obtain the differential inequality
\be \la{profor2}
\fr{d}{dt} \|A^{\fr{1}{2} + \fr{\alpha}{4} - \delta}u\|_{L^2}^2
+ \|A^{1 + \fr{\alpha}{4} - \delta}u\|_{L^2}^2 
\le C\| A^{\fr{\alpha}{4} - \delta} \PP(qRq)\|_{L^2}^2
+  C\|A^{\fr{\alpha}{4} - \delta}\PP (u \cdot \na u)\|_{L^2}^2 
+ C\|A^{\fr{\alpha}{4} - \delta}f\|_{L^2}^2
\ee due to Young's inequality. We bound the nonlinear term in $q$ by using the fractional product estimate \eqref{prodfor7} with $s = \fr{\alpha}{2}$ and obtain 
\be 
\| A^{\fr{\alpha}{4} - \delta} \PP(qRq)\|_{L^2}^2 
\le C\|q\|_{L^{\infty}}^2\|\l^{\fr{\alpha}{2}}q\|_{L^2}^2.
\ee As for the nonlinear term in $u$, we use the divergence-free condition obeyed by $u$ and estimate
\be \la{profor3}
\|A^{\fr{\alpha}{4} - \delta}\PP (u \cdot \na u)\|_{L^2}^2 
\le C\|\na \cdot (u \otimes u)\|_{H^{\fr{\alpha}{2} - 2\delta}}^2
\le C\| u\|_{H^{1 + \fr{\alpha}{2} - 2\delta}}^4
\le C\|A^{\fr{1}{2} + \fr{\alpha}{4} - \delta} u\|_{L^2}^4,
\ee where the second inequality follows from the fact that $H^{1+ \fr{\alpha}{2} - 2\delta}$ is a Banach Algebra for a sufficiently small $\delta$, and the last inequality uses the continuous embedding of $\mathcal{D}(A^{\fr{1}{2} + \fr{\alpha}{4} - \delta})$ into $H^{1 + \fr{\alpha}{2} - 2\delta}$ (\cite{GS}). Putting \eqref{profor2}--\eqref{profor3} together gives 
\be 
\fr{d}{dt} \|A^{\fr{1}{2} + \fr{\alpha}{4} - \delta}u\|_{L^2}^2
+ \|A^{1 + \fr{\alpha}{4} - \delta}u\|_{L^2}^2 
\le C\|A^{\fr{1}{2} + \fr{\alpha}{4} - \delta} u\|_{L^2}^4
+ C\|q\|_{L^{\infty}}^2\|\l^{\fr{\alpha}{2}}q\|_{L^2}^2
+ C\|A^{\fr{\alpha}{4} - \delta}f\|_{L^2}^2.
\ee Since 
\be 
\int_{t}^{t+1} \|A^{\fr{1}{2} + \fr{\alpha}{4} - \delta} u(s)\|_{L^2}^2 ds
\le C\int_{t}^{t+1} \|Au(s)\|_{L^2}^2 ds,
\ee the conditions of the uniform Gronwall Lemma \ref{gron} are satisfied for all $t \ge t_5$ in view of \eqref{profor111}, \eqref{absbal88}, and \eqref{proo}. Therefore, we deduce the existence of a time $t_6 \ge t_5$ depending on $\|\omega_0\|_{\mathcal{V}}$, and a radius $R_2$ depending only on $f$ such that 
\be  \la{proforr}
 \|A^{\fr{1}{2} + \fr{\alpha}{4} - \delta}u(t)\|_{L^2}^2 + \int_{t}^{t+1} \|A^{1 + \fr{\alpha}{4} - \delta}u(s)\|_{L^2}^2 ds \le R_2
\ee for all $t \ge t_6$. 

\textbf{Step 4. Velocity gradient $L^p$ bounds.} The velocity $u$ can be represented as
\be 
u(t) = e^{-(t-\tau)A}u(\tau) - \int_{\tau}^{t} e^{-(t-s)A} \left(B(u,u)+ \PP(qRq) - f\right)(s) ds
\ee for any $t \in [\tau, \infty)$. In view of the Stokes semi-group estimate
\be 
\|A^{\fr{1}{2}} e^{-tA} v\|_{L^p} \le Ct^{-\fr{1}{2}} \|v\|_{L^p}
\ee that holds for $p \in (1, \infty)$ (see \cite[Proposition 1.2]{GM}), we have 
\beg{align}
\|A^{\fr{1}{2}} u\|_{L^p}
&\le C\fr{\|u(\tau)\|_{L^{p}}}{\sqrt{t-\tau}} + \int_{\tau}^{t} \fr{1}{\sqrt{t-s}} \|B(u,u)\|_{L^p} ds \nonumber
\\&\quad\quad+ \int_{\tau}^{t} \fr{1}{\sqrt{t-s}} \|\PP(qRq)\|_{L^p} ds
+ \int_{\tau}^{t} \fr{1}{\sqrt{t-s}} \|f\|_{L^p} ds.
\end{align} We estimate 
\beg{align} 
\|B(u,u)\|_{L^p}
&\le C\|u \cdot \na u\|_{L^p}
\le C\|u\|_{L^{\infty}} \|\na u\|_{L^p}
\le C\|u\|_{L^{\infty}} \|\na u\|_{L^2}^{\fr{1}{p}}\|\na u\|_{L^{2p-2}}^{\fr{p-1}{p}} \nonumber
\\&\le C\|u\|_{L^{\infty}} \|\na u\|_{L^2}^{\fr{1}{p}}\|\na u\|_{L^\infty}^{\fr{p-1}{p}}
\le C\|u\|_{H^{1+\epsilon}} \|u\|_{H^{2+\epsilon}}^{\fr{p-1}{p}} \|\na u\|_{L^2}^{\fr{1}{p}}
\end{align} via interpolation of $L^p$ spaces and use of Sobolev embeddings. Due to H\"older's inequality with exponents $\fr{2p}{p+1}$ and $\fr{2p}{p-1}$, we obtain
\beg{align}
&\int_{\tau}^{t} \fr{1}{\sqrt{t-s}} \|B(u,u)\|_{L^p} ds 
\\&\le C\sup\limits_{s \in [\tau, t]} \left(\|\na u(s)\|_{L^2}^{\fr{1}{p}}  \|u(s)\|_{H^{1+\epsilon}}  \right) \left(\int_{\tau}^{t} \fr{1}{(t-s)^{\fr{p}{p+1}}} ds \right)^{\fr{p+1}{2p}}  \left(\int_{\tau}^{t} \|u\|_{H^{2+\epsilon}}^2 ds \right)^{\fr{p-1}{2p}} \nonumber
\\&\le C\sup\limits_{s \in [\tau, t]} \left(\|\na u\|_{L^2}^{\fr{1}{p}} \|u\|_{H^{1+\epsilon}}  \right) \left(\int_{\tau}^{t} \|u\|_{H^{2+\epsilon}}^2 ds \right)^{\fr{p-1}{2p}} (t-\tau)^{\fr{1}{2p}}.
\end{align} As for the nonlinear term in $q$, we have 
\be 
\int_{\tau}^{t} \fr{1}{\sqrt{t-s}} \|\PP(qRq)\|_{L^p} ds
\le C\sup\limits_{s \in [\tau, t]} \|q(s)\|_{L^{\infty}}^2 \sqrt{t-\tau}
\ee due to the finiteness of the domain size and the boundedness of the Dirichlet Riesz transform on $L^p$ spaces. We obtain the bound
\beg{align} 
\|A^{\fr{1}{2}} u(t)\|_{L^p} 
&\le C\fr{\|u(\tau)\|_{L^{p}}}{\sqrt{t-\tau}}
+ C\sup\limits_{s \in [\tau, t]} \left(\|\na u\|_{L^2}^{\fr{1}{p}}  \|u\|_{H^{1+\epsilon}}  \right) \left(\int_{\tau}^{t} \|u\|_{H^{2+\epsilon}}^2 ds \right)^{\fr{p-1}{2p}} (t-\tau)^{\fr{1}{2p}} \nonumber
\\&\quad\quad+ C\sup\limits_{s \in [\tau, t]} \|q(s)\|_{L^{\infty}}^2 \sqrt{t-\tau}
+ C\|f\|_{L^p} \sqrt{t-\tau}
\end{align} for any $\tau > 0$ and $t \ge  \tau$. Fix a nonnegative integer $k \ge 0$. Taking $\tau = t_6 + k$ and $t \in [t_6+k+1, t_6+k+2]$, and noting that $1 \le \sqrt{t-\tau} \le \sqrt{2}$, we have 
\beg{align} 
\|A^{\fr{1}{2}} u(t)\|_{L^p} 
&\le C\|u(t_6 + k)\|_{H^1}
+ CR_2^{\fr{1}{p}+1} \left(\int_{t_6+k}^{t_6 + k +2} \|u\|_{H^{2+\epsilon}}^2 ds \right)^{\fr{p-1}{2p}} 
+  \sqrt{2}C
+ \sqrt{2} C\|f\|_{L^p}. 
\end{align}
Due to the boundedness of the local in time integral $\int_{t_6+k}^{t_6 + k +2} \|u\|_{H^{2+\epsilon}}^2 ds$ independently of $t_6$ and $k$,  we infer the existence of a radius $R_3>0$ depending only on the body forces such that  
\be 
\sup\limits_{t \in [t_6 + k+1, t_6 + k + 2]}  \|A^{\fr{1}{2}} u(t)\|_{L^p}
\le R_3.
\ee This is true for any integer $k \ge 0$, thus
\be 
\|A^{\fr{1}{2}} u(t)\|_{L^p} \le R_3
\ee for any $t \ge t_7$ where $t_7:=t_6 + 1$. 

\textbf{Step 5. Charge Density $H^1$ bounds.} The evolution of the $L^2$ norm of $\na q$ described by \eqref{t214} is not of type \eqref{gron1} due to the absence of the local in time integrability property \eqref{gron33} for $y = \|\na q\|_{L^2}^2$. Hence Lemma~\ref{gron} does not apply in this case. In order to show that the $L^2$ norm of $\na q$ is uniformly bounded for large times, independently of the initial data, we need to estimate the nonlinear term differently. Indeed, $\|\na q\|_{L^2}$ obeys the differential inequality 
\be 
\fr{1}{2} \fr{d}{dt} \|\na q\|_{L^2}^2 + \|\l^{1 + \fr{\alpha}{2}}q\|_{L^2}^2
\le \int_{\Omega} |\na u| |\na q|^2 dx.
\ee
By the Brezis-Mironescu interpolation inequality, we have
\be 
\|\na q\|_{L^{2+\alpha}} \le C\|q\|_{L^{\infty}}^{\fr{\alpha}{2+\alpha}} \|\l^{1 + \fr{\alpha}{2}} q\|_{L^2}^{\fr{2}{2+\alpha}}.
\ee By making use of H\"older's inequality with exponents $\fr{2+\alpha}{\alpha}, 2 + \alpha, 2+\alpha$, we obtain
\be 
\fr{1}{2} \fr{d}{dt} \|\na q\|_{L^2}^2 + \|\l^{1 + \fr{\alpha}{2}}q\|_{L^2}^2
\le C\|\na u\|_{L^{\fr{2+\alpha}{\alpha}}} \|q\|_{L^{\infty}}^{\fr{2\alpha}{2+\alpha}} \|\l^{1 + \fr{\alpha}{2}} q\|_{L^2}^{\fr{4}{2+\alpha}},
\ee which, followed by an application of Young's inequality with exponents $\fr{2+\alpha}{2}$ and $\fr{2+\alpha}{\alpha}$, reduces to
\be 
\fr{d}{dt} \|\na q\|_{L^2}^2 + \|\l^{1 + \fr{\alpha}{2}}q\|_{L^2}^2
\le C\|\na u\|_{L^{\fr{2+\alpha}{\alpha}}}^{\fr{2+\alpha}{\alpha}} \|q\|_{L^{\infty}}^{2}.
\ee Since $\|\na u\|_{L^p} \le C\|A^{\fr{1}{2}}u\|_{L^p}$ for any $p \in (1, \infty)$ (\cite[Proposition 1.4]{GM}), it holds that
\be 
\int_{t}^{t+1}  \|\na u(s)\|_{L^{\fr{2+\alpha}{\alpha}}}^{\fr{2+\alpha}{\alpha}} \|q(s)\|_{L^{\infty}}^{2} ds
\le \rho
\ee for some $\rho$ depending only on the body forces and the power $\alpha$, thus the uniform Gronwall Lemma \ref{gron} is applicable and yields the existence of a time $t_8 \ge t_7$ depending only on $\|\na q_0\|_{L^2}$, $\|\na u_0\|_{L^2}$, and a radius $R_4 > 0$ depending only on $f$ such that 
\be 
\|\na q(t)\|_{L^2} + \int_{t}^{t+1} \|\l^{1 + \fr{\alpha}{2}}q(s) \|_{L^2}^2 ds \le R_4
\ee holds for all times $t \ge t_8$.

\textbf{Step 6. Velocity $H^2$ bounds.} The following energy inequality 
\be 
\fr{d}{dt} \|Au\|_{L^2}^2 + \|A^{\fr{3}{2}}u\|_{L^2}^2
\le C\|\na q\|_{L^2}^2 \|\l^{1+\fr{\alpha}{2}} q \|_{L^2}^2 + C\|Au\|_{L^2}^4 + C\|\na f\|_{L^2}^2
\ee holds and is of type \eqref{gron1}. In view of Lemma \ref{gron} with $F_1 = \|\na q\|_{L^2}^2 \|\l^{1+\fr{\alpha}{2}} q \|_{L^2}^2$, $F_2 = 1$ and $n=2$, we obtain \eqref{absball1}. This ends the proof of Proposition \ref{absball}.

\beg{rem} As a consequence of Proposition \ref{absball}, we infer the existence of a positive time $T$ such that 
\be 
\mathcal{S}(t) \mathcal{B}_{\rho} \subset \mathcal{B}_{\rho}
\ee for all $t \ge T$.
\end{rem}

\beg{rem} The case of $\TT^2$ with periodic boundary conditions is simpler: for any $\alpha \in (0,1]$, there exists a radius $R>0$ depending only on $\|\na f\|_{L^2}$, such that for each $\omega_0 = (q_0, u_0) \in \mathcal{V}$, there exists a time $T$ depending only on $\|\na q_0\|_{L^2}$ and $\|\na u_0\|_{L^2}$ and universal constants such that 
\be \la{perfor1}
\mathcal{S}(t) \omega_0 \in \mathcal{B}_{R} := \left\{(q, u) \in \mathcal{V}: \|\l^{1+\fr{\alpha}{2}} q\|_{L^2} + \|\Delta u\|_{L^2} \le R \right\}
\ee for all $t \ge T$. Indeed, for any $s \in [\fr{\alpha}{2},1 + \fr{\alpha}{2}]$, the energy equality 
\be 
\fr{1}{2} \fr{d}{dt} \|\l^{s}q\|_{L^2}^2 + \|\l^{s+\fr{\alpha}{2}}q\|_{L^2}^2
= - \int_{\TT^2} \left[\l^{s} (u \cdot \na q) - u \cdot \na \l^{s} q \right] \l^{s} q dx
\ee holds and yield 
\be 
\fr{1}{2} \fr{d}{dt} \|\l^{s}q\|_{L^2}^2 + \|\l^{s+\fr{\alpha}{2}}q\|_{L^2}^2
\le C\|\Delta u\|_{L^2} \|\l^{s + \fr{\alpha}{2}} q\|_{L^2} \|\l^{s}q\|_{L^2}
\ee due to standard periodic commutator estimates. An application of Young's inequality gives the differential inequality 
\be \la{perfor}
\fr{d}{dt} \|\l^{s}q\|_{L^2}^2 + \|\l^{s+\fr{\alpha}{2}}q\|_{L^2}^2
\le C\|\Delta u\|_{L^2}^2 \|\l^{s}q\|_{L^2}^2.
\ee Choosing $s = \fr{\alpha}{2}$ and using \eqref{profor111} and \eqref{absbal88}, we deduce that \eqref{perfor} is of type \eqref{gron1} and obtain good control of both $\|\l^{\fr{\alpha}{2}}q\|_{L^2}$ and $\int_{t}^{t+1} \|\l^{\alpha}q\|_{L^2}^2 ds$. Then we take $s = \alpha$, repeat the same argument, and obtain control of $\|\l^{\alpha}q\|_{L^2}$ and $\int_{t}^{t+1} \|\l^{\fr{3\alpha}{2}}q\|_{L^2}^2 ds$. A bootstrapping argument yields the existence of a time $T_1$ depending on $\|\omega_0\|_{\mathcal{V}}$ and a radius $R'$ depending only on $f$ such that 
\be 
\|\l^{1+\fr{\alpha}{2}}q(t)\|_{L^2}^2 + \int_{t}^{t+1} \|\l^{1+\alpha}q(s)\|_{L^2}^2 ds \le R'
\ee for all $t \ge T_1$. Therefore, we deduce \eqref{perfor1} and obtain the absorbing ball $B_{R}$ which is compact in the strong norm of $\mathcal{V}$. That is not the case on bounded smooth domains $\Omega$ with homogeneous Dirichlet boundary conditions, due to the more complicated commutator estimates. 
\end{rem}

\beg{rem} One of the main elements of the proof of Proposition \ref{absball} is the boundedness of the velocity gradient in $L^p$ spaces for all $p \in (2, \infty)$. The maximal $L^p$ regularity has been studied in the literature for the Stokes equations (\cite{HiSa,Mar,ShSh} and references therein), the Navier-Stokes equations (\cite{FaKoWe,GeHeHi,GM,kim,Tol,ToWa} and references therein), and parabolic evolution equations (\cite{DeHiPr,HiPr} and references therein) on bounded and unbounded domains,  under various regularity conditions imposed on the boundaries, and equipped with different types of boundary conditions. 
\end{rem}

\subsection{Continuity Properties of the Solution Map.} We investigate the instantaneous continuity of the solution map $\mathcal{S}(t)$ at each fixed positive time $t$. 

\begin{prop} \la{conts}
Let $w_1^0 = (q_1^0, u_1^0), w_2^0 = (q_2^0, u_2^0) \in \mathcal{V}$. Let $t >0$. 
There exist functions $K_1(t)$, $K_2(t)$, and $K_3(t)$, locally uniformly bounded as functions of $t \geq 0$, and locally bounded as initial data $w_1^0, w_2^0$ are varied in $\mathcal{V}$, such that 
$\mathcal{S}(t)$ is Lipschitz continuous in $\mathcal{H}$ obeying 
\be \la{Conth}
\|\mathcal{S}(t)w_1^0 - \mathcal{S}(t)w_2^0\|_{\mathcal{H}}^2 \leq K_1(t) \|w_1^0 - w_2^0\|_{\mathcal{H}}^2,
\ee
$\mathcal{S}(t)$ is Lipschitz continuous in $\mathcal{V}$ obeying 
\be \la{ContV}
\|\mathcal{S}(t)w_1^0 - \mathcal{S}(t)w_2^0\|_{\mathcal{V}}^2 \leq K_2(t) \|w_1^0 - w_2^0\|_{\mathcal{V}}^2,
\ee and $\mathcal{S}(t)$ is Lipschitz continuous from $\mathcal{H}$ to $\mathcal{V}$ obeying
\be \la{Conthv}
t^{\delta} \|\mathcal{S}(t)w_1^0 - \mathcal{S}(t)w_2^0\|_{\mathcal{V}}^2 \leq K_3(t) \|w_1^0 - w_2^0\|_{\mathcal{H}}^2
\ee for any $\delta > \fr{2}{\alpha}.$
\end{prop} 

\noindent \textbf{Proof.} We set $(q_1, u_1) = \mathcal{S}(t) (q_1^0, u_1^0)$ and $(q_2, u_2) = \mathcal{S}(t) (q_2^0, u_2^0)$. The differences $q = q_1 - q_2$ and $u = u_1 - u_2$ obey the system
\be \la{con1}
\beg{cases}
\pa_t q + \l^{\alpha}q = - u_1 \cdot \na q - u \cdot \na q_2,
\\\pa_t u + Au = - B(u_1, u) - B(u, u_2) - \PP(q_1 Rq) - \PP(qRq_2).
\end{cases}
\ee The following differential inequality 
\beg{align} \la{con4}
&\fr{d}{dt} \left(\|q\|_{L^2}^2 + \|u\|_{L^2}^2\right) 
+ \|\l^{\fr{\alpha}{2}}q\|_{L^2} + \|\na u\|_{L^2}^2 \nonumber
\\&\quad\le C\left(\|\l^{1 + \fr{\alpha}{2}}q_2\|_{L^{2}}^2 + \|\na u_2\|_{L^2}^2 + \|q_1\|_{L^{4}}^2 + 1\right) \left(\|q\|_{L^2}^2 + \|u\|_{L^2}^2\right) 
\end{align} holds, as shown in \eqref{co}. Consequently, the Lipschitz continuity of $\mathcal{S}(t)$ in the norm of $\mathcal{H}$, given by  \eqref{Conth},  follows with 
\be 
K_1(t) = \exp \left\{C\int_{0}^{t} \left(\|\l^{1 + \fr{\alpha}{2}}q_2\|_{L^{2}}^2 + \|\na u_2\|_{L^2}^2 + \|q_1\|_{L^{4}}^2 + 1\right) ds  \right\}.
\ee 
In order to study the Lipschitz continuity of $\mathcal{S}(t)$ in the norm of $\mathcal{V}$, we take the $L^2$ inner product of the first and second equations in \eqref{con1} with $-\Delta q$ and $Au$ respectively, add the resulting energy equalities, and estimate. We obtain 
\beg{align} \la{co22}
&\fr{1}{2} \fr{d}{dt} \left(\|\na q\|_{L^2}^2 + \|\na u\|_{L^2}^2 \right)
+ \|\l^{1 + \fr{\alpha}{2}} q\|_{L^2}^2 + \|A u\|_{L^2}^2 \nonumber
\\&\quad\le \int_{\Omega} |\na u_1| |\na q|^2 dx 
+ \int_{\Omega} |\na u| |\na q_2| |\na q| dx
+ \int_{\Omega} |u_1| |\na u| |Au| dx \nonumber
\\&\quad\quad+ \int_{\Omega}|u| |\na u_2| |Au| dx 
+ \int_{\Omega} |q_1| |Rq| |Au| dx 
+ \int_{\Omega} |q| |Rq_2| |Au| dx 
\end{align} after integration by parts. By making use of H\"older's inequality, continuous Sobolev embeddings, the boundedness of the Riesz transform on $L^4$, and the ellipticity of the Stokes operator, \eqref{co22} reduces to
\beg{align} \la{co2}
&\fr{1}{2} \fr{d}{dt} \left(\|\na q\|_{L^2}^2 + \|\na u\|_{L^2}^2 \right)
+ \|\l^{1 + \fr{\alpha}{2}} q\|_{L^2}^2 + \|A u\|_{L^2}^2 \nonumber
\\&\quad\le C\|\Delta u_1\|_{L^2} \|\na q\|_{L^2} \|\l^{1+\fr{\alpha}{2}}q\|_{L^2}
+ C\|A u\|_{L^2} \|\na q\|_{L^2} \|\l^{1+\fr{\alpha}{2}}q_2\|_{L^2}
+ C\|u_1\|_{L^4} \|\na u\|_{L^2}^{\fr{1}{2}} \|Au\|_{L^2}^{\fr{3}{2}} \nonumber
\\&\quad\quad+ C\|\na u\|_{L^2} \|\na u_2\|_{L^4}\|Au\|_{L^2}  
+ C\|q_1\|_{L^4} \|q\|_{L^4} \|Au\|_{L^2}
+ C \|q\|_{L^4} \|q_2\|_{L^4} \|Au\|_{L^2}.
\end{align} An application of Young's inequality yields 
\beg{align} \la{con6}
&\fr{d}{dt} \left(\|\na q\|_{L^2}^2 + \|\na u\|_{L^2}^2 \right) 
+ \|\l^{1 + \fr{\alpha}{2}} q\|_{L^2}^2 + \|\Delta u\|_{L^2}^2 \nonumber
\\&\quad\le C\left(\|\Delta u_1\|_{L^2}^2 + \|\Delta u_2\|_{L^2}^2 + \|u_1\|_{L^4}^4 + \|q_1\|_{L^4}^2 +  \|\l^{1+\fr{\alpha}{2}}q_2\|_{L^2}^2 \right) \left(\|\na q\|_{L^2}^2 + \|\na u\|_{L^2}^2 \right),
\end{align} which gives the desired $\mathcal{V}$-Lipschitz continuity property \eqref{ContV}, with 
\be 
K_2(t) = \exp \left\{C\int_{0}^{t} \left(\|\Delta u_1\|_{L^2}^2 + \|\Delta u_2\|_{L^2}^2 + \|u_1\|_{L^4}^4 + \|q_1\|_{L^4}^2 +  \|\l^{1+\fr{\alpha}{2}}q_2\|_{L^2}^2 \right) ds  \right\}.
\ee Finally, we prove the Lipschitzianity property \eqref{Conthv}. We seek a differential inequality of the form 
\be \la{con5}
\fr{d}{dt} \left(t^{\delta} \|\omega\|_{\mathcal{V}}^2 \right) 
\le Ct^{-\beta}  \|\omega\|_{\mathcal{H}}^2 + C\delta t^{\delta - 1} \|\na u\|_{L^2}^2  + Z(t) \left(t^{\delta} \|\omega\|_{\mathcal{V}}^2\right) 
\ee for some $\beta < 1$ and a locally integrable function in time $Z(t)$.  Solving \eqref{con5}, integrating \eqref{con4} in time from $0$ to $t$, and using \eqref{Conth}, we obtain 
\be 
t^{\delta} \|\omega\|_{\mathcal{V}}^2  \le K_3(t) \|\omega_0\|_{\mathcal{H}}^2
\ee where 
\be 
K_3(t) = C \left[\fr{t^{1-\beta}}{1-\beta} K_1(t)  + t^{\delta - 1} \ln K_1(t) \int_{0}^{t} K_1(s) ds + 1\right]\exp \left\{\int_{0}^{t} Z(s) ds\right\}.
\ee Indeed, the $L^2$ norm of $t^{\delta} \|\omega\|_{\mathcal{V}}^2$ obeys
\beg{align}  \la{dissneed}
&\fr{d}{dt} \left(t^{\delta} \|\omega\|_{\mathcal{V}}^2\right) + t^{\delta} \left(\|\l^{1+\fr{\alpha}{2}} q\|_{L^2}^2 + \|\Delta u\|_{L^2}^2 \right) \nonumber
\\&\quad\le \delta t^{\delta - 1} \|\omega\|_{\mathcal{V}}^2
+ C\left(\|\Delta u_1\|_{L^2}^2 + \|\Delta u_2\|_{L^2}^2 + \|u_1\|_{L^4}^4 + \|q_1\|_{L^4}^2 +  \|\l^{1+\fr{\alpha}{2}}q_2\|_{L^2}^2 \right) \left(t^{\delta}\|\omega\|_{\mathcal{V}}^2 \right) \nonumber
\\&\quad\le \delta t^{\delta - 1} \|\na q\|_{L^2}^2 +  \delta t^{\delta - 1} \|\na u\|_{L^2}^2 \nonumber
\\&\quad\quad\quad\quad\quad\quad+ C\left(\|\Delta u_1\|_{L^2}^2 + \|\Delta u_2\|_{L^2}^2 + \|u_1\|_{L^4}^4 + \|q_1\|_{L^4}^2 +  \|\l^{1+\fr{\alpha}{2}}q_2\|_{L^2}^2 \right) \left(t^{\delta}\|\omega\|_{\mathcal{V}}^2 \right)
\end{align} due to \eqref{con6}. This latter energy inequality is of type \eqref{con5} provided that we have good control of the term $\delta t^{\delta - 1} \|\na q\|_{L^2}^2$. By the Brezis-Mironescu inequality, we have 
\be 
\|\na q\|_{L^2}^2 \le C\|q\|_{L^2}^{\fr{2\alpha}{2+\alpha}} \|\l^{1+\fr{\alpha}{2}}q\|_{L^2}^{\fr{4}{2+\alpha}},
\ee hence
\be 
\delta t^{\delta - 1} \|\na q\|_{L^2}^2
\le C\delta t^{\delta - 1 - \fr{2\delta}{2+\alpha}} \|q\|_{L^2}^{\fr{2\alpha}{2+\alpha}} \left(t^{\fr{2\delta}{2+\alpha}}\|\l^{1+\fr{\alpha}{2}}q\|_{L^2}^{\fr{4}{2+\alpha}}\right)
\le \fr{1}{2}t^{\delta} \|\l^{1+\fr{\alpha}{2}}q\|_{L^2}^2 + C_{\alpha}\delta^{\fr{2+\alpha}{\alpha}} t^{\delta - 1 - \fr{2}{\alpha}} \|q\|_{L^2}^2
\ee via use of Young's inequality for sums with exponents $\fr{2+\alpha}{2}$ and $\fr{2+\alpha}{\alpha}$. Since $\delta$ is greater than $\fr{2}{\alpha}$,  \eqref{con5} holds with $\beta = 1 + \fr{2}{\alpha} - \delta$ and 
\be 
Z(t) = C\left(\|\Delta u_1\|_{L^2}^2 + \|\Delta u_2\|_{L^2}^2 + \|u_1\|_{L^4}^4 + \|q_1\|_{L^4}^2 +  \|\l^{1+\fr{\alpha}{2}}q_2\|_{L^2}^2 \right),
\ee completing the proof of Proposition \ref{conts}.

\beg{rem} The choice of $\delta$ in \eqref{Conthv} results from the need to control the local in time integrals of $\|\na q\|_{L^2}^2$ by constant multiples of the difference of the initial data in $\mathcal{H}$. However, the energy equality \eqref{con4} gives such a boundedness only for  $\int_{0}^{t} \|\l^{\fr{\alpha}{2}} q\|_{L^2}^2 ds$. The remedy is interpolation and control by the dissipation of the energy inequality in hand \eqref{dissneed}, which imposes restrictions on the power $\delta$ of the time singularity. 
\end{rem}

\subsection{Injectivity of the Solution Map.} We obtain injectivity of the solution map $\mathcal{S}(t)$ by adapting the approach of \cite{cfbook} to the system \eqref{system}. 

\beg{prop} \la{backun} Let $\omega_1^0 = (q_1(0), u_1(0)), \omega_2^0 = (q_2(0), u_2(0)) \in \mathcal{V}$. Suppose there exists a time $T>0$ such that $\mathcal{S}(T) \omega_1^0 = \mathcal{S}(T) \omega_2^0$. Then $\omega_1^0 = \omega_2^0$.
\end{prop}

\noindent \textbf{Proof.} The proof is divided into main steps.

\textbf{Step 1. Time analyticity of solutions.} Suppose $(q_0, u_0) \in \mathcal{V}$, and denote the solution of \eqref{system} at time $t$ by $(q(t), u(t))$. We complexify all functional spaces and operators, fix an angle $\theta \in (-\fr{\pi}{2}, \fr{\pi}{2})$, and take $t = se^{i\theta}$ for $s > 0$. We have
\beg{align}
&\fr{d}{ds} \left[\|\na q(se^{i\theta})\|_{L^2}^2 + \|A^{\fr{1}{2}} u(se^{i\theta})\|_{L^2}^2 \right] \nonumber
\\&= \fr{d}{ds} \left[(q(se^{i\theta}), -\Delta q(se^{i\theta}))_{L^2} + (u(se^{i\theta}), Au(se^{i\theta}))_{L^2} \right] \nonumber
\\&= \left(e^{i\theta} \fr{dq}{dt} (se^{i\theta}), -\Delta q(se^{i\theta})\right)_{L^2} 
+ \left(q(se^{i\theta}), - e^{i\theta}  \Delta  \fr{dq}{dt} (se^{i\theta})\right)_{L^2}  \nonumber
\\&\quad\quad+ \left(e^{i\theta}\fr{du}{dt} (se^{i\theta}), Au(se^{i\theta})\right)_{L^2}
+ \left(u(se^{i\theta}),e^{i\theta} A\fr{du}{dt} (se^{i\theta})\right)_{L^2} \nonumber
\\&= 2Re \left[e^{i\theta} \left(\fr{dq}{dt} (se^{i\theta}) , - \Delta q (se^{i\theta}) \right)_{L^2} + e^{i\theta} \left(\fr{du}{dt} (se^{i\theta}), Au (se^{i\theta}) \right)_{L^2} \right]
\end{align} where $Re(z)$ denotes the real part of a complex number $z \in \C$, and $(\cdot, \cdot)_{L^2}$ is the complexified $L^2$ inner product. Thus, the evolution of the norm $\|\na q(se^{i\theta})\|_{L^2}^2 + \|A^{\fr{1}{2}} u(se^{i\theta})\|_{L^2}^2$ is described by 
\beg{align} \la{backun1}
&\fr{1}{2} \fr{d}{ds} \left[\|\na q(se^{i\theta})\|_{L^2}^2 + \|A^{\fr{1}{2}} u(se^{i\theta})\|_{L^2}^2 \right] 
+ \cos \theta \left[\|\l^{1+ \fr{\alpha}{2}} q(se^{i\theta})\|_{L^2}^2 + \|A u(se^{i\theta})\|_{L^2}^2 \right]
\\&= Re \left[e^{i\theta} (u\cdot \na q, -\Delta q)_{L^2} - e^{i\theta}(B(u,u), Au)_{L^2} - e^{i\theta} (\PP(qRq), Au)_{L^2} + e^{i\theta}(f, Au)_{L^2} \right].
\end{align} We estimate 
\beg{align} 
&\left|(u\cdot \na q, -\Delta q)_{L^2} \right| 
\le C\|\na u\|_{L^{\fr{4}{\alpha}}} \|\na q\|_{L^2} \|\na q\|_{L^{\fr{4}{2-\alpha}}} \le C\|\na u\|_{L^2}^{\fr{\alpha}{2}} \|Au\|_{L^2}^{\fr{2-\alpha}{2}} \|\na q\|_{L^2} \|\l^{1+ \fr{\alpha}{2}}q\|_{L^2} \nonumber
\\&\le \fr{\cos \theta}{8} \left[\|\l^{1+ \fr{\alpha}{2}} q(se^{i\theta})\|_{L^2}^2 + \|A u(se^{i\theta})\|_{L^2}^2  \right] + \fr{C}{(\cos \theta) ^{\fr{4-\alpha}{\alpha}}} \|\na u(se^{i\theta})\|_{L^2}^2 \|\na q(se^{i\theta})\|_{L^2}^{\fr{4}{\alpha}},
\end{align}
\be 
\left|(B(u,u), Au)_{L^2} \right|
\le C\|\na u\|_{L^2}^{\fr{3}{2}} \|Au\|_{L^2}^{\fr{3}{2}}
\le \fr{\cos \theta}{8} \|Au(se^{i\theta})\|_{L^2}^2 + \fr{C}{(\cos \theta)^3}\|\na u(se^{i\theta})\|_{L^2}^6,
\ee 
\be 
\left|(\PP(qRq), Au)_{L^2} \right| 
\le \fr{\cos \theta}{8} \|Au(se^{i\theta})\|_{L^2}^2 + \fr{C}{\cos \theta} \|\na q(se^{i\theta})\|_{L^2}^4,
\ee and
\be \la{backun2}
\left|(f, Au)_{L^2} \right|
\le \fr{\cos \theta}{8} \|Au(se^{i\theta})\|_{L^2}^2 + \fr{C}{\cos \theta} \|f\|_{L^2}^2
\ee using the H\"older, Gagliardo-Nirenberg, and Young inequalities, and continuous Sobolev embeddings. Combining \eqref{backun1}--\eqref{backun2}, we obtain the differential inequality
\beg{align}
&\fr{d}{ds} \left[\|\na q(se^{i\theta})\|_{L^2}^2 + \|A^{\fr{1}{2}} u(se^{i\theta})\|_{L^2}^2 \right] 
\le \fr{C}{(\cos \theta) ^{\fr{4-\alpha}{\alpha}}} \|\na u(se^{i\theta})\|_{L^2}^2 \|\na q(se^{i\theta})\|_{L^2}^{\fr{4}{\alpha}} \nonumber
\\&\quad\quad+  \fr{C}{(\cos \theta)^3}\|\na u(se^{i\theta})\|_{L^2}^6 +  \fr{C}{\cos \theta} \|\na q(se^{i\theta})\|_{L^2}^4 + \fr{C}{\cos \theta} \|f\|_{L^2}^2,
\end{align} from which we conclude that 
\be 
\|\na q(se^{i\theta})\|_{L^2}^2 + \|A^{\fr{1}{2}} u(se^{i\theta})\|_{L^2}^2
\le 2 \left[\|\na q_0\|_{L^2}^2 + \|\na u_0\|_{L^2}^2  + 1\right]
\ee provided that 
\be \la{backun4}
s \left(\fr{C}{\cos \theta} \|f\|_{L^2}^2 + \fr{C}{\cos \theta} + \fr{C}{(\cos \theta)^3} + \fr{C}{(\cos \theta) ^{\fr{4-\alpha}{\alpha}}}  \right) \le \Gamma_0
\ee Here $C$ is a positive universal constant,  and $\Gamma_0$ is a positive constant depending only on $\|\omega_0\|_{\mathcal{V}}$. Therefore $(q,u)$ is locally time analytic  on the region $\mathcal{R}$ consisting of complex times $t = se^{i\theta}$ obeying \eqref{backun4}. Due to the uniform-in-time boundedness of $(q,u)$ in the norm of $\mathcal{V}$, the time analyticity becomes global. 

\textbf{Step 2. Backward uniqueness.} Since $\mathcal{S}(T)w_1^0 = \mathcal{S}(T) w_2^0$, then $\mathcal{S}(t)w_1^0 = \mathcal{S}(t) w_2^0$ for all times $t \ge T$ due to the uniqueness of solutions in $\mathcal{V}$. From the time analyticity property derived in Step 1, we conclude that $\mathcal{S}(T)w_1^0$ and $\mathcal{S}(T) w_2^0$ coincides for all positive times. Consequently $w_1^0 = w_2^0$, ending the proof of Proposition \ref{backun}.

\subsection{Decay of Volume Elements.}

Let $\phi$ be a smooth function defined on open set $\Omega \subset \RR^N$, $N \ge 1$, and taking values in $\mathcal{V}$. Let $\Sigma_t = S(t) \phi(\Omega)$. The volume element in $\Sigma_t$ is given by 
\[
\left|\fr{\pa}{\pa \alpha_1} (S(t) \phi(\alpha)) \wedge \dots \wedge \fr{\pa}{\pa \alpha_N} (S(t) \phi(\alpha)) \right| d\alpha_1 \dots d\alpha_N,
\] where  $d\alpha_1 \dots d\alpha_N$ is the volume element in $\RR^N$. The functions 
\be 
\omega_i = \fr{\pa}{\pa \alpha_i} S(t) \phi(\alpha), \hspace{1cm} i = 1, \dots, N
\ee solve the linearized system
\be \la{vol1}
\pa_t (q,u) + \mathcal{A} (q,u) + L(\bar{\omega}) (q,u) = 0
\ee along $(\bar{q}(t), \bar{u}(t)) := \bar{\omega} (t) = \mathcal{S}(t) \phi(\alpha)$, where 
\be 
\mathcal{A}(q,u) = (\l^{\alpha}q, Au)
\ee and
\be 
L(\bar{\omega}) (q,u) = \left(\bar{u} \cdot \na q + u \cdot \na \bar{q}, B(u, \bar{u}) + B(\bar{u},u) + \PP(qR\bar{q} + \bar{q}Rq) \right). 
\ee We address the time evolution of the volume element of the $N$-dimensional surface $\phi(\Omega)$ transported by $\mathcal{S}(t)$. For that purpose, we consider the norm
\be 
V_N(t) = \|\omega_1(t) \wedge \dots \wedge \omega_N(t)\|_{\Lambda^N \mathcal{H}}
\ee where $\omega_1, \dots, \omega_N$ solves \eqref{vol1} along some $\bar{\omega}(t)= \mathcal{S}(t) \bar{\omega}_0$, and $\Lambda^N \mathcal{H}$ is the $N$-th exterior product of $\mathcal{H}$ with the following scalar product 
\[
(\omega_1 \wedge \dots \wedge \omega_N ; y_1 \wedge \dots \wedge y_N)_{\Lambda^N \mathcal{H}} = det (w_i, y_j)_{\mathcal{H}}.
\]

\beg{prop} \la{vol} There exists a positive time $t_0$ depending only on $\|\bar{\omega}_0\|_{\mathcal{V}}$, a positive integer $N_0$ depending only on $\|f\|_{L^2}$, and a positive constant $c$ depending on $\alpha$, such that the following decaying estimate
\be \la{vol6}
V_N(t) \le V_N(0) e^{-cN^{1 + \fr{\alpha}{2}}t}
\ee holds for all $t \ge t_0$ and $N \ge N_0$. 
\end{prop}

\noindent \textbf{Proof of Proposition \ref{vol}.} We have 
\be 
\pa_t (\omega_1 \wedge \dots \wedge \omega_N) + (\mathcal{A}+ L(\bar{\omega}))_{N} (\omega_1 \wedge \dots \wedge \omega_N) = 0,
\ee where 
\be 
(\mathcal{A} + L(\bar{\omega}))_{N} = (\mathcal{A} + L(\bar{\omega})) \wedge \mathcal{I} \wedge \dots \wedge \mathcal{I} + \dots + \mathcal{I} \wedge \dots \wedge \mathcal{I} \wedge (\mathcal{A} + L(\bar{\omega})),
\ee $\mathcal{I}$ being the identity operator. Consequently, we obtain the evolution equation
\be 
\fr{d}{dt} V_N + \mathrm{Trace}((\mathcal{A} + L(\bar{\omega}))Q_N) V_N = 0
\ee where $Q_N$ is the orthogonal projection in $\mathcal{H}$ onto the space spanned by $\omega_1, \dots, \omega_N$. By Gronwall's inequality, we obtain 
\be \la{vol4}
V_N(t) \le V_N(0) \exp \left\{-\int_{0}^{t} \mathrm{Trace}((\mathcal{A} + L(\bar{\omega}))Q_N)  ds \right\}.
\ee For each time $t > 0$, we let $\phi_i =  (r_i, v_i)$, $i = 1, \dots, N$, be an orthonormal family of functions in $\mathcal{H}$ spanning the linear span of $\omega_1, \dots, \omega_N$. Then, we have 
\be 
\mathrm{Trace}((\mathcal{A} + L(\bar{\omega}))Q_N)  = \sum\limits_{i=1}^{N} (\mathcal{A} \phi_i, \phi_i)_{L^2} +  \sum\limits_{i=1}^{N} (L(\bar{\omega}) \phi_i, \phi_i)_{L^2}.
\ee In view of Lemma \ref{voll}, we obtain the lower estimate
\be 
\sum\limits_{i=1}^{N} (\mathcal{A} \phi_i, \phi_i)_{L^2}
= \sum\limits_{i=1}^{N} \left[(\l^{\alpha} r_i, r_i)_{L^2} + (Av_i, v_i)_{L^2} \right]
\ge \mu_1 + \dots + \mu_N \ge CN^{1 + \fr{\alpha}{2}}
\ee where $\mu_1, \dots, \mu_N$ are the first $N$ eigenvalues of $\mathcal{A}$. Now we show that 
\be \la{vol3}
|\mathrm{Trace} (L(\bar{\omega})Q_N) |
\le C\left(\|\bar{\omega}\|_{{\mathcal{V}}}^2 + \|\bar \omega\|_{\mathcal{V}}^{\fr{4}{\alpha}} + 1 \right)N + \fr{1}{2} \mathrm{Trace}(\mathcal{A}Q_N). 
\ee Indeed, the trace of $L(\bar{\omega})Q_N$ can be estimated as follows,
\beg{align} \la{tra1}
&|\mathrm{Trace} (L(\bar{\omega})Q_N) |
= \left|\sum\limits_{i=1}^{N} (L(\bar{\omega}) \phi_i, \phi_i)_{L^2} \right| \nonumber
\\&\le \sum\limits_{i=1}^{N} \left|(\bar{u} \cdot \na r_i + v_i \cdot \na \bar{q}, r_i)_{L^2} \right|
+ \sum\limits_{i=1}^{N} \left|(B(v_i, \bar{u}) + B(\bar{u}, v_i) + \PP(r_i R\bar{q} + \bar{q} Rr_i), v_i)_{L^2} \right| \nonumber
\\&\le C\sum\limits_{i=1}^{N} \left[\|v_i\|_{L^{\fr{4}{\alpha}}} \|\na \bar{q}\|_{L^2} \|r_i\|_{L^{\fr{4}{2-\alpha}}} + \|v_i\|_{L^4} \|\na \bar{u}\|_{L^2} \|v_i\|_{L^4} + \|r_i\|_{L^2} \|\bar{q}\|_{L^4} \|v_i\|_{L^4} \right].
\end{align} Here the boundedness of the Riesz transform on $L^p$ spaces is exploited. In view of the continuous embedding of $\mathcal{D}(\l^{\fr{\alpha}{2}})$ into $L^{\fr{4}{2-\alpha}}$ and the Gagliardo-Nirenberg interpolation inequalities, we bound
\beg{align}  \la{tra2}
&\|v_i\|_{L^{\fr{4}{\alpha}}} \|\na \bar{q}\|_{L^2} \|r_i\|_{L^{\fr{4}{2-\alpha}}}
\le C\|v_i\|_{L^2}^{\fr{\alpha}{2}} \|\na v_i\|_{L^2}^{\fr{2-\alpha}{2}} \|\na \bar{q}\|_{L^2} \|\l^{\fr{\alpha}{2}} r_i\|_{L^2} \nonumber
\\&\quad\quad\le \fr{1}{4} \left[\|\l^{\fr{\alpha}{2}} r_i\|_{L^2}^2  + \|A^{\fr{1}{2}} v_i\|_{L^2}^2 \right] + C\|\na \bar{q}\|_{L^2}^{\fr{4}{\alpha}} \|v_i\|_{L^2}^2.
\end{align} Applications of Ladyzhenskaya's interpolation inequality and Young's inequality give
\beg{align} \la{tra3}
&\|v_i\|_{L^4} \|\na \bar{u}\|_{L^2} \|v_i\|_{L^4} + \|r_i\|_{L^2} \|\bar{q}\|_{L^4} \|v_i\|_{L^4} \nonumber
\\&\quad\quad\le \fr{1}{4}\|A^{\fr{1}{2}} v_i\|_{L^2}^2  + C\left(\|\na \bar{u}\|_{L^2}^2 + \|\bar{q}\|_{L^4}^2 + 1 \right) \left(\|v_i\|_{L^2}^2 + \|r_i\|_{L^2}^2\right).
\end{align} Putting \eqref{tra1}--\eqref{tra3} together, and using the normalization $\|\phi_i\|_{L^2}^2 = \|r_i\|_{L^2}^2 + \|v_i\|_{L^2}^2 = 1$, we obtain the desired estimate \eqref{vol3}, from which we infer that
\beg{align}
\int_{0}^{t} \mathrm{Trace}((\mathcal{A} + L(\bar{\omega}))Q_N)(s) ds
&\ge \fr{1}{2} \int_{0}^{t} \mathrm{Trace}(\mathcal{A}Q_N)ds  - CN \int_{0}^{t} \left(\|\bar \omega\|_{\mathcal{V}}^2 + \|\bar \omega\|_{\mathcal{V}}^{\fr{4}{\alpha}} + 1 \right) ds \nonumber
\\&\ge CN t \left(N^{\fr{\alpha}{2}} - \fr{1}{t} \int_{0}^{t} \left(\|\bar \omega\|_{\mathcal{V}}^2 + \|\bar \omega\|_{\mathcal{V} }^{\fr{4}{\alpha}} + 1 \right) ds \right).
\end{align} We apply Proposition \ref{absball} and obtain the existence of a time $t_0$ depending only on $\|\bar{\omega}_0\|_{\mathcal{V}}$ and a radius $\rho_f$ depending only on $\|f\|_{L^2}$ such that 
\be 
\|\bar \omega\|_{\mathcal{V}}^2 + \|\bar \omega\|_{\mathcal{V}}^{\fr{4}{\alpha}} + 1\le \rho_f
\ee for all $t \ge t_0$. Consequently, it holds that 
\be \la{vol5}
\int_{0}^{t} \mathrm{Trace}((\mathcal{A} + L(\bar{\omega}))Q_N)(s) ds \ge CN^{1+\fr{\alpha}{2}}t
\ee for all $t \ge t_0$, provided that $N^{\fr{\alpha}{2}} \ge 2 \rho_f$. Putting \eqref{vol4} and \eqref{vol5} together, we obtain \eqref{vol6}, completing the proof of Proposition \ref{vol}. 

\subsection{Existence of a Finite Dimensional Global Attractor.}

As a consequence of the connectedness and compactness of the absorbing ball $\mathcal{B}_{\rho}$ in $\mathcal{H}$, the continuity and injectivity of the solution map $\mathcal{S}(t)$, and the exponential time decay of volume elements, we conclude that the model \eqref{system} has a finite-dimensional global attractor. We refer the reader to  \cite[Chapter 14]{cfbook} for a detailed proof of the analogous result for the two-dimensional forced Navier-Stokes equations. 

\beg{rem} \la{potzero} We note that Propositions \ref{absball}, \ref{conts}, \ref{backun}, and \ref{vol} hold in the presence of a time independent potential $\Phi$. In this latter case, the radius of the absorbing ball depends on the size of the body forces $f$ in the fluid and the potential $\Phi$. This gives Theorem \ref{att}. 
\end{rem}

\section{Regularity of the Global Attractor for $\alpha=1$: Proof of Theorem \ref{attp}} \la{s6}

In this section, we address the forced electroconvection model \eqref{system}, where $\alpha$ is taken to be 1. We prove the existence of an absorbing ball, compact in the strong norm of $\mathcal{V}$. 

\beg{prop} \la{cabsball} Suppose $\alpha=1$ and $f \in \mathcal{D}(A^{\fr{1}{2}})$. Then there exists a radius $\tilde \rho > 0$ depending only on $f$ and some universal constants such that for each $\omega_0 = (q_0, u_0) \in \mathcal{V}$, there exists a time $\tilde T_0$ depending only on $\|\na q_0\|_{L^2}$ and $\|\na u_0\|_{L^2}$ and universal constants such that 
\be \la{cabsball1}
\mathcal{S}(t) \omega_0 \in  \mathcal{B}_{\tilde \rho} := \left\{(q, u) \in \mathcal{V}: \|\l^{\fr{3}{2}} q\|_{L^2} + \|\Delta u\|_{L^2} \le \tilde \rho \right\}
\ee for all $t \ge \tilde T_0$.
\end{prop}

\noindent \textbf{Proof.} There exists a radius $R>0$ depending only on the body forces $f$, such that for any $\omega_0 \in \mathcal{V}$, there is a positive time $t_0$ depending only on the $\mathcal{V}$-norm of $\omega_0 = (q_0, u_0)$ such that the solution $(q,u)$ of \eqref{system}, with initial datum $\omega_0$, obeys
\be \la{ft0}
\|q(t)\|_{L^{\infty}} + \|\na q(t)\|_{L^2} + \|\Delta u(t)\|_{L^2} + \int_{t}^{t+1} \|\na \Delta u(s)\|_{L^2}^2 ds \le R
\ee at any time $t \ge t_0$, a fact that follows from the proof of Proposition \ref{absball}. Moreover, there exists a positive time $t_1 \ge t_0$ depending only on $\|\omega_0\|_{\mathcal{V}}$ such that for $q(t_1) \in \mathcal{D}(\l^{\fr{3}{2}})$  with a size dependency only on the forces $f$, and such that $q \in L^{\infty}(t_1,T; \mathcal{D}(\l^{\fr{3}{2}}) \cap L^2(t_1, T; \mathcal{D}(\l^2))$ for any $T \ge t_1$. We seek bounds for the charge density in those latter Lebesgue spaces, independent of the initial datum but depending only on the size of the body forces. 
 
The $L^2$ norm of $\l^{\fr{3}{2}}q$ evolves according to the energy equality
\beg{align} \la{jji}
\fr{1}{2} \fr{d}{dt} \|\l^{\fr{3}{2}} q\|_{L^2}^2
+ \|\Delta q\|_{L^2}^2
= -\int_{\Omega} \na \l^{\fr{1}{2}} (u \cdot \na q) \cdot \na \l^{\fr{1}{2}} q dx,
\end{align} which is equivalent to 
\beg{align} \la{ft}
\fr{1}{2} \fr{d}{dt} \|\l^{\fr{3}{2}} q\|_{L^2}^2
+ \|\Delta q\|_{L^2}^2
&= -\int_{\Omega} \left[\na \l^{\fr{1}{2}} (u \cdot \na q) - \l^{\fr{1}{2}} \na (u \cdot \na q) \right] \cdot \na \l^{\fr{1}{2}} q dx \nonumber
\\&\quad\quad -\int_{\Omega} \l^{\fr{1}{2}} \na (u \cdot \na q) \cdot \na \l^{\fr{1}{2}} q dx.
\end{align} We set 
\be \la{newt33}
A := \int_{\Omega} \left[\na \l^{\fr{1}{2}} (u \cdot \na q) - \l^{\fr{1}{2}} \na (u \cdot \na q) \right] \cdot \na \l^{\fr{1}{2}} q dx
\ee and 
\be \la{newt3}
B := \int_{\Omega} \l^{\fr{1}{2}} \na (u \cdot \na q) \cdot \na \l^{\fr{1}{2}} q dx.
\ee We decompose $B$ into a sum of three spatial integrals $B_1, B_2$ and $B_3$, where
\be
B_1 = \int_{\Omega} (\l^{\fr{1}{2}} (u \cdot \na \na q) - u \cdot \l^{\fr{1}{2}} \na \na q) \cdot \na \l^{\fr{1}{2}}q dx ,
\ee
\be 
B_2 = \int_{\Omega} u \cdot (\l^{\fr{1}{2}} \na \na q - \na \na \l^{\fr{1}{2}}q) \cdot \na \l^{\fr{1}{2}}q dx,
\ee 
\be 
B_3 = \int_{\Omega} (\l^{\fr{1}{2}} (\na u \cdot \na q) - \na u \cdot \l^{\fr{1}{2}} \na q )\cdot \na \l^{\fr{1}{2}}q dx,
\ee 
\be 
B_4 = \int_{\Omega} (\na u \cdot \l^{\fr{1}{2}} \na q  - \na u \cdot \na \l^{\fr{1}{2}}  q )\cdot \na \l^{\fr{1}{2}}q dx, 
\ee and
\be 
B_5 = \int_{\Omega} (\na u \cdot \na \l^{\fr{1}{2}} q) \cdot \na \l^{\fr{1}{2}}q dx.
\ee
In view of Corollary \ref{EXIU} with $\alpha=1$ and $p=9$, the embedding of $\mathcal{D}(\l^{\fr{1}{2}})$ in $L^4$, and $L^p$ interpolation inequalities,  we have
\beg{align} 
|A| 
&\le \|[\na, \l^{\fr{1}{2}}] (u \cdot \na q)\|_{L^{\fr{4}{3}}} \|\l^{\fr{3}{2}}q\|_{L^4}
\le C\|\Delta u\|_{L^2} \|\na q\|_{L^9} \|\Delta q\|_{L^2} \nonumber
\\&\le C\|\Delta u\|_{L^2}\|\na q\|_{L^2}^{\fr{1}{8}}\|\na q\|_{L^{18}}^{\fr{7}{8}} \|\Delta q\|_{L^2}
\le C\|\Delta u\|_{L^2}\|\na q\|_{L^2}^{\fr{1}{8}}\|\Delta q\|_{L^2}^{1+\fr{7}{8}}\nonumber
\\&\le \fr{1}{8} \|\Delta q\|_{L^2}^2 + C\|\Delta u\|_{L^2}^{16}\|\na q\|_{L^2}^{2}.
\end{align} We note that the very regular dissipation $\|\Delta q\|_{L^2}$ is exploited for the sake of interpolation. 
We estimate 
\beg{align}
|B_1|
&\le \|[\l^{\fr{1}{2}}, u] \na \na q \|_{L^2}\|\l^{\fr{3}{2}}q\|_{L^2}
\le C\|\na u\|_{L^{\infty}} \|\na \na q\|_{L^2} \|\l^{\fr{3}{2}}q\|_{L^2} \nonumber
\\&\le C\|\Delta u\|_{L^2}^{\fr{1}{2}} \|\na \Delta u\|_{L^2}^{\fr{1}{2}} \|\na q\|_{L^2}^{\fr{1}{2}} \|\Delta q\|_{L^2}^{\fr{3}{2}}
\le \fr{1}{8} \|\Delta q\|_{L^2}^2 + C\|\Delta u\|_{L^2}^{2} \|\na \Delta u\|_{L^2}^2\|\na q\|_{L^2}^{2}
\end{align} by appealing to Proposition \ref{coes1} with $s = \fr{1}{2}$. 
By making use of Corollary \ref{exiu55} with $\alpha = 1$ and $p=9$, we bound
\beg{align}
|B_2| &\le \|u \cdot [\l^{\fr{1}{2}}, \na \na]q \|_{L^{\fr{4}{3}}} \|\l^{\fr{3}{2}}q\|_{L^4}
\le C\|\Delta u\|_{L^2} \|\na q\|_{L^9} \|\Delta q\|_{L^2} \nonumber
\\&\le C\|\Delta u\|_{L^2} \|\na q\|_{L^2}^{\fr{1}{8}} \|\Delta q\|_{L^2}^{\fr{15}{8}}
\le  \fr{1}{8} \|\Delta q\|_{L^2}^2 + C\|\Delta u\|_{L^2}^{16}\|\na q\|_{L^2}^{2}.
\end{align} 
Another application of Proposition \ref{coes1} yields
\be 
|B_3| \le \|[\l^{\fr{1}{2}}, \na u] \na q \|_{L^2}\|\l^{\fr{3}{2}}q\|_{L^2}
\le C\|\na \Delta u\|_{L^{2}} \|\na q\|_{L^2} \|\l^{\fr{3}{2}}q\|_{L^2}
\le \fr{1}{8} \|\Delta q\|_{L^2}^2 + C\|\na \Delta u\|_{L^2}^{2}\|\na q\|_{L^2}^{2}
\ee after making use of standard Sobolev embedding. 
By Corollary \ref{exiu55i} with $\alpha=1$ and $p=9$, we have
\be 
|B_4| \le \|\na u \cdot [\l^{\fr{1}{2}}, \na] q\|_{L^{\fr{4}{3}}}\|\l^{\fr{3}{2}}q\|_{L^4}
\le C\|\Delta u\|_{L^2} \|\na q\|_{L^9}  \|\Delta q\|_{L^2}
\le \fr{1}{8} \|\Delta q\|_{L^2}^2 + C\|\Delta u\|_{L^2}^{16}\|\na q\|_{L^2}^{2}.
\ee Finally, we estimate 
\be \la{yi}
|B_5| \le C\|\na u\|_{L^4} \|\l^{\fr{3}{2}}q\|_{L^4} \|\l^{\fr{3}{2}}q\|_{L^2}
\le \fr{1}{8} \|\Delta q\|_{L^2}^2 + C\|\Delta u\|_{L^2}^4 \|\na q\|_{L^2}^2
\ee via a direct application of H\"older's inequality. 
Putting \eqref{ft}--\eqref{yi} together, we conclude that 
\be \la{newt1}
\fr{d}{dt} \|\l^{\fr{3}{2}} q\|_{L^2}^2
+ \|\Delta q\|_{L^2}^2
\le C\left(\|\Delta u\|_{L^2}^{16} + \|\Delta u\|_{L^2}^4 + \|\Delta u\|_{L^2}^{2} \|\na \Delta u\|_{L^2}^2\right) \|\na q\|_{L^2}^2.
\ee In view of \eqref{ft0}, the property \eqref{cabsball1} holds, ending the proof of Proposition \ref{cabsball}.

Now we address the smoothness of the attractor $\tilde{X}$:

\beg{prop} \la{smooglo} Let $\omega_0 \in \mathcal{V}$ and $f \in \mathcal{D}(A^{\fr{k+1}{2}})$. Suppose there exists a time $t_k^0 > 0$ depending on $\|\omega_0\|_{\mathcal{V}}$, and a radius $R_k>0$ depending only on $\|f\|_{H^{k}}$ such that the estimate 
\be \la{smooglo1}
\|\l^{k} q(t)\|_{L^2} + \|A^{\fr{k}{2} + \fr{1}{2}} u(t)\|_{L^2} + \int_{t}^{t+1} \|A^{\fr{k}{2} + 1} u(s)\|_{L^2}^2 ds \le R_k
\ee holds for all $t \ge t_k^0$. Moreover, suppose there is a time $t_k \ge t_k^0$ such that the $L^2$ norm of $\l^{k+\fr{1}{2}} q(t_k)$ is bounded by some constant depending only on  $\|f\|_{H^{k}}$ and $\|\omega_0\|_{\mathcal{V}}$. Then there exists a time $\tilde t_{k+1}^0 > 0$ depending on $\|\omega_0\|_{\mathcal{V}}$, and a radius $\tilde R_{k+1}>0$ depending only on $\|f\|_{H^{k+1}}$ such that the estimate 
\be \la{smooglo2} 
\|\l^{k+1} q(t)\|_{L^2} + \|A^{\fr{k}{2} + 1} u(t)\|_{L^2} + \int_{t}^{t+1} \|A^{\fr{k+3}{2}} u(s)\|_{L^2}^2 ds \le \tilde R_{k+1}
\ee holds for all $t \ge \tilde t_k^0$. Moreover, there is a time $t_{k+1} \ge \tilde t_k^0$ such that the $L^2$ norm of $\l^{k+\fr{3}{2}} q(t_{k+1})$ is bounded by some constant depending only on  $\|f\|_{H^{k+1}}$ and $\|\omega_0\|_{\mathcal{V}}$.
\end{prop}

\noindent \textbf{Proof.} The $L^2$ norm of $\l^{k+\fr{1}{2}} q$ evolves according to the energy equality
\be \la{newt2}
\fr{1}{2} \fr{d}{dt} \|\l^{k+\fr{1}{2}} q\|_{L^2}^2
+ \|\l^{k+1} q\|_{L^2}^2
= -\int_{\Omega} \l^{k+\fr{1}{2}} (u \cdot \na q) \l^{k+\fr{1}{2}}q dx.
\ee We let
\be 
\mathcal{N} := \int_{\Omega} \l^{k+\fr{1}{2}} (u \cdot \na q) \l^{k+\fr{1}{2}}q dx
\ee and we distinguish three different cases: $k = 1$, $k \ge 2$ even, and $k \ge 3$ odd. 

If $k=1$, the equality \eqref{newt2} reduces to 
\be
\fr{d}{dt} \|\l^{\fr{3}{2}} q\|_{L^2}^2
+ \|\Delta q\|_{L^2}^2
\le C\left(\|\Delta u\|_{L^2}^{16} + \|\Delta u\|_{L^2}^4 + \|\Delta u\|_{L^2}^{2} \|\na \Delta u\|_{L^2}^2\right) \|\na q\|_{L^2}^2,
\ee as shown in \eqref{newt1}.

Now we fix an even integer $k \ge 2$ and decompose the nonlinear term $\mathcal{N}$ as a sum 
\be 
\mathcal{N} = \mathcal{N}_1 + \mathcal{N}_2,
\ee where 
\be 
\mathcal{N}_1 = \int_{\Omega} \l^{\fr{1}{2}} \left[\l^k(u \cdot \na q) - u \cdot \na \l^kq \right] \l^{k + \fr{1}{2}}q dx
\ee and 
\be 
\mathcal{N}_2 = \int_{\Omega} \l^{\fr{1}{2}} \left[ u \cdot \na \l^k q \right] \l^{k + \fr{1}{2} }q dx.
\ee
We estimate the term $\mathcal{N}_1$
\beg{align}
|\mathcal{N}_1| 
&\le \|\l^k (u \cdot \na q) - u \cdot \na \l^k q  \|_{L^2} \|\l^{k+1}q\|_{L^2} \nonumber
\le C\|u\|_{H^{k+1}} \|\l^{k+\fr{1}{2}}q\|_{L^2} \|\l^{k+1}q\|_{L^2} \nonumber
\\&\le C\| u\|_{H^{k+1}} \|\l^{k}q\|_{L^2}^{\fr{1}{2}} \|\l^{k+1}q\|_{L^2}^{\fr{3}{2}} 
\le \fr{1}{16} \|\l^{k+1}q\|_{L^2}^2 + C\| u\|_{H^{k+1}}^4 \|\l^{k}q\|_{L^2}^2
\end{align}
by appealing to Proposition \ref{pr3}. Due to the incompressibility of the fluid, we can decompose $\mathcal{N}_2$ as the sum
\be 
\mathcal{N}_2  = \mathcal{N}_{2,1} + \mathcal{N}_{2,2}
\ee where
\be 
\mathcal{N}_{2,1} = \int_{\Omega} \left[\l^{\fr{1}{2}} \left(u \cdot \na \l^k q \right) - u \cdot \l^{\fr{1}{2}} \na \l^k q \right] \l^{k+\fr{1}{2}}q dx
\ee and 
\be 
\mathcal{N}_{2,2} = \int_{\Omega} u \cdot \left[\l^{\fr{1}{2}} \na \l^k q - \na \l^{\fr{1}{2}} \l^k q \right]  \l^{k+\fr{1}{2}}q dx. 
\ee  In view of the commutator estimate \cite[Theorem (2.2)]{CN2}, the following pointwise commutator estimate
\be 
|[\l^{\fr{1}{2}}, \na]\tilde{q}| \le Cd(x)^{-\fr{1}{2} - 1-\fr{2}{9}} \|\tilde{q}\|_{L^9}
\ee holds for any $\tilde{q} \in C_0^{\infty}(\Omega)$, thus the operator $u \cdot [\na, \l^{\fr{1}{2}}]$ extends from $C_0^{\infty}$ to $L^9$ such that the estimate
\be 
\|u \cdot [\l^{\fr{1}{2}}, \na] \tilde{q}\|_{L^{\fr{4}{3}}} 
\le C\|u d(\cdot)^{-1}\|_{L^{180}} \|d(\cdot)^{- \fr{1}{2} - \fr{2}{9}} \|_{L^{\fr{90}{67}}} \|\tilde{q}\|_{L^9}
\le C\|\na u\|_{L^{180}}  \|\tilde{q}\|_{L^9}
\le C\|\Delta u\|_{L^2} \|\tilde{q}\|_{L^9}
\ee holds for any $\tilde{q} \in L^9$ due to Rellich's inequality. As a consequence, the term $N_{2,2}$ can be bounded as follows,
\beg{align}
|N_{2,2}|
&\le C\|u \cdot [\l^{\fr{1}{2}}, \na] \l^k{q}\|_{L^{\fr{4}{3}}} \|\l^{k+\fr{1}{2}} q\|_{L^4}
\le C\|\Delta u\|_{L^2} \|\l^k q\|_{L^9} \|\l^{k+1} q\|_{L^2} \nonumber
\\&\le C\|\Delta u\|_{L^2} \|\l^{k} q\|_{L^2}^{\fr{1}{8}}\|\l^{k}q\|_{L^{18}}^{\fr{7}{8}}\|\l^{k+1} q\|_{L^2}
\le C\|\Delta u\|_{L^2} \|\l^{k} q\|_{L^2}^{\fr{1}{8}}\|\l^{k+1} q\|_{L^2}^{\fr{15}{8}} \nonumber
\\&\le \fr{1}{16}\|\l^{k+1}q\|_{L^2}^2  + C\|\Delta u\|_{L^2}^{16}\|\l^k q\|_{L^2}^2.
\end{align}

We estimate $\mathcal{N}_{2,1}$
\beg{align} 
|\mathcal{N}_{2,1}| 
&\le \|[\l^{\fr{1}{2}}, u]\na \l^k q\|_{L^2} \|\l^{k+\fr{1}{2}}q\|_{L^2} \nonumber
\\&\le C\|\na u\|_{L^{\infty}} \|\na \l^k q\|_{L^2} \|\l^{k+\fr{1}{2}}q\|_{L^2} \nonumber
\\&\le C\|\na \Delta u\|_{L^2} \|\l^{k+1}q\|_{L^2}^{\fr{3}{2}} \|\l^k q\|_{L^2}^{\fr{1}{2}} \nonumber
\\&\le \fr{1}{16} \|\l^{k+1}q\|_{L^2}^2 + C\|\na \Delta u\|_{L^2}^4 \|\l^k q\|_{L^2}^2
\end{align} by making use of Proposition \ref{coes1}. Therefore, the energy equality \eqref{newt2} yields 
\be 
\fr{d}{dt} \|\l^{k+\fr{1}{2}}q\|_{L^2}^2 + \|\l^{k+1}q\|_{L^2}^2
\le C\left(\|\Delta u\|_{L^2}^{16} + \|A^{\fr{k+1}{2}} u\|_{L^2}^4  \right)\|\l^{k}q\|_{L^2}^2.
\ee 
As for the last case, we fix an odd integer $k \ge 3$, rewrite $\mathcal{N}$ as 
\be 
\mathcal{N} = \int_{\Omega} \na \l^{\fr{1}{2}} \l^{k-1} (u \cdot \na q) \cdot \na \l^{k-\fr{1}{2}}q dx,
\ee and decompose it into the sum of three terms 
\be 
\mathcal{N} 
= \tilde{\mathcal{N}}_1
+ \tilde{\mathcal{N}}_2
\ee where
\be 
\tilde{\mathcal{N}}_1 = \int_{\Omega} \na \l^{\fr{1}{2}} \left[\l^{k-1} (u \cdot \na q) - u \cdot \na \l^{k-1}q \right] \cdot \na \l^{k-\fr{1}{2}}q dx,
\ee and
\be 
\tilde{\mathcal{N}}_2 = \int_{\Omega} \na \l^{\fr{1}{2}} (u \cdot \na \l^{k-1}q )  \cdot \na \l^{k-\fr{1}{2}}q dx.
\ee 
The term $\tilde{\mathcal{N}}_2$ has the same structure as the nonlinear term on the right-hand side of \eqref{jji} with $q$ replaced by $\l^{k-1} q$, thus it bounds as 
\beg{align}
|\tilde{\mathcal{N}}_2| 
&\le \fr{1}{16} \|\l^{k+1}q\|_{L^2}^2
+ C\left(\|\Delta u\|_{L^2}^{16} + \|\Delta u\|_{L^2}^4 + \|\Delta u\|_{L^2}^2 \|\na \Delta u\|_{L^2}^2 \right) \|\l^kq\|_{L^2}^2.
\end{align}  As for the term $\tilde{\mathcal{N}}_1$, we can rewrite it as 
\be 
\tilde{\mathcal{N}}_1 = \int_{\Omega} \na \left[\l^{k-1} (u \cdot \na q) - u \cdot \na \l^{k-1}q \right] \cdot \na \l^k q dx
\ee after integrating by parts several times, and then we decompose it as a sum
\be 
\tilde{\mathcal{N}}_1 = \tilde{\mathcal{N}}_{1,1} + \tilde{\mathcal{N}}_{1,2}
\ee where
\be 
\tilde{\mathcal{N}}_{1,1} =  \int_{\Omega} \left[\l^{k-1}\na (u \cdot \na q) - u \cdot \na \l^{k-1} \na q \right] \cdot \na \l^k q dx
\ee and 
\be 
\tilde{\mathcal{N}}_{1,2} = \int_{\Omega} \left[ u \cdot \na \l^{k-1} \na q - \na (u \cdot \na \l^{k-1}q) \right] \cdot \na \l^k q dx.
\ee The decomposition above uses the fact that $\na$ and $\l^{k-1}$ commutes when $k$ is odd. By expanding $\na (u \cdot \na \l^{k-1} q)$, the term $\tilde{\mathcal{N}}_{1,2}$ reduces to 
\be
\tilde{\mathcal{N}}_{1,2} = - \int_{\Omega} \left[\na u \cdot \na \l^{k-1} q \right] \cdot \na \l^k q dx
\ee and is bounded by 
\be 
|\tilde{\mathcal{N}}_{1,2}| \le C\|\na u\|_{L^{\infty}} \|\l^k q\|_{L^2} \|\l^{k+1}q\|_{L^2}
\le \fr{1}{16} \|\l^{k+1}q\|_{L^2}^2 + C\|\na \Delta u\|_{L^2}^2 \|\l^k q\|_{L^2}^2.
\ee In view of the commutator estimate \eqref{pr32}, we estimate 
\beg{align}
|\tilde{\mathcal{N}}_{1,1}|
&\le \|[\na \l^{k-1}, u \cdot \na ]q\|_{L^2} \|\na \l^k q\|_{L^2} \nonumber
\\&\le C\| u\|_{H^{k+1}} \|\l^{k+\fr{1}{2}}q\|_{L^2} \|\l^{k+1}q\|_{L^2} \nonumber
\\&\le \fr{1}{16} \|\l^{k+1}q\|_{L^2}^2 + C\|A^{\fr{k+1}{2}}u\|_{L^2}^4 \|\l^k q\|_{L^2}^2,
\end{align} where the last bound follows from interpolation and Young's inequality. Therefore, the energy inequality 
\be 
\fr{d}{dt} \|\l^{k+\fr{1}{2}} q\|_{L^2}^2 + \|\l^{k+1} q\|_{L^2}^2 
\le C\left[\|A^{\fr{k+1}{2}} u\|_{L^2}^4 + \|\na \Delta u\|_{L^2}^2 + \|\Delta u\|_{L^2}^{16} \right]\|\l^k q\|_{L^2}^2
\ee holds when $k \ge 3$ is odd. In all three cases, and due to the assumption \eqref{smooglo1} and the Gronwall Lemma \ref{gron}, we obtain a time $T_{k,1} \ge t_k^0$ depending on the size of the initial datum in $\mathcal{V}$, and a radius $\rho_{k,1}>0$ depending only on $f$ such that the estimate 
\be \la{smooglo3}
\|\l^{k+ \fr{1}{2}} q(t)\|_{L^2} + \int_{t}^{t+1} \|\l^{k+1} q(s)\|_{L^2}^2 ds \le \rho_{k,1}
\ee holds for all $t \ge T_{k,1}$. From \eqref{smooglo3}, we infer the existence of a time $T_{k,2} \ge T_{k,1}$ such that the $L^2$ norm of $\l^{k+1} q(T_{k,2})$ is bounded by some constant depending only on  $\|f\|_{H^{k}}$ and $\|\omega_0\|_{\mathcal{V}}$.

The $L^2$ norm of $\l^{k+1} q$ obeys 
\be 
\fr{d}{dt} \|\l^{k+1} q\|_{L^2}^2 + \|\l^{k+\fr{3}{2}} q\|_{L^2}^2 
\le C\|A^{\fr{k}{2} + 1}u\|_{L^2}^2 \|\l^{k+1}q\|_{L^2}^2
\ee as shown in \eqref{t26}. By using the assumption \eqref{smooglo1} and applying the uniform Gronwall Lemma \ref{gron}, we obtain a time $T_{k,3} \ge T_{k,2}$ depending only on $f$ and $\|\omega_0\|_{\mathcal{V}}$, and a radius $\rho_{k,2}>0$ depending only on $f$ such that the estimate 
\be \la{smooglo4}
\|\l^{k+1} q(t)\|_{L^2} + \int_{t}^{t+1} \|\l^{k+\fr{3}{2}} q(s)\|_{L^2}^2 ds \le \rho_{k,2}
\ee holds for all $t \ge T_{k,3}$. As for the $L^2$ norm of $A^{\fr{k}{2}+1}u$, we have
\be 
\fr{d}{dt} \|A^{\fr{k}{2}+1}u\|_{L^2}^2 + \|A^{\fr{k}{2} + \fr{3}{2}} u\|_{L^2}^2 
\le C\|\l^{k+1} q\|_{L^2}^4 + C\|A^{\fr{k}{2} + 1} u\|_{L^2}^4 + C\|A^{\fr{k}{2}} f\|_{L^2}^2
\ee as shown in \eqref{t23}. A use of \eqref{smooglo1} and \eqref{smooglo4} shows that the assumptions of the Gronwall Lemma \ref{gron} are satisfied and consequently, we obtain a time $T_{k,4} \ge T_{k,3}$ depending only on $\|\omega_0\|_{\mathcal{V}}$, and a radius $\rho_{k,3}>0$ depending only on $f$ such that 
\be \la{smooglo5}
\|A^{\fr{k}{2} + 1}u(t)\|_{L^2} + \int_{t}^{t+1} \|A^{\fr{k+3}{2}}u(s)\|_{L^2}^2 ds \le \rho_{k,3}
\ee holds for all $t \ge T_{k,4}$. Going back to \eqref{smooglo4}, we infer the existence of a time $T_{k,5} \ge T_{k,4}$ such that the $L^2$ norm of $\l^{k + \fr{3}{2}}q(T_{k,5})$ is bounded uniformly in $f$ and $\|\omega_0\|_{\mathcal{V}}$. We have thus completed the proof of Proposition~\ref{smooglo}.

We end this section by the proof of Theorem \ref{attp}:

\textbf{Proof of Theorem \ref{attp}.} The existence of the global attractor $\tilde{X}$ is based on the compactness of the absorbing ball $\mathcal{B}_{\tilde \rho}$ in the strong norm of $\mathcal{V}$ (Proposition \ref{cabsball}), the continuity of the solution map (Proposition \ref{conts}) and the injectivity of the solution map (Proposition \ref{backun}). The finite fractal dimensionality in $\mathcal{V}$ is a consequence the decay of volume elements (Proposition \ref{vol}) and the Lipschitz continuity property \eqref{Conthv}. The smoothness of the attractor follows from Proposition \ref{smooglo}.

\section{Global Gevrey Regularity in the Periodic Case: Proof of Theorem \ref{t3}} \la{s7}

The proof of Theorem \ref{t3} is based on the method of \cite{FT}, adapted for fractional dissipation. 

We need the following propositions:

\beg{prop} \la{pr5} Let $\tau \ge 0$ and $m > 2 $. Suppose $u \in \mathcal{D} (e^{\tau \l^{\fr{\alpha}{2}}}  \l^{\fr{m}{2} + 1})$ and $q \in \mathcal{D} (e^{\tau \l^{\fr{\alpha}{2}}}  \l^{\fr{m}{2} + \fr{\alpha}{2}})$. There exists a positive constant $C$ depending only on $m$ and $\alpha$ such that the following estimate 
\beg{align}
&|(e^{\tau \l^{\fr{\alpha}{2}}} \l^{\fr{m}{2}} (u \cdot \na q), e^{\tau \l^{\fr{\alpha}{2}}}\l^{\fr{m}{2}}q)_{L^2}| \nonumber
\\&\quad\quad\le C\|e^{\tau \l^{\fr{\alpha}{2}}}  \l^{\fr{m}{2} + 1} u\|_{L^2} \|e^{\tau \l^{\fr{\alpha}{2}}} \l^{\fr{m}{2}} q\|_{L^2} \left( \|e^{\tau \l^{\fr{\alpha}{2}}} \l^{\fr{m}{2}} q\|_{L^2}   + \tau \|e^{\tau \l^{\fr{\alpha}{2}}} \l^{\fr{m}{2} + \fr{\alpha}{2}} q\|_{L^2}\right) 
\end{align} holds.
\end{prop}

\noindent \textbf{Proof.} Let
\be \la{pr51}
u = \sum\limits_{j \in \ZZ^2 \setminus \left\{0\right\}} u_j e^{ij \cdot x},
\ee and 
\be \la{pr52}
q = \sum\limits_{j \in \ZZ^2 \setminus \left\{0\right\}} q_j e^{ij \cdot x}
\ee be the Fourier series expansions of $u$ and $q$ respectively. Denoting $e^{\tau \l^{\fr{\alpha}{2}}}u$ an $e^{\tau \l^{\fr{\alpha}{2}}}q$ by $u^*$ and $q^*$, we have 
\be \la{pr53}
u^* = \sum\limits_{j \in \ZZ^2 \setminus \left\{0\right\}} u_j^* e^{ij \cdot x}, u_j^* = e^{\tau|j|^{\fr{\alpha}{2}}} u_j,
\ee and 
\be \la{pr54}
q^* = \sum\limits_{j \in \ZZ^2 \setminus \left\{0\right\}} q_j^* e^{ij \cdot x}, q_j^* = e^{\tau|j|^{\fr{\alpha}{2}}} q_j. 
\ee
In view of the divergence-free condition $\na \cdot u = 0$, the $L^2$ cancellation
\be 
(u \cdot \na e^{\tau \l^{\fr{\alpha}{2}}} \l^{\fr{m}{2}}q, e^{\tau \l^{\fr{\alpha}{2}}} \l^{\fr{m}{2}}q)_{L^2} = 0 
\ee holds, hence 
\beg{align} \la{pr55}
&|(e^{\tau \l^{\fr{\alpha}{2}}} \l^{\fr{m}{2}} (u \cdot \na q), e^{\tau \l^{\fr{\alpha}{2}}}\l^{\fr{m}{2}}q)_{L^2}| 
= |(e^{\tau \l^{\fr{\alpha}{2}}} \l^{\fr{m}{2}} (u \cdot \na q) - u \cdot \na e^{\tau \l^{\fr{\alpha}{2}}} \l^{\fr{m}{2}}q, e^{\tau \l^{\fr{\alpha}{2}}}\l^{\fr{m}{2}}q)_{L^2}| \nonumber
\\&\quad\quad= (2\pi)^2 \left|\sum\limits_{j+k+ l = 0} (|l|^{\fr{m}{2}} e^{\tau |l|^{\fr{\alpha}{2}}} - |k|^{\fr{m}{2}} e^{\tau |k|^{\fr{\alpha}{2}}})  |l|^{\fr{m}{2}} e^{\tau |l|^{\fr{\alpha}{2}}}(u_j \cdot k) q_k q_{l} \right|.
\end{align} 
Applying the mean value theorem to the function $f(x) = x^{\fr{m}{2}} e^{\tau x^{\fr{\alpha}{2}}}$, whose derivative is given by $f'(x) = \fr{m}{2} x^{\fr{m}{2}-1} e^{\tau x^{\fr{\alpha}{2}}} + \fr{1}{2} \tau \alpha x^{\fr{m}{2} + \fr{\alpha}{2} - 1} e^{\tau x^{\fr{\alpha}{2}}}$, we bound the difference
\be 
\left||l|^{\fr{m}{2}} e^{\tau |l|^{\fr{\alpha}{2}}} - |k|^{\fr{m}{2}} e^{\tau |k|^{\fr{\alpha}{2}}}\right|
\le \left[\fr{m}{2} M^{\fr{m}{2} - 1}e^{\tau M^{\fr{\alpha}{2}} } + \fr{1}{2} \tau \alpha M^{\fr{m}{2} + \fr{\alpha}{2} - 1} e^{\tau M^{\fr{\alpha}{2}}}\right] ||k| - |l||
\ee for any $m \ge 2$, where $M := \max \left\{|k|, |l| \right\}.$ Consequently, we bound the sum \eqref{pr55} by 
\beg{align}
&|(e^{\tau \l^{\fr{\alpha}{2}}} \l^{\fr{m}{2}} (u \cdot \na q), e^{\tau \l^{\fr{\alpha}{2}}}\l^{\fr{m}{2}}q)_{L^2}| \nonumber
\\&\quad\le (2\pi)^2 \sum\limits_{j+k+ l = 0} \left[\fr{m}{2} M^{\fr{m}{2} - 1}e^{\tau M^{\fr{\alpha}{2}} } + \fr{1}{2} \tau \alpha M^{\fr{m}{2} + \fr{\alpha}{2} - 1} e^{\tau M^{\fr{\alpha}{2}}}\right] ||k| - |l|| |l|^{\fr{m}{2}} e^{\tau |l|^{\fr{\alpha}{2}}}|u_j \cdot k| |q_k| |q_{l}| \nonumber
\\&\quad\le (2\pi)^2 \sum\limits_{j+k+ l = 0} \left[\fr{m}{2} M^{\fr{m}{2} - 1}e^{\tau M^{\fr{\alpha}{2}} } + \fr{1}{2} \tau \alpha M^{\fr{m}{2} + \fr{\alpha}{2} - 1} e^{\tau M^{\fr{\alpha}{2}}}\right] |j| |l|^{\fr{m}{2}} e^{\tau |l|^{\fr{\alpha}{2}}}|u_j| |k| |q_k| |q_{l}|,
\end{align} where the last inequality uses the relation $j+k+l = 0$, that implies the inequality $||k| - |l|| \le |j|$. We split this latter series into the sum $S_1 + S_2 + S_3 + S_4$, where
\be 
S_1 = 2\pi^2m  \sum\limits_{j+k+ l = 0, |l| \le |k|}  M^{\fr{m}{2} - 1}e^{\tau M^{\fr{\alpha}{2}} } |j| |l|^{\fr{m}{2}} e^{\tau |l|^{\fr{\alpha}{2}}}|u_j| |k| |q_k| |q_{l}|,
\ee
\be 
S_2 =  2\pi^2m \sum\limits_{j+k+ l = 0, |k| < |l|}  M^{\fr{m}{2} - 1}e^{\tau M^{\fr{\alpha}{2}} } |j| |l|^{\fr{m}{2}} e^{\tau |l|^{\fr{\alpha}{2}}}|u_j| |k| |q_k| |q_{l}|,
\ee 
\be
S_3 = 2\pi^2  \alpha \tau  \sum\limits_{j+k+ l = 0, |l| \le |k|}  M^{\fr{m}{2} + \fr{\alpha}{2} - 1} e^{\tau M^{\fr{\alpha}{2}}} |j| |l|^{\fr{m}{2}} e^{\tau |l|^{\fr{\alpha}{2}}}|u_j| |k| |q_k| |q_{l}|,
\ee and 
\be 
S_4 = 2\pi^2  \alpha \tau  \sum\limits_{j+k+ l = 0, |k| \le |l|}  M^{\fr{m}{2} + \fr{\alpha}{2} - 1} e^{\tau M^{\fr{\alpha}{2}}} |j| |l|^{\fr{m}{2}} e^{\tau |l|^{\fr{\alpha}{2}}}|u_j| |k| |q_k| |q_{l}|.
\ee In view of H\"older, Young, and Plancherel inequalities, we estimate the sum $S_1$ as follows,
\beg{align} \la{pr56}
S_1
&= 2\pi^2m  \sum\limits_{j+k+ l = 0, |l| \le |k|}  |k|^{\fr{m}{2} - 1}e^{\tau |k|^{\fr{\alpha}{2}} } |j| |l|^{\fr{m}{2}} e^{\tau |l|^{\fr{\alpha}{2}}}|u_j| |k| |q_k| |q_{l}| \nonumber
\\&\le 2\pi^2m  \sum\limits_{j+k+ l = 0}   |j|  |u_j| |k|^{\fr{m}{2}} |q_k^*||l|^{\fr{m}{2}} |q_{l}^*| \nonumber
\\&\le C\||j||u_j|\|_{\ell^1(\ZZ^2 \setminus \left\{0\right\})} \||k|^{\fr{m}{2}}  |q_k^*| \|_{\ell^2(\ZZ^2 \setminus \left\{0\right\})} \||l|^{\fr{m}{2}} |q_{l}^*|  \|_{\ell^2(\ZZ^2 \setminus \left\{0\right\})} \nonumber
\\&\le C\|\l^{2 + \epsilon} u\|_{L^2} \|e^{\tau \l^{\fr{\alpha}{2}}} \l^{\fr{m}{2}} q\|_{L^2}^2
\end{align} for any $\epsilon > 0$. In order to bound the sum $S_2$, we use the relations $|l|^{\fr{\alpha}{2}} \le |k|^{\fr{\alpha}{2}}  + |j|^{\fr{\alpha}{2}} $ and $|l|^{\fr{m}{2}} \le 2^{\fr{m}{2}} \left(|k|^{\fr{m}{2}} + |j|^{\fr{m}{2}}\right)$ and obtain 
\beg{align}
S_2
&= 2\pi^2m \sum\limits_{j+k+ l = 0, |k| < |l|}  |l|^{\fr{m}{2} - 1}e^{\tau |l|^{\fr{\alpha}{2}} } |j| |l|^{\fr{m}{2}} e^{\tau |l|^{\fr{\alpha}{2}}}|u_j| |k| |q_k| |q_{l}| \nonumber
\\&\le 2^{\fr{m}{2} + 1}\pi^2 m \sum\limits_{j+k+ l = 0}  |l|^{\fr{m}{2}}  |j| \left(|k|^{\fr{m}{2}} + |j|^{\fr{m}{2}} \right) |u_j^*| |q_k^*| |q_{l}^*| \nonumber
\\&\le C\||j||u_j^*|\|_{\ell^1(\ZZ^2 \setminus \left\{0\right\})} \||k|^{\fr{m}{2}}  |q_k^*| \|_{\ell^2(\ZZ^2 \setminus \left\{0\right\})} \||l|^{\fr{m}{2}} |q_{l}^*|  \|_{\ell^2(\ZZ^2 \setminus \left\{0\right\})} \nonumber
\\&\quad\quad+ C\||j|^{\fr{m}{2}+1}|u_j^*|\|_{\ell^2(\ZZ^2 \setminus \left\{0\right\})} \|q_k^* \|_{\ell^1(\ZZ^2 \setminus \left\{0\right\})} \||l|^{\fr{m}{2}} |q_{l}^*|  \|_{\ell^2(\ZZ^2 \setminus \left\{0\right\})} \nonumber
\\&\le C\|e^{\tau \l^{\fr{\alpha}{2}}}  \l^{2+\epsilon} u\|_{L^2} \|e^{\tau \l^{\fr{\alpha}{2}}}  \l^{\fr{m}{2}} q\|_{L^2}^2  \nonumber
\\&\quad\quad+ C\|e^{\tau \l^{\fr{\alpha}{2}}}  \l^{\fr{m}{2} + 1} u\|_{L^2} \|e^{\tau \l^{\fr{\alpha}{2}}}  \l^{\fr{m}{2}} q\|_{L^2}\|e^{\tau \l^{\fr{\alpha}{2}}}  \l^{1 + \epsilon} q\|_{L^2}
\end{align} for any $\epsilon > 0$. In the same manner, we estimate 
\beg{align}
S_3 &= 2\pi^2  \alpha \tau  \sum\limits_{j+k+ l = 0, |l| \le |k|}  |k|^{\fr{m}{2} + \fr{\alpha}{2} - 1} e^{\tau |k|^{\fr{\alpha}{2}}} |j| |l|^{\fr{m}{2}} e^{\tau |l|^{\fr{\alpha}{2}}}|u_j| |k| |q_k| |q_{l}| \nonumber
\\&\le 2\pi^2 \alpha \tau \sum\limits_{j+k+ l = 0}   |j| |u_j| |k|^{\fr{m}{2} + \fr{\alpha}{2}}  |q_k^*| |l|^{\fr{m}{2}}  |q_{l}^*| \nonumber
\\&\le C\tau \|\l^{2 + \epsilon} u\|_{L^2} \|e^{\tau \l^{\fr{\alpha}{2}}} \l^{\fr{m}{2}} q\|_{L^2}\|e^{\tau \l^{\fr{\alpha}{2}}} \l^{\fr{m}{2} + \fr{\alpha}{2}} q\|_{L^2}
\end{align} for any $\epsilon > 0$, and 
\beg{align} \la{pr57}
S_4 &=  2\pi^2  \alpha \tau  \sum\limits_{j+k+ l = 0, |k| \le |l|}  |l|^{\fr{m}{2} + \fr{\alpha}{2} - 1} e^{\tau |l|^{\fr{\alpha}{2}}} |j| |l|^{\fr{m}{2}} e^{\tau |l|^{\fr{\alpha}{2}}}|u_j| |k| |q_k| |q_{l}| \nonumber
\\&\le C\tau \sum\limits_{j+k+ l = 0} |j| \left(|k|^{\fr{m}{2}} + |j|^{\fr{m}{2}} \right) |l|^{\fr{m}{2} + \fr{\alpha}{2}} |u_j^*| |q_k^*| |q_{l}^*|  \nonumber
\\&\le C\tau \|e^{\tau \l^{\fr{\alpha}{2}}}  \l^{2+\epsilon} u\|_{L^2} \|e^{\tau \l^{\fr{\alpha}{2}}}  \l^{\fr{m}{2}} q\|_{L^2} \|e^{\tau \l^{\fr{\alpha}{2}}}  \l^{\fr{m}{2} + \fr{\alpha}{2}} q\|_{L^2}  \nonumber
\\&\quad\quad+ C\tau \|e^{\tau \l^{\fr{\alpha}{2}}}  \l^{\fr{m}{2} + 1} u\|_{L^2} \|e^{\tau \l^{\fr{\alpha}{2}}}  \l^{1 + \epsilon} q\|_{L^2} \|e^{\tau \l^{\fr{\alpha}{2}}}  \l^{\fr{m}{2} + \fr{\alpha}{2}}  q\|_{L^2}
\end{align} for any $\epsilon > 0$. 
Adding \eqref{pr56}--\eqref{pr57} and choosing $m > 2$, we infer that
\beg{align}
&|(e^{\tau \l^{\fr{\alpha}{2}}} \l^{\fr{m}{2}} (u \cdot \na q), e^{\tau \l^{\fr{\alpha}{2}}}\l^{\fr{m}{2}}q)_{L^2}| 
\le C\left(\|\l^{\fr{m}{2}+ 1} u\|_{L^2} + \|e^{\tau \l^{\fr{\alpha}{2}}}  \l^{\fr{m}{2} + 1} u\|_{L^2}\right)\|e^{\tau \l^{\fr{\alpha}{2}}} \l^{\fr{m}{2}} q\|_{L^2}^2 \nonumber
\\&+ C\tau \left(\|\l^{\fr{m}{2} + 1} u\|_{L^2} + \|e^{\tau \l^{\fr{\alpha}{2}}}  \l^{\fr{m}{2} + 1} u\|_{L^2} \right)\|e^{\tau \l^{\fr{\alpha}{2}}} \l^{\fr{m}{2}} q\|_{L^2}\|e^{\tau \l^{\fr{\alpha}{2}}} \l^{\fr{m}{2} + \fr{\alpha}{2}} q\|_{L^2},
\end{align} finishing the proof of Proposition \ref{pr5}.

\beg{prop} \la{pr6} Let $\tau \ge 0$ and $m > 2 $. Suppose $u \in \mathcal{D} (e^{\tau \l^{\fr{\alpha}{2}}}  \l^{\fr{m}{2} + 2})$ and $q \in \mathcal{D} (e^{\tau \l^{\fr{\alpha}{2}}}  \l^{\fr{m}{2}})$. There exists a positive constant $C$ depending only on $m$  such that the following estimate 
\beg{align} \la{pr61}
|(e^{\tau \l^{\fr{\alpha}{2}}} \l^{\fr{m}{2}+1} (qRq), e^{\tau \l^{\fr{\alpha}{2}}}\l^{\fr{m}{2}+1}u)_{L^2}| 
\le C\|e^{\tau \l^{\fr{\alpha}{2}}} \l^{\fr{m}{2}+2}u\|_{L^2} \|e^{\tau \l^{\fr{\alpha}{2}}} \l^{\fr{m}{2}}q\|_{L^2}^2
\end{align} holds.
\end{prop}

\noindent \textbf{Proof.} We set $u, q, u^*$ and $q^*$ as in \eqref{pr51}--\eqref{pr54}. The Fourier series expansion of $Rq = \na \l^{-1}q$ is given by 
\be 
Rq = \sum\limits_{j \in \ZZ^2 \setminus \left\{0\right\}} i \fr{j}{|j|} q_j e^{ij \cdot x}.
\ee Thus, we have 
\beg{align}
&|(e^{\tau \l^{\fr{\alpha}{2}}} \l^{\fr{m}{2}+1} (qRq), e^{\tau \l^{\fr{\alpha}{2}}}\l^{\fr{m}{2}+1}u)_{L^2}| 
= \left|(2\pi)^2i \sum\limits_{j + k + l = 0} |l|^{m+2} e^{2\tau |l|^{\fr{\alpha}{2}}}q_j (u_l \cdot k)|k|^{-1} q_k \right| \nonumber
\\&\quad\le C\sum\limits_{j + k + l = 0}  |l|^{\fr{m}{2} + 2} (|k|^{\fr{m}{2}} + |j|^{\fr{m}{2}}) |u_l^*| |q_k^*| |q_j^*| \nonumber
\\&\quad\le C\|q_j^* \|_{\ell^1(\ZZ^2 \setminus \left\{0\right\})} \||k|^{\fr{m}{2}} |q_k^*| \|_{\ell^2(\ZZ^2 \setminus \left\{0\right\})} \||l|^{\fr{m}{2}+2} |u_l^*| \|_{\ell^2(\ZZ^2 \setminus \left\{0\right\})} \nonumber
\\&\quad\le C\|e^{e^{\tau \l^{\fr{\alpha}{2}}}} \l^{1+\epsilon} q\|_{L^2}  \|e^{\tau \l^{\fr{\alpha}{2}}} \l^{\fr{m}{2}}q\|_{L^2} \|e^{\tau \l^{\fr{\alpha}{2}}} \l^{\fr{m}{2}+2}u\|_{L^2}
\end{align} for any $\epsilon > 0$. Therefore, we obtain \eqref{pr61} provided that $m > 2$.

\beg{prop} \la{pr7} Let $\tau \ge 0$ and $m > 2 $. Suppose $u \in \mathcal{D} (e^{\tau \l^{\fr{\alpha}{2}}}  \l^{\fr{m}{2} + 2})$. There exists a positive constant $C$ depending only on $m$  such that the following estimate 
\beg{align} \la{pr71}
|(e^{\tau \l^{\fr{\alpha}{2}}} \l^{\fr{m}{2}+1} u \cdot \na u, e^{\tau \l^{\fr{\alpha}{2}}}\l^{\fr{m}{2}+1}u)_{L^2}| 
\le C\|e^{\tau \l^{\fr{\alpha}{2}}} \l^{\fr{m}{2}+2}u\|_{L^2} \|e^{\tau \l^{\fr{\alpha}{2}}} \l^{\fr{m}{2}+1}u\|_{L^2}^2
\end{align} holds.
\end{prop}

\noindent \textbf{Proof.} Setting $u$ and $u^*$ as in \eqref{pr51} and \eqref{pr53} respectively, we estimate 
\beg{align}
&|(e^{\tau \l^{\fr{\alpha}{2}}} \l^{\fr{m}{2}+1} u \cdot \na u, e^{\tau \l^{\fr{\alpha}{2}}}\l^{\fr{m}{2}+1}u)_{L^2}| 
= \left|\sum\limits_{j+k+l = 0} e^{2 \tau |l|^{\fr{\alpha}{2}}} |l|^{m+2} [(u_j \cdot k) u_k] \cdot u_l \right| \nonumber
\\&\quad\le C\sum\limits_{j+k+l = 0}  |l|^{\fr{m}{2} + 2} \left(|k|^{\fr{m}{2}} + |j|^{\fr{m}{2}} \right) |k| |u_j^*||u_k^*| |u_l^*| \nonumber
\\&\quad\le C\|\l^{1+\epsilon} u^*\|_{L^2} \|\l^{\fr{m}{2} + 1} u^* \|_{L^2} \|\l^{\fr{m}{2} + 2} u^*\|_{L^2}
+ C \|\l^{2+\epsilon} u^*\|_{L^2} \|\l^{\fr{m}{2}} u^*\|_{L^2} \|\l^{\fr{m}{2} + 2} u^* \|_{L^2}
\end{align} for any $\epsilon > 0$, yielding \eqref{pr71}. 

We end this section by proving Theorem \ref{t3}:

\noindent \textbf{Proof of Theorem \ref{t3}.} The proof is divided into two main steps:

\textbf{Step 1. Local Gevrey Regularity.} We take the scalar product in $\mathcal{D}(e^{\tau(t) \l^{\fr{\alpha}{2}}})$ of the equation \eqref{eq1} obeyed by the charge density $q$ with $\l^{m} q$. We obtain the energy equality
\beg{align} \la{t32}
&\fr{1}{2} \fr{d}{dt} \|e^{\tau(t) \l^{\fr{\alpha}{2}}}  \l^{\fr{m}{2}}q \|_{L^2}^2 
- \tau'(t)\|e^{\tau(t) \l^{\fr{\alpha}{2}}} \l^{\fr{m}{2} + \fr{\alpha}{4}} q\|_{L^2}^2
+ \|e^{\tau(t) \l^{\fr{\alpha}{2}}} \l^{\fr{m}{2} + \fr{\alpha}{2}} q\|_{L^2}^2 \nonumber
\\&\quad\quad= - (e^{\tau(t) \l^{\fr{\alpha}{2}}} \l^{\fr{m}{2}} (u \cdot \na q),e^{\tau(t) \l^{\fr{\alpha}{2}}} \l^{\fr{m}{2}}q )_{L^2}.
\end{align} 
We estimate the nonlinear term by making use of Proposition \ref{pr5}, Young's inequality, and the boundedness of $\tau$ by $1$, yielding  
\beg{align} 
&|(e^{\tau(t) \l^{\fr{\alpha}{2}}} \l^{\fr{m}{2}} (u \cdot \na q),e^{\tau(t) \l^{\fr{\alpha}{2}}} \l^{\fr{m}{2}}q )_{L^2}| \nonumber
\\&\quad\quad\le \fr{1}{4} \|e^{\tau \l^{\fr{\alpha}{2}}} \l^{\fr{m}{2} + \fr{\alpha}{2}} q\|_{L^2}^2
+ C\left(\|e^{\tau \l^{\fr{\alpha}{2}}}  \l^{\fr{m}{2} + 1} u\|_{L^2}^2 + \|e^{\tau \l^{\fr{\alpha}{2}}}  \l^{\fr{m}{2} + 1} u\|_{L^2} \right) \|e^{\tau \l^{\fr{\alpha}{2}}} \l^{\fr{m}{2}} q\|_{L^2}^2.
\end{align} 
Since $\tau'(t) \le \fr{1}{4}$, we obtain the differential inequality
\beg{align} \la{t35}
&\fr{d}{dt} \|e^{\tau(t) \l^{\fr{\alpha}{2}}}  \l^{\fr{m}{2}}q \|_{L^2}^2 
+ \|e^{\tau(t) \l^{\fr{\alpha}{2}}} \l^{\fr{m}{2} + \fr{\alpha}{2}} q\|_{L^2}^2 \nonumber
\\&\quad\quad\le C\left(\|e^{\tau \l^{\fr{\alpha}{2}}}  \l^{\fr{m}{2} + 1} u\|_{L^2}^2 + \|e^{\tau \l^{\fr{\alpha}{2}}}  \l^{\fr{m}{2} + 1} u\|_{L^2} \right) \|e^{\tau \l^{\fr{\alpha}{2}}} \l^{\fr{m}{2}} q\|_{L^2}^2.
\end{align}
Now we take the scalar product in $\mathcal{D}(e^{\tau(t) \l^{\fr{\alpha}{2}}})$ of the equation \eqref{eq11} obeyed by the velocity $u$ with $\l^{m + 2} u$. Due to the divergence-free condition obeyed by $u$, we have the cancellation
\be 
(e^{\tau(t) \l^{\fr{\alpha}{2}}} \na p, e^{\tau(t) \l^{\fr{\alpha}{2}}} \l^{m+2} u)_{L^2} = 0.
\ee Hence the $L^2$ norm of $e^{\tau(t) \l^{\fr{\alpha}{2}}} \l^{\fr{m}{2} + 1} u$ evolves according to 
\beg{align}
&\fr{1}{2} \fr{d}{dt} \|e^{\tau(t) \l^{\fr{\alpha}{2}}} \l^{\fr{m}{2} + 1} u\|_{L^2}^2
- \tau'(t) \|e^{\tau(t) \l^{\fr{\alpha}{2}}} \l^{\fr{m}{2} + \fr{\alpha}{4} + 1} u\|_{L^2}^2
+ \|e^{\tau(t) \l^{\fr{\alpha}{2}}} \l^{\fr{m}{2} + 2} u\|_{L^2}^2 \nonumber
\\&= -(e^{\tau(t) \l^{\fr{\alpha}{2}}} \l^{\fr{m}{2} + 1}(u \cdot \na u), e^{\tau(t) \l^{\fr{\alpha}{2}}} \l^{\fr{m}{2} + 1} u)_{L^2}
- (e^{\tau(t) \l^{\fr{\alpha}{2}}} \l^{\fr{m}{2} + 1} (qRq), e^{\tau(t) \l^{\fr{\alpha}{2}}} \l^{\fr{m}{2} + 1} u)_{L^2}.
\end{align} In view of Propositions \ref{pr6} and \ref{pr7} followed by applications of Young's inequality for products, we obtain the differential inequality
\beg{align} \la{t36}
\fr{d}{dt} \|e^{\tau(t) \l^{\fr{\alpha}{2}}} \l^{\fr{m}{2} + 1} u\|_{L^2}^2
+ \|e^{\tau(t) \l^{\fr{\alpha}{2}}} \l^{\fr{m}{2} + 2} u\|_{L^2}^2
\le C\|e^{\tau(t) \l^{\fr{\alpha}{2}}} \l^{\fr{m}{2} + 1} u\|_{L^2}^4 + C\|e^{\tau(t) \l^{\fr{\alpha}{2}}} \l^{\fr{m}{2}} q\|_{L^2}^4 .
\end{align}
We add \eqref{t35} and \eqref{t36}. Setting
\be 
y(t) = \|e^{\tau(t) \l^{\fr{\alpha}{2}}} \l^{\fr{m}{2}} q(t)\|_{L^2}^2 + \|e^{\tau(t) \l^{\fr{\alpha}{2}}} \l^{\fr{m}{2} + 1} u(t)\|_{L^2}^2,
\ee we have 
\be 
y'(t) \le C y(t)^2
\ee for all $t \ge 0$. where  Dividing both sides by $y(t)^2$ and integrating in time from $0$ to $t$, we obtain
\be 
\fr{1}{y(t)} \ge \fr{1}{y(0)} - C t \ge \fr{1}{2y(0)}
\ee provided that 
\be 
t \le \fr{1}{2 y(0)} := T_0.
\ee Therefore, 
\be \la{t58}
\|e^{\tau(t) \l^{\fr{\alpha}{2}}} \l^{\fr{m}{2}} q(t)\|_{L^2}^2 + \|e^{\tau(t) \l^{\fr{\alpha}{2}}} \l^{\fr{m}{2} + 1} u(t)\|_{L^2}^2
\le 2\|\l^{\fr{m}{2}} q_0\|_{L^2}^2 + 2\|\l^{\fr{m}{2} + 1} u_0\|_{L^2}^2
\ee for all $t \in [0, T_0]$. 

\textbf{Step 2. Extension of the local analyticity property.} For a fixed real number $m > 2$, we prove that the charge density $q$ is bounded in $L^{\infty}(0,\infty, H^{\fr{m}{2}}(\TT^2))$ and the velocity $u$ in bounded in $L^{\infty}(0,\infty, H^{\fr{m}{2}+ 1} (\TT^2))$, from which we can conclude that the Gevrey regularity \eqref{t58} propagates from the short time interval $(0,T_0)$ into $(0, \infty)$. For that objective, we show that
\be \la{extension}
\|\l^{\fr{m}{2}} q(t)\|_{L^2}^2 + \|\l^{\fr{m}{2} + 1}u(t)\|_{L^2}^2 \le C(\|\l^{\fr{m}{2}}q_0\|_{L^2}, \|\l^{\fr{m}{2}}u_0\|_{L^2}) e^{-ct}
\ee for all $t \ge 0$.
Indeed, the norm $\|\l^{\fr{m}{2}}q\|_{L^2}^2 + \|\l^{\fr{m}{2} + 1}u\|_{L^2}^2$ obeys
\beg{align} \la{extension1}
&\fr{1}{2} \fr{d}{dt} \left[\|\l^{\fr{m}{2}}q\|_{L^2}^2 + \|\l^{\fr{m}{2} + 1}u\|_{L^2}^2 \right]
+ \|\l^{\fr{m+ \alpha}{2}}q\|_{L^2}^2 + \|\l^{\fr{m}{2} + 2}u\|_{L^2}^2 \nonumber
\\&\quad\quad= - (\l^{\fr{m}{2}} (u \cdot \na q), \l^{\fr{m}{2}}q)_{L^2}
- (\l^{\fr{m}{2}}(qRq), \l^{\fr{m}{2} + 2}u)_{L^2}
- (\l^{\fr{m}{2} + 1} (u \cdot \na u), \l^{\fr{m}{2} + 1} u)_{L^2}
\end{align} 
We estimate 
\beg{align} 
&|(\l^{\fr{m}{2}} (u \cdot \na q), \l^{\fr{m}{2}}q)_{L^2} |
= |(\l^{\fr{m}{2}} (u \cdot \na q) - u \cdot \na \l^{\fr{m}{2}}q, \l^{\fr{m}{2}}q)_{L^2} |  \nonumber 
\\&\quad\quad\le C\left[\|\na u\|_{L^{\infty}}\|\l^{\fr{m}{2}}q\|_{L^2} + \|\l^{\fr{m}{2}}u\|_{L^{\fr{4}{\alpha}}} \|\na q\|_{L^{\fr{4}{2-\alpha}}}  \right]\|\l^{\fr{m}{2}}q\|_{L^2} \nonumber
\\&\quad\quad\le C\left[\|\na \Delta u\|_{L^{2}}\|\l^{\fr{m}{2}}q\|_{L^2} + \|\l^{\fr{m}{2} + 1}u\|_{L^2} \|\l^{1+\fr{\alpha}{2}} q\|_{L^2}  \right]\|\l^{\fr{m}{2}}q\|_{L^2} \nonumber
\\&\quad\quad\le \fr{1}{4} \|\l^{\fr{m}{2} + 2}u\|_{L^2}^2 +C\left[\|\l^{1+\fr{\alpha}{2}} q\|_{L^2}^2 + \|\na \Delta u\|_{L^2}^2 \right]  \|\l^{\fr{m}{2}}q\|_{L^2} ^2, 
\end{align}
\beg{align}
&|(\l^{\fr{m}{2}}(qRq), \l^{\fr{m}{2} + 2}u)_{L^2} | 
\le \|\l^{\fr{m}{2} + 2}u\|_{L^2} \|\l^{\fr{m}{2}}(qRq) \|_{L^2} \nonumber
\\&\quad\quad\le C\|\l^{\fr{m}{2} + 2}u\|_{L^2}\left[\|\l^{\fr{m}{2}}q\|_{L^2}\|Rq\|_{L^{\infty}} + \|R\l^{\fr{m}{2}}q\|_{L^2} \|q\|_{L^{\infty}} \right] \nonumber
\\&\quad\quad\le C\|\l^{\fr{m}{2} + 2}u\|_{L^2} \|\l^{\fr{m}{2}}q\|_{L^2}\|\l^{1+ \fr{\alpha}{2}}q\|_{L^2} \nonumber
\\&\quad\quad\le \fr{1}{8} \|\l^{\fr{m}{2} + 2}u\|_{L^2}^2 + C\|\l^{1+ \fr{\alpha}{2}}q\|_{L^2}^2 \|\l^{\fr{m}{2}}q\|_{L^2}^2,
\end{align} and 
\beg{align} \la{extension2}
&|(\l^{\fr{m}{2} + 1} (u \cdot \na u), \l^{\fr{m}{2} + 1} u)_{L^2} | 
 = |(\l^{\fr{m}{2} + 1} (u \cdot \na u) - u \cdot \na \l^{\fr{m}{2} + 1} u , \l^{\fr{m}{2} + 1} u)_{L^2} | \nonumber
 \\&\quad\quad\le C\|\na u\|_{L^{\infty}} \|\l^{\fr{m}{2} + 1}u\|_{L^2}^2
 \le \fr{1}{8} \|\l^{\fr{m}{2} + 2}u\|_{L^2}^2 + C\|\na \Delta u\|_{L^{2}}^2 \|\l^{\fr{m}{2} + 1}u\|_{L^2}^2
\end{align} by making use of the continuous Sobolev embeddings $H^{\fr{\alpha}{2}}(\TT^2) \subset L^{\fr{4}{2-\alpha}}(\TT^2)$, $H^1(\TT^2) \subset L^{\fr{4}{\alpha}}(\TT^2)$, and $H^{1+\epsilon}(\TT^2) \subset L^{\infty}(\TT^2)$ that hold for any $\epsilon > 0$, the boundedness of the Riesz transform on Sobolev spaces, periodic fractional product and commutator estimates \cite[Appendix A]{cvt}, and Young's inequality for products. Putting \eqref{extension1}--\eqref{extension2} together, we obtain the energy inequality
\be 
\fr{d}{dt} \left[\|\l^{\fr{m}{2}}q\|_{L^2}^2 + \|\l^{\fr{m}{2} + 1}u\|_{L^2}^2 \right]
\le C\left[\|\na \Delta u\|_{L^2}^2 + \|\l^{1 + \fr{\alpha}{2}}q\|_{L^2}^2 \right] \left[\|\l^{\fr{m}{2}}q\|_{L^2}^2 + \|\l^{\fr{m}{2} + 1}u\|_{L^2}^2 \right].
\ee Since 
\be 
\int_{0}^{\infty} \left[\|\l^{1+\fr{\alpha}{2}}q(s)\|_{L^2}^2 + \|\na \Delta u(s)\|_{L^2}^2\right] ds \le C(\|\l q_0\|_{L^2}, \|\l^2 u_0\|_{L^2})
\ee holds for all $t \ge 0$ (see Theorems \ref{Existence} and \ref{tt2}), we conclude that $(q,u)$ satisfies \eqref{extension}.  We have thus finished the proof of Step 2, completing the proof of Theorem \ref{t3}.

\appendix

\section{Spectral Lemma} \la{s9}

We present a lemma describing the asymptotic behavior of eigenvalues associated with a vector-valued operator:

\beg{lem} \la{voll} Let $\tilde{H}$ be a Hilbert space. Suppose $A_1$ and $A_2$ are operators defined on $\mathcal{D}(A_1) \subset \tilde{H}$ and $\mathcal{D}(A_2) \subset \tilde{H}$ respectively
\be 
A_1: \mathcal{D}(A_1) \subset \tilde{H} \mapsto \tilde{H},
\ee
\be 
A_2: \mathcal{D}(A_2) \subset \tilde{H} \mapsto \tilde{H},
\ee
such that $A_1$ and $A_2$ are strictly positive and injective, with compact inverses, $A_1^{-1}$ and $A_2^{-1}$, in $\tilde{H}$. Let $\tilde{A}$ be the operator defined on $\mathcal{D}(A_1) \times \mathcal{D}(A_2)$ by 
\be 
\tilde{A} (a_1, a_2) = (A_1a_1, A_2a_2).
\ee  Then $\tilde{A}, A_1$ and $A_2$ have unbounded increasing sequences of eigenvalues, $\left\{\mu_j \right\}_{j=1}^{\infty}$, $\left\{\lambda_j^1 \right\}_{j=1}^{\infty}$ and $\left\{\lambda_j^2 \right\}_{j=1}^{\infty}$ respectively, such that 
\be 
\left\{\mu_j \right\}_{j=1}^{\infty} = \left\{\lambda_j^1 \right\}_{j=1}^{\infty} \cup \left\{\lambda_j^2 \right\}_{j=1}^{\infty}.
\ee If $\lambda_j^1 \ge c_1 j^{\beta_1}$ and $\lambda_j^2 \ge c_1 j^{\beta_2}$ for all nonnegative integers $j$, then 
\be 
\mu_j \ge \fr{\min \left\{c_1, c_2\right\}}{2^{1 + \min \left\{\beta_1, \beta_2 \right\}}} j^{\min \left\{\beta_1, \beta_2 \right\}}
\ee for all integers $j \ge 0$. Consequently, the sum of the first $N$ eigenvalues of $\tilde{A}$ obeys
\be \la{volll}
\mu_1 + \dots + \mu_N \ge C_{\beta_1, \beta_2} \min\left\{c_1, c_2\right\} N^{1 + \min \left\{\beta_1, \beta_2 \right\}} 
\ee for some positive constants $C_{\beta_1, \beta_2}$ depending only $\beta_1$ and $\beta_2$.
\end{lem} 

\noindent \textbf{Proof.} The operators $A_1^{-1}$ and $A_2^{-1}$ are self-adjoint, injective, and compact, with ranges $\mathcal{D}(A_1)$ and $\mathcal{D}(A_2)$ respectively. By the spectral theory for Hilbert spaces, there are orthonormal bases of $\tilde{H}$, $\left\{\xi_j^1 \right\}_{j=1}^{\infty}$ and $\left\{\xi_j^2\right\}_{j=1}^{\infty}$, consisting of eigenfunctions of the operators $A_1$ and $A_2$ respectively, such that 
\be 
A_1 \xi_j^1 = \lambda_j^1 \xi_j^1,
\ee
\be 
A_2 \xi_j^2 = \lambda_j^2 \xi_j^2,
\ee with $0 < \lambda_1^1 \le \lambda_2^1 \le \dots \le \lambda_j^1 \le \lambda_{j+1}^1 \le \dots \rightarrow \infty$ and $0 < \lambda_1^2 \le \lambda_2^2 \le \dots \le \lambda_j^2 \le \lambda_{j+1}^2 \le \dots \rightarrow \infty$.
The operator $\tilde{A}^{-1}$ is also self-adjoint, injective, and compact in $\tilde{H} \times \tilde{H}$, so there is an orthonormal basis of $\tilde{H} \times \tilde{H}$ consisting of eigenvectors $\left\{\xi_j \right\}_{j=1}^{\infty}$ of $\tilde{A}$, such that 
\be 
\tilde{A}\xi_j = \mu_j \xi_j,
\ee with $0 < \mu_1 \le \mu_2 \le \dots \le \mu_j \le \mu_{j+1} \le \dots \rightarrow \infty$. The eigenvalues of $\tilde{A}$ are precisely the collection of eigenvalues of $A_1$ and $A_2$, counted with multiplicity. For $N \ge 1$, we have 
\be 
\left\{\mu_i : i =1, \dots, N \right\} = \left\{\lambda_i^1 : i = 1, \dots, j\right\} \cup \left\{\lambda_i^2: i= 1, \dots, k \right\}
\ee for some nonnegative integers $j$ and $k$ obeying $N = j+k$. If $\mu_N = \lambda_j^1$, then $\mu_N \ge c_1j^{\beta_1}$ and $\mu_N \ge \lambda_k^2 \ge c_2 k^{\beta_2}$. If $\mu_N = \lambda_k^2$, then $\mu_N \ge c_2 k^{\beta_2}$ and $\mu_N \ge \lambda_j^1 \ge c_1 j^{\beta_1}$. Thus, we infer that
\beg{align}
\mu_N &= \fr{1}{2} \mu_N + \fr{1}{2} \mu_N
\ge \fr{c_1}{2} j^{\beta_1} + \fr{c_2}{2} k^{\beta_2}
\ge \fr{1}{2} \min\left\{c_1, c_2 \right\} \left[j^{\min \left\{\beta_1, \beta_2\right\}} + k^{\min \left\{\beta_1, \beta_2\right\}} \right] \nonumber
\\&\quad\quad\ge \fr{\min \left\{c_1, c_2\right\}}{2^{1+ \min\left\{\beta_1, \beta_2 \right\}}} (j+k)^{\min \left\{\beta_1, \beta_2\right\}}
= \fr{\min \left\{c_1, c_2\right\}}{2^{1+ \min\left\{\beta_1, \beta_2 \right\}}} N^{\min \left\{\beta_1, \beta_2\right\}}.
\end{align} As a consequence of these latter lower bounds, we obtain \eqref{volll}. This ends the proof of Lemma \ref{voll}.

\section{Uniform Gronwall Lemma} \la{s8}

We present a Gronwall Lemma that will be used to study the time asymptotic behavior of solutions. 

\beg{lem} \la{gron} Let $y(t)$ be a nonnegative function of time $t$ that solves the differential inequality
\be  \la{gron1}
\fr{d}{dt} y + cy \le C_1 + C_2F_1 + C_3F_2y^n,
\ee where $c>0$ is a positive real number, $C_1, C_2$ and $C_3$ are nonnegative real numbers, $n$ is a nonnegative integer, and $F_1$ and $F_2$ are nonegative functions of time $t$. Suppose there exists a time $t_0$ and a positive number $R$ such that $y(t_0) < \infty$ and, for any $t \ge t_0$, it holds that 
\be \la{gron2}
\int_{t}^{t+1} F_1(s) ds \le R
\ee if $C_3 = 0$, and 
\be \la{gron33}
\int_{t}^{t+1} \left[F_1(s) + F_2(s)y^{n-1}(s) + y(s) \right] ds \le R
\ee if $C_3 \ne 0$ and $n \ge 1$.  Then  
\be \la{gron4}
y(t) \le \left(c^{-1} C_1 + 2C_2R + 2R\right) e^{2C_3R}
\ee for all times $t \ge t_0 + 1$. 
\end{lem}

\noindent \textbf{Proof.} We distinguish two cases: $C_3 \ne 0$, $n \ge 1$ and $C_3 = 0$. In the first case, we fix two times $s$ and $t$ such that $t_0 \le s \le t$. We multiply both sides of the inequality \eqref{gron1} by $e^{ct - C_3\int_{s}^{t} F_2 y^{n-1}(\tau) d\tau}$ and integrate in time from $s$ to $t$. We obtain the bound 
\be 
y(t) \le \left(y(s) + \fr{C_1}{c} + C_2\int_{s}^{t} F_1(\tau) d\tau \right) \exp \left\{C_3 \int_{s}^{t} F_2 (\tau) y^{n-1}(\tau) d\tau \right\}.
\ee In view of \eqref{gron33}, we have 
\be 
\int_{t_0 + k}^{t_0 + k + 1} y(\tau) d\tau \le R
\ee for any nonnegative integer $k \ge 0$. Thus, for each integer $k \ge 0$, there exists a time $\tilde t_k \in [t_0 + k, t_0 + k + 1]$ such that  
\be 
y(\tilde t_k) \le 2R.
\ee We note that the distance between two consecutive times $\tilde t_k$ and $\tilde t_{k+1}$ does not exceed two. By making use of \eqref{gron33}, we infer that 
\beg{align} \la{gron7}
y(t) 
&\le \left(y(\tilde t_k) +\fr{ C_1}{c} + C_2 \int_{\tilde t_k}^{\tilde t_{k+1}} F_1(\tau) d\tau  \right)  \exp \left\{C_3 \int_{\tilde t_k}^{\tilde t_{k+1}} F_2 (\tau) y^{n-1}(\tau) d\tau \right\} \nonumber
\\&\le \left(2R + \fr{C_1}{c} + 2C_2R\right) e^{2C_3R}
\end{align}  for any $t \in [\tilde t_k, \tilde t_{k+1}]$. Therefore, \eqref{gron7} holds on the time interval $[\tilde t_0, \infty)$, yielding the desired bound \eqref{gron4}. In the case where $C_3$ vanishes, the estimate \eqref{gron4} holds as a consequence of \cite[Lemma 1]{AI}. 

\vspace{.2in}

{\bf{Acknowledgment.}} The work of M.I. was partially supported by NSF grant DMS 2204614.\\

{\bf{Data Availability Statement.}} The research does not have any associated data.\\

{\bf{Conflict of Interest.}} The authors declare that they have no conflict of interest.\\


\begin{thebibliography}{99}


\bibitem{AI} E.~Abdo, M.~Ignatova, \emph{Long time dynamics of a model of electroconvection}, Trans. Amer. Math. Soc.~{\bf{374}}, 5849--5875 (2021).



\bibitem{BW} A. Bonito, P. Wei, \emph{Electroconvection of thin liquid crystals: Model reduction and numerical simulations}, Journal of Computational Physics~{\bf{405}} (2020). 

\bibitem{BM} H.~Brezis, P.~Mironescu, \emph{Gagliardo-Nirenberg inequalities and non-inequalities: The full story}, Annales de l'Institut Henri Poincar\'e C, Analyse non linéaire~{\bf{35}} (5), 355--1376 (2018).




\bibitem{CS} L.A. Caffarelli, P.R. Stinga, \emph{Fractional elliptic equations, Caccioppoli estimates and regularity}, Annales de l'Institut Henri Poincar\'e C, Analyse non linéaire~{\bf{33}} (3), 767--807 (2016).

\bibitem{caf} L.A. Caffarelli, A.~Vasseur, \emph{Drift diffusion equations with fractional diffusion and the quasi-geostrophic equation},  Ann. of Math.~{\bf{171}}(3), 1903--1930 (2010).



\bibitem{ceiv}  P.~Constantin, T.~Elgindi, M.~Ignatova, V.~Vicol, \emph{On some electroconvection models}, Journal of Nonlinear Science~{\bf{27}}, 197--211 (2017).

\bibitem{cfbook} P. Constantin, C. Foias, \emph{Navier-Stokes equations}, Chicago University Press, Chicago (1988).


\bibitem{CI} P. Constantin, M. Ignatova, \emph{Critical SQG in bounded domains}, M. Ann. PDE~{\bf{2}} (8) (2016).


\bibitem{CI3} P. Constantin, M. Ignatova, \emph{Estimates near the boundary for critical SQG}, Ann. PDE~{\bf{6}} (1) (2020).

\bibitem{CMT} P. Constantin, A. Majda, and E. Tabak, \emph{Formation of strong fronts in the 2-D quasi-geostrophic thermal active scalar}, Nonlinearity~{\bf{7}}, 1495--1533 (1994).

\bibitem{CN1} P. Constantin, H. Nguyen, \emph{Global weak solutions for generalized SQG in bounded domains}, Anal. PDE~{\bf{11}} (4), 1029--1047 (2018).


\bibitem{CN2} P. Constantin, H. Nguyen, \emph{Global weak solutions for SQG in bounded domains}, Communications on Pure and Applied Mathematics~{\bf{71}} (11),  2323--2333 (2017). 

\bibitem{CN3} P. Constantin, H. Nguyen, \emph{Local and global strong solutions for SQG in bounded domains}, Physica D: Nonlinear Phenomena~{\bf{376--377}}, 195--203 (2018). 

\bibitem{cvt} P. Constantin, A. Tarfulea, V. Vicol, \emph{Long time dynamics of forced critical SQG}, Communications in Mathematical Physics~{\bf{335}} (1), 93--141 (2015).

\bibitem{CV} P. Constantin, V. Vicol, \emph{Nonlinear maximum principles for dissipative linear nonlocal operators and applications},  Geometric and Functional Analysis~{\bf{22}} (5), 1289--1321 (2012).

\bibitem{CW1} P.~Constantin, J.~Wu, \emph{Behavior of solutions of 2D quasi-geostrophic equations}, SIAM J. Math. Anal.~{\bf{30}}, 937--948 (1999).

\bibitem{CW} P. Constantin, J. Wu, \emph{Regularity of H\"older continuous solutions of the supercritical quasi-geostrophic equation}, Ann. Inst. H. Poincaré Anal. Non Linéaire~{\bf{25}} (6), 1103--1110 (2008).

\bibitem{CZV} M. Coti Zelati, V. Vicol, \emph{On the global regularity for the supercritical SQG equation}, Indiana University Mathematics Journal~{\bf{65}} (2), 535--552 (2016).


\bibitem{DDMB} Z.A.~Daya, V.B.~Deyirmenjian, S.W.~Morris, J.R.~de Bruyn, \emph{Annular electroconvection with shear}, Phys. Rev. Lett.~{\bf{80}}, 964--967 (1998).

\bibitem{DeHiPr} R. Denk, M. Hieber, J. Prüss, \emph{Optimal $L^p-L^q$ estimates for parabolic boundary value problems with inhomogeneous data}, Math. Z.~{\bf{257}}, 193--224 (2007). 


\bibitem{DM} X-T. Duong, \emph{The $L^p$ boundedness of Riesz transforms associated with divergence form operators}, Proc. Centre Math. Appl.~{\bf{37}}, 15--25 (1999).

\bibitem{EH} D.E. Edmunds, R. Hurri-Syrjänen, \emph{Weighted Hardy inequalities}, Journal of Mathematical Analysis and Applications~{\bf{310}} (2), 424--435 (2005).

\bibitem{FaKoWe} R. Farwig, H. Kozono, D. Wegmann, \emph{Maximal regularity of the Stokes operator in an exterior domain with moving boundary and application to the Navier-Stokes equations}, Mathematische Annalen~{\bf{375}}, 949--972 (2019)
 
\bibitem{FT} C. Foias, R. Temam, \emph{Gevrey class regularity for the solutions of the Navier-Stokes equations}, Journal of Functional Analysis~{\bf{87}} (2), 359--369 (1989).

\bibitem{GeHeHi} M. Geissert, H. Heck, M. Hieber, \emph{$L^p$-theory of the Navier-Stokes flow in the exterior of a moving or rotating obstacle}, J. reine angew. Math.~{\bf{596}}, 45--62 (2006). 


\bibitem{GM} Y. Giga, T. Miyakwa, \emph{Solutions in $L^r$ of the Navier-Stokes initial value problem}, Archive for Rational Mechanics and Analysis~{\bf{89}}, 267–-281 (1985).

\bibitem{GS} J-L. Guermond, A. Salgado, \emph{A note on the Stokes operator and its powers}, Journal of Applied Mathematics and Computing volume~{\bf{36}}, 241--250 (2011).

\bibitem{HiPr} M. Hieber, J. Pruss, \emph{Heat kernels and maximal $L^p-L^q$ estimates for parabolic evolution equations}, Communications in Partial Differential Equations~{\bf{22}}, 1647--1669 (1997). 


\bibitem{HiSa} M. Hieber, J. Saal, \emph{The Stokes Equation in the $L^p$-Setting: Well-Posedness and Regularity Properties}, In: Giga, Y., Novotný, A. (eds) Handbook of Mathematical Analysis in Mechanics of Viscous Fluids,  Springer, Cham (2018).

\bibitem{HK} T. Hmidi, S. Keraani, \emph{Global solutions of the super-critical 2D quasi-geostrophic equation in Besov spaces}, Advances in Mathematics~{\bf{214}} (2), 618--638 (2007).

\bibitem{I} M. Ignatova, \emph{Construction of solutions of the critical SQG equation in bounded domains}, Advances in Mathematics~{\bf{351}}, 1000--1023 (2019).

\bibitem{JK} D. Jerison, C. Kenig, \emph{The Inhomogeneous Dirichlet Problem in Lipschitz Domains}, Journal of Functional Analysis~{\bf{130}}, 161--219 (1995). 

\bibitem{kim} J-M. Kim, \emph{On regularity criteria of the Navier-Stokes equations in bounded domains}, Journal of Mathematical Physics~{\bf{51}}, 053102 (2010). 

\bibitem{knv} A.~Kiselev, F.~Nazarov, and A.~Volberg, \emph{Global well-posedness for the critical 2{D} dissipative quasi-geostrophic equation}, Invent. Math.~{\bf{167}} (3), 445--453 (2007).

\bibitem{MLCAW} M-C. MA, G. Li, X. Chen, L-A. Archer, J. Wan, \emph{Suppression of dendrite growth by cross-flow in microfluidics}, Science Advances~{\bf{7}} (8) (2021). 

\bibitem{mani} A. Mani, K. M. Wang, Electroconvection Near Electrochemical Interfaces: Experiments, Modeling, and Computation,
Ann. Review of Fluid Mech. {\bf{52}}, 509-529 (2020).

\bibitem{Mar} P. Maremonti, \emph{On the $L^p-L^q$ estimates of the gradient of solutions to the Stokes problem}, Journal of Evolution Equations~{\bf{19}}, 645--676 (2019).

\bibitem{S} Z. Shen, \emph{Bounds of Riesz transforms on $L^p$ spaces for second order elliptic operators}, Annales de l'Institut Fourier~{\bf{55}} (1), 173--197 (2005).


\bibitem{ShSh} Y. Shibata, S. Shimizu, \emph{On the maximal $L^p-L^q$ regularity of the Stokes problem with first order boundary condition; model problems}, J. Math. Soc. Japan~{\bf{64}} (2), 561--626 (2012). 



\bibitem{SV} L.F. Stokols, A.F. Vasseur, \emph{H\"older Regularity up to the Boundary for Critical SQG on Bounded Domains}, Archive for Rational Mechanics and Analysis~{\bf{236}}, 1543--1591 (2020). 

\bibitem{TR} J. Tan, E. Ryan, \emph{Computational study of electro-convection effects on dendrite growth in batteries}, Journal of Power Sources~{\bf{323}}, 67--77 (2016).

\bibitem{Tol} P. Tolksdorf, \emph{On the $L^p$-theory of the Navier-Stokes equations on three-dimensional bounded Lipschitz domains}, Math. Ann.~{\bf{371}}, 445--460 (2018).

\bibitem{ToWa} P. Tolksdorf, K. Watanabe, \emph{The Navier–Stokes equations in exterior Lipschitz domains: $L^p$-theory}, Journal of Differential Equations~{\bf{269}} (7), 5765--5801 (2020). 


\bibitem{DDMT} P.~Tsai, Z.A.~Daya, V.B.~Deyirmenjian, S.W.~Morris, \emph{Direct numerical simulation of supercritical annular electroconvection}, Phys. Rev E~{\bf{76}}, 1--11 (2007). 

\bibitem{ZYZZHLYHLBZZ} Y.~Zhang, X.~Yang, Y.~Zhan, Y.~Zhang, J.~He, P.~Lv, D.~Yuan, X.~Hu, D.~Liu, D-J.~Broer, G.~Zhou, W.~Zhao, \emph{Electroconvection in Zwitterion-Doped Nematic Liquid Crystals and Application as Smart Windows}, Advanced Optical Materials~{\bf{9}} (3) (2020). 

\end{thebibliography}
\end{document}